\pdfoutput=1    

\documentclass[reqno, 12pt]{amsart}

\usepackage{cmap}

\usepackage{amsmath,mathtools}  
\usepackage{amsfonts}
\usepackage{amssymb}
\usepackage{graphicx}

\usepackage{comment}  

\usepackage{url}

\usepackage{ifthen}
\newboolean{BKMRK}
\setboolean{BKMRK}{true}    

\ifthenelse {\boolean{BKMRK}}
  { \usepackage{hyperref} }
  {  }

\ifthenelse {\boolean{BKMRK}}
  { \hypersetup{colorlinks=true, linkcolor=red, citecolor=red, pdfstartview=FitH, linktocpage=false
  ,pdflang={en-US}, pdftitle={Fun with Fourier series}, pdfauthor={Robert Baillie}, pdfsubject={Mathematics}, pdfkeywords={Fourier series, sinc function, polylogarithms}
               }
  }
  {  }

\def\picDblWidth{3.20in}    
\def\picSingleWidth{4.0in}

\newtheorem{theorem}{Theorem}
\theoremstyle{plain}

\numberwithin{equation}{section}

\def\figureMinipageWidth{.49\textwidth}  

\usepackage[top = .75in, bottom = .75in, left = .75in, right = .75in]{geometry}

\input{glyphtounicode}
\begin{document}

\title{Fun With Fourier Series}
\author{Robert Baillie}

\date{\today}
\subjclass[2020]{Primary 40-01; Secondary 42-01}

\keywords{Fourier series, sinc function, polylogarithms}

\begin{abstract}
By using computers to do experimental manipulations on Fourier series, we construct additional series with interesting properties.
We construct several series whose sums remain unchanged when the $n$\textsuperscript{th} term is multiplied by $\sin(n)/n$.
One example is this classic series for $\pi/4$:
\begin{align*}
\frac{\pi}{4}
  & = 1 - \frac{1}{3} + \frac{1}{5} - \frac{1}{7} + \dots \\
  & = 1 \cdot \frac{\sin(1)}{1} - \frac{1}{3} \cdot \frac{\sin(3)}{3} + \frac{1}{5} \cdot \frac{\sin(5)}{5} - \frac{1}{7} \cdot \frac{\sin(7)}{7} + \dots \, .
\end{align*}
Another example is
\[
\sum_{n=1}^{\infty} \frac{\sin(n)}{n} = \sum_{n=1}^{\infty} \left(\frac{\sin(n)}{n}\right)^2.
\]

\end{abstract}

\maketitle

\setcounter{tocdepth}{2}  
\smaller[2]    
\tableofcontents
\larger[2]

\setlength{\parskip}{10pt plus 1pt minus 1pt}  
\setlength{\parindent}{0pt}  

\section{Introduction}  


Mathematics can be \emph{fun}!

One of the fun things in mathematics is discovering new things.
It's fun to discover something you didn't know before, even if that fact is already known to others.
One way to make discoveries in mathematics is to try things out to see what happens.
For example, consider any of the many infinite series for $\pi$.
What happens if you square each term in the series?
What do you get - just some random number, or a new expression that involves $\pi$ or other constants?

Computers remove the drudgery from doing arithmetic, algebra, and calculus.
With little effort, a computer can add millions of terms of an infinite series.
Furthermore, modern computer algebra systems make it possible to graph, simplify, and manipulate very complicated algebraic expressions.
This allows us to answer exploratory ``What if?'' questions quickly and with little chance of human error.
The practice of using computers to explore mathematics and try new things even has a name: \emph{Experimental Mathematics}.

In this paper, we will use computers to experiment with Fourier series.
This leads to some surprising variations on standard formulas for $\pi$.

The reader may wish to glance at Section \ref{S:CollectionOfFormulas} to see some of the formulas we will discover and prove.

The author has always thought that Fourier series were amazing.
For example, Fourier series provide an easy way to prove surprising facts like
\begin{equation} \label{E:PiSqOver6}
  \frac{\pi^2}{6} = \frac{1}{1^2} + \frac{1}{2^2} + \frac{1}{3^2} + \frac{1}{4^2} + \frac{1}{5^2} + \dots \, .
\end{equation}
($\pi$ is related to circles.
What does $\pi$ have to do with these squares of integers?
Videos \cite{Mathologer} and \cite{3BlueOneBrown} give two different explanations of the connection!)

Here is one proof of Equation \eqref{E:PiSqOver6}: over the interval $0 \leq x \leq 2 \pi$, the polynomial $(3 x^2 - 6 \pi x + 2 \pi^2)/12$ has the Fourier series
\begin{equation} \label{E:Zeta2}
\frac{3 x^2 - 6 \pi x + 2 \pi^2}{12} =
  \sum_{n=1}^{\infty} \frac{\cos(n x)}{n^2} \, .
\end{equation}
This is Equation 13.8 on page 31 of Tolstov's book \cite{Tolstov}.
We will prove this later (Equation \eqref{E:FS-SpecialQuadratic}).

If we substitute $x = 0$ into Equation \eqref{E:Zeta2}, the left side reduces to $\pi^2/6$.
Also, every numerator $\cos(n x)$ is $\cos(0) = 1$, which gives us the series on the right side of \eqref{E:PiSqOver6}.

\begin{figure}[ht]
  \mbox{\includegraphics[width=\picSingleWidth]{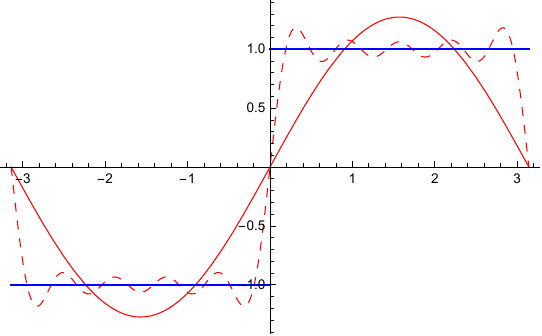}}
  \caption{Square wave (blue), first term (solid red), and sum of 5 terms (dashed)}
  \label{fig:squareWave}
\end{figure}
Not only that, but it is also surprising that a ``square wave'', which is horizontal for the \emph{uncountable} number of $x$ values in an entire interval,
can be produced by adding a \emph{countable} number of sine functions, each of which is curved.
Here is the Fourier series for a square wave with period $2 \pi$:
\begin{equation}
\frac{4}{\pi} \left( \sin(x) +  \frac{\sin(3x)}{3} + \frac{\sin(5x)}{5} \ldots \right) =
\frac{4}{\pi} \sum_{n=1}^{\infty} \frac{\sin((2n-1)x)}{2n-1} =
  \begin{cases}
   -1  &\text{if $-\pi < x < 0$,} \\
    0  &\text{if $x = 0$,} \label{E:unitSquareWave} \\
    1  &\text{if $0 < x < \pi$.}
  \end{cases}
\end{equation}
This is Example 1 on p.\ 497 of Kaplan's book, \cite{Kaplan}.
Figure \ref{fig:squareWave} shows this square wave in blue, together with the first term of its series (solid red) and the sum of the first 5 terms (dashed).

We will find a number of infinite series that have the amusing property that the $n$\textsuperscript{th} term can be multiplied by $\sin(n)/n$ without changing the sum.
One series with this property is the classic Gregory-Leibniz series \cite{BGL}.
\begin{equation*}
\frac{\pi}{4} = 1 - \frac{1}{3} + \frac{1}{5} - \frac{1}{7} \dots \, .
\end{equation*}
We will prove that
\begin{equation*}  
\frac{\pi}{4} = 1 - \frac{1}{3} + \frac{1}{5} - \frac{1}{7} \dots \,
= 1 \cdot \frac{\sin(1)}{1} - \frac{1}{3} \cdot \frac{\sin(3)}{3} + \frac{1}{5}  \cdot \frac{\sin(5)}{5} - \frac{1}{7}  \cdot \frac{\sin(7)}{7} \dots \, .
\end{equation*}
Another series with this property is
\begin{equation}\label{E:sincn}
\sum_{n=1}^{\infty} \frac{\sin(n)}{n} = \sum_{n=1}^{\infty} \left(\frac{\sin(n)}{n}\right)^2 \, .
\end{equation}
Every term on the right is different from the corresponding term on the left,
so the reader may wonder how Equation \eqref{E:sincn} is even possible.

\begin{figure}[ht]
  \mbox{\includegraphics[width=\picSingleWidth]{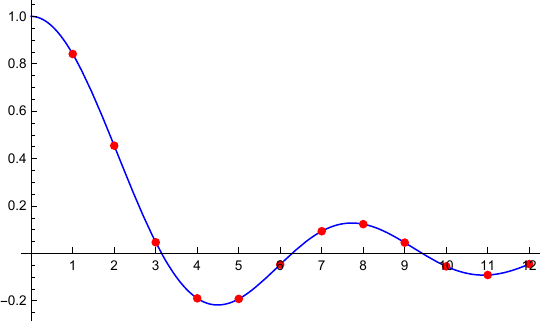}}
  \caption{$\sin(x)/x$ for $0 < x \leq 12$ with a dot at each $\sin(n)/n$}
  \label{fig:sincPlot}
\end{figure}

The first few values of $\sin(n)/n$ are the dots on the graph in Figure \ref{fig:sincPlot}.
Some numeric values of $\sin(n)/n$ and $(\sin(n)/n)^2$ are shown in Table \ref{Ta:SincTable}.

Some values of $\sin(n)/n$ are positive and some are negative, so, there is some cancellation when we add them.
Also, every $\sin(n)/n$ is between $-1$ and $+1$, so $(\sin(n)/n)^2$ is positive, but closer to zero (i.e., tinier) than $\sin(n)/n$.
So, Equation \eqref{E:sincn} works because the cancellation in the sum on the left of is exactly balanced out by the fact that the terms on the right are all positive, but tinier.


\begin{table}[ht]
 \begin{center}
  \begin{tabular}{ r r r r }
  \rule[-8pt]{0pt}{22pt}  
   $n$  &  $\sin(n)$   &   $\frac{\sin(n)}{n}$   &   $\left( \frac{\sin(n)}{n} \right) ^2$  \\ \hline

   1  & $  0.84147 $ & $   0.84147 $ & $ 0.70807 $ \\
   2  & $  0.90930 $ & $   0.45465 $ & $ 0.20671 $ \\
   3  & $  0.14112 $ & $   0.04704 $ & $ 0.00221 $ \\
   4  & $ -0.75680 $ & $  -0.18920 $ & $ 0.03580 $ \\
   5  & $ -0.95892 $ & $  -0.19178 $ & $ 0.03678 $ \\
   6  & $ -0.27942 $ & $  -0.04657 $ & $ 0.00217 $ \\
   7  & $  0.65699 $ & $   0.09386 $ & $ 0.00881 $ \\
   8  & $  0.98936 $ & $   0.12367 $ & $ 0.01529 $ \\
   9  & $  0.41212 $ & $   0.04579 $ & $ 0.00210 $ \\
  10  & $ -0.54402 $ & $  -0.05440 $ & $ 0.00296 $ \\
  11  & $ -0.99999 $ & $  -0.09091 $ & $ 0.00826 $ \\
  12  & $ -0.53657 $ & $  -0.04471 $ & $ 0.00200 $ \\

  \end{tabular}
  \caption{$\sin(n)$, $\sin(n)/n$, and $(\sin(n)/n)^2$}
  \label{Ta:SincTable}
 \end{center}
\end{table}

We even find a series where we can multiply the $n$\textsuperscript{th} term by $\sin(n)/n$, $(\sin(n)/n)^2$, and $(\sin(n)/n)^3$, and all \emph{four} series have the same sum!
See Equation \eqref{E:FourEqualSums}.

Sometimes we will need to work backwards from the Fourier series of an unknown function, to express that function in terms of polynomials.
This involves trial and error and detective work and sometimes fails, as we shall see in Section \ref{S:Recovering}.
Polylogarithms are also useful here; see Section \ref{S:UsingPolylogs}.

Most of the results in this paper were first found by doing computer experiments.
The author used the \textit{Mathematica} computer algebra system to try out variations of known series just to see what
would happen, to plot complicated functions, and to do the integrations needed to calculate coefficients of Fourier series.
Often, \textit{Mathematica} output would make it clear that a conjectured equality was false.
This means that we don't need to waste time trying to prove something that isn't true.

Unfortunately, most published papers in mathematics (including some written by this author) fail to explain how results were discovered.
Instead, they merely give theorems and proofs, as if the theorems were found by magic.

So, our goal here is not just to list and prove theorems.
Instead, the goal is to show \emph{how} computers can help us discover and prove some interesting results.
The results we obtain illustrate some of the ``low-hanging fruit'' that can be discovered without using advanced mathematics.
It is hoped that the reader will find these results interesting, and will be inspired to do computer experiments of their own.

\textbf{If you don't have \textit{Mathematica}.}
Even if the reader does not have access to \textit{Mathematica}, it is still possible to use \textit{Mathematica} to verify at least \emph{some} of the results.
For example, consider Equation \eqref{E:TwoPieceIntegral}.
The \emph{second} integral in that equation is
\[
\int_{1}^{2\pi - 1} \frac{\pi - x}{2} \sin(nx) \, dx \, .
\]
Here is the \textit{Mathematica} code to do this calculation and simplify the result:
\begin{verbatim}
  i2tmp = Integrate[Sin[n x] * (Pi - x)/2, {x, 1, 2 Pi - 1}]
  i2 = Simplify[i2tmp, Assumptions -> Element[n, Integers]]
\end{verbatim}
The \textit{Wolfram Alpha} website \url{https://www.wolframalpha.com/} can compute the integral.
Just copy
\begin{verbatim}
  Integrate[Sin[n x] * (Pi - x)/2, {x, 1, 2 Pi - 1}]
\end{verbatim}
and paste it into that webpage.
The result of the integration will be displayed:
\[
\frac{\cos (\pi  n) (\sin (n-\pi  n)+(\pi -1) n \cos (n-\pi  n))}{n^2} \, .
\]
However, \textit{Wolfram Alpha} is unable within the alloted time to do the simplification specified in the second line of code (\verb+i2 = ...+), which would give
\[
\frac{(\pi -1) n \cos (n) + \sin (n)}{n^2} \, .
\]
So, without access to \textit{Mathematica}, the reader will have to do more of the work by hand.



\subsection{What's new in the April, 2026 version of this paper.} \label{S:WhatsNew}

Section \ref{S:AnotherTrilogExample} was added to provide a simple proof that
\begin{equation*}
\sum_{n = 1}^{\infty} (-1)^{n+1} \left( \frac{\sin(n)}{n} \right)^3 = \frac{1}{2} \, .
\end{equation*}
Equation \eqref{E:3sumsOneHalf} was added to emphasize a rather amusing identity.

\textbf{What's new in the June, 2023 version of this paper.}

Above, we added links to videos that explain Equation \eqref{E:PiSqOver6}.
At the end of the proof of Theorem \ref{T:Thm1}, we added a comment that polylogarithms can be used to directly derive the value of the sums in the Theorem.
Page 56 had two typos: ``Diuviding by $n$ and summing over $n$'' was replaced with ``Dividing by $n^2$ and summing over $n$''.
In file \verb+FS.m+, the comment in line 3 was changed to give the correct link.

\textbf{What's new in the May, 2023 version of this paper.}

Many parts of this paper were rewritten for (hopefully) increased clarity.

In Equation \eqref {E:SumAndSumCubed}, we derive an example of a series where $\sum a_{n} = \sum a_{n}^{3}$.
This equation has also been added to the collection of interesting formulas in the Appendix, Section \ref{S:CollectionOfFormulas}.

Section \ref{S:Package} describes a \textit{Mathematica} package, \verb+FS.m+, that makes it easy to compute Fourier series for many types of functions.
That section was rewritten to include more examples.
You can enter a function like $x^2 - x$, and \textit{Mathematica} will compute expressions for the coefficients of $\cos(nx)$ and $\sin(nx)$ in the Fourier series for your function.
\textit{Mathematica} also draws graphs of your function and the sum of terms in the Fourier series.
The package also allows you to enter piecewise continuous functions like the square wave in Figure \ref{fig:squareWave}, or like $g(x)$ in Equation \eqref{E:gxTentative2}.
You can also enter piecewise linear functions by specifying just the ($x$, $y$) coordinates at endpoints, corners, or jumps.

Several functions to evaluate Parseval's formula (Equation \eqref{E:ParsevalsEquation}) were added to \verb+FS.m+.

A bug in \verb+sinSeries[ ]+ and \verb+cosSeries[ ]+ was fixed so these functions now work in recent versions of \textit{Mathematica}.

The text file \verb+FS.m+ is an ``ancillary file'' at the arXiv link to this paper at \\
\url{https://arxiv.org/abs/0806.0150} . \\
The direct link to this version of \verb+FS.m+ is \\
\url{https://arxiv.org/src/0806.0150/anc/FS.m} .

The small package that had been included in Section \ref{S:Package} was merely a subset of \verb+FS.m+, so it has been removed.

\section{Preliminaries and prerequisites} \label{S:Preliminaries}

The most advanced prerequisite needed here is a basic knowledge of Fourier series.
Chapter 7 of Kaplan's \textit{Advanced Calculus} \cite{Kaplan} is sufficient.
Tolstov's book \cite{Tolstov} is another good source of information and examples.

We will use the following facts from analytic geometry: suppose we have two points $(x_1, y_1)$ and $(x_2, y_2)$, with $x_1 \neq x_2$.
Then the equation of the line between these points is given by $y = m x + b$, where
\begin{equation}\label{E:slope}
m = \frac{y_2 - y_1}{x_2 - x_1}
\end{equation}
and
\begin{equation}\label{E:intercept}
b = y_1 - m x_1 \, .
\end{equation}
$m$ is the slope of the line. $b$ is the $y$-intercept, that is, the $y$ value where the line intersects the $y$ axis.

We will be working with functions that are ``piecewise very smooth'' \cite[page 496]{Kaplan}.
This means they are piecewise continuous, have at most a finite number of jump discontinuities over a relevant interval (say, $[-\pi, \pi]$),
and the ``pieces'' have continuous first and second derivatives.

Such functions $F(x)$ have the Fourier series representation
\begin{equation} \label{E:FxRepPeriod2Pi}
F(x) \approx \frac{a_0}{2} + \sum_{n=1}^{\infty} a_n \cos(n x) + b_n \sin(n x) \, .
\end{equation}

We used ``$\approx$'' instead of ``$=$'' for the following reason: for those $x$ where $F(x)$ is continuous, both sides are equal.
But for values of $x$ at which $F(x)$ has a jump discontinuity from $y = y_0$ to $y = y_1$,
the sum on the right converges to the midpoint of the $y$ values, that is $(y_0 + y_1)/2$.

The Fourier coefficients $a_n$, for $n \geq 0$, and $b_n$, for $n \geq 1$, can be computed as follows.
If $f(x)$ is periodic, and one period is from $x = -\pi$ to $x = \pi$, then we say that $f(x)$ is periodic from $x = -\pi$ to $x = \pi$.
In this case, $a_n$ and $b_n$ are given by
\begin{equation}\label{E:an}
  a_n = \frac{1}{\pi} \int_{-\pi}^{\pi} F(x) \cos(nx) \, dx
\end{equation}
and
\begin{equation}\label{E:bn}
  b_n = \frac{1}{\pi} \int_{-\pi}^{\pi} F(x) \sin(nx) \, dx \, .
\end{equation}

If the function is periodic from $0$ to $2\pi$ instead of from $-\pi$ to $\pi$, then in Equations \eqref{E:an} and \eqref{E:bn}, we integrate from $0$ to $2\pi$ instead of from $-\pi$ to $\pi$.

Sometimes, we will work with functions that have period different from $2\pi$,
or are defined over intervals other than $(-\pi, \pi)$ or $(0, 2 \pi)$.
If the interval is from $x = a$ to $x = b$ and the period is $p = b - a$, then equations \eqref{E:an} and \eqref{E:bn} become
\begin{equation}\label{E:anPeriodPab}
  a_n = \frac{2}{p} \int_{a}^{b} F(x) \cos \left( \frac{2 \pi}{p} n x \right) \, dx
\end{equation}
and
\begin{equation}\label{E:bnPeriodPab}
  b_n = \frac{2}{p} \int_{a}^{b} F(x) \sin \left( \frac{2 \pi}{p} n x \right) \, dx \, .
\end{equation}

If the period is $p$, then the function has the Fourier series
\begin{equation} \label{E:FxRepPeriodP}
F(x) \approx \frac{a_0}{2} + \sum_{n=1}^{\infty} a_n \cos \left( \frac{2 \pi}{p} n x \right) + b_n \sin \left( \frac{2 \pi}{p} n x \right) \, ,
\end{equation}
which is equivalent to the simpler \eqref{E:FxRepPeriod2Pi} if $p = 2 \pi$.

\textbf{Sine and cosine series:}

In some of our examples, instead of functions periodic over $(-\pi, \pi)$, we will work with \emph{even} functions over $(0, \pi)$.
An \emph{even} function is a function such that $F(-x) = F(x)$.
For example, $x^2$ and $\cos(x)$ are even.

By extending the function to negative $x$ values with $F(-x) = F(x)$, the resulting function over $(-\pi, \pi)$ has period $2\pi$.
In this case, the $b_n$ are all 0, and only the $a_n$ need to be computed.
The resulting series is called the ``Fourier cosine series''.
For $n \geq 0$, the $a_n$ can be computed from
\begin{equation}\label{E:anEven}
  a_n = \frac{2}{\pi} \int_{0}^{\pi} F(x) \cos(nx) \, dx \, .
\end{equation}
Comparing with Equation \eqref{E:an}, note that we integrate over half the period, then multiply by 2.

An \emph{odd} function is a function such that $F(-x) = -F(x)$.
For example, $x$ and $\sin(x)$ are odd.
In this case, the $a_n$ are all 0, and only the $b_n$ need to be computed.
The series is then called the ``Fourier sine series''.
If the function is odd over the interval $(-\pi, \pi)$, then
for $n \geq 1$, the $b_n$ can be computed from
\begin{equation}\label{E:bnOdd}
  b_n = \frac{2}{\pi} \int_{0}^{\pi} F(x) \sin(nx) \, dx \, .
\end{equation}

Suppose $F(x)$ is \emph{even} over the interval $0 < x < b$ but the period $p$ is $2b$.
Then for $-b < x < 0$, we have $-x > 0$, and $F(x) = F(-x)$.
For $n \geq 0$, the Fourier cosine coefficients are given by this generalization of Equation \eqref{E:anEven}
\begin{equation}\label{E:anEvenPeriodP}
  a_n = \frac{2}{b} \int_{0}^{b} F(x) \cos \left( \frac{2 \pi}{p} n x \right) \, dx \, .
\end{equation}

Suppose $F(x)$ is \emph{odd} over $0 < x < b$, and $p = 2b$.
For $n \geq 1$, the Fourier sine coefficients are given by
\begin{equation}\label{E:bnOddPeriodP}
  b_n = \frac{2}{b} \int_{0}^{b} F(x) \sin \left( \frac{2 \pi}{p} n x \right) \, dx \, .
\end{equation}

\textbf{Integrating polynomials with trig functions.}

Most of the Fourier coefficients we'll encounter can be computed by integrating products of polynomials
and trig functions.
For this, integration by parts is the standard technique to use.

Here, $u$ and $v$ are functions of $x$:
\begin{equation} \label{E:ByParts}
\int u \, dv \, = u v - \int v \, du \, .
\end{equation}

A typical example might involve computing an integral like
\[
\int x^2 \sin(n x)  \, dx \, .
\]
To evaluate this, we would take $u = x^2$ and $dv = \sin(n x) \, dx$.
Then $du = 2 x \, dx$ and $v = (-1/n) \cos(n x)$, so
\[
\int x^2 \sin(n x)  \, dx = -\frac{x^2}{n} \cos(n x) + \frac{2}{n} \int x \cos(n x)  \, dx \, .
\]
This reduces the exponent of $x$ by 1.
We now repeat the process with $u = x$ and $dv = \cos(n x) \, dx$:
Then $du = dx$ and $v = (1/n) \sin(n x)$, so
\[
\int x \cos(n x)  \, dx = \frac{x}{n} \sin(n x) - \frac{1}{n} \int \sin(n x)  \, dx
= \frac{x}{n} \sin(n x) + \frac{\cos(n x)}{n^2} \, .
\]
Therefore,
\[
\int x^2 \sin(n x)  \, dx \ =
-\frac{x^2}{n} \cos(n x) + \frac{2}{n} \left( \frac{x}{n} \sin(n x) + \frac{\cos(n x)}{n^2} \right) \, .
\]

\textbf{An easy example.}

It is a standard textbook exercise to calculate that
\begin{equation}
\sum_{n = 1}^{\infty} \frac{\sin(n x)}{n} = \frac{\pi - x}{2}  \quad \text{ for $0 < x < 2 \pi$}  \label{E:five1} \, .
\end{equation}
For example, see \cite[Equation (13.7) on page 31]{Tolstov} or \cite[Equation (c) on page 506]{Kaplan}.

We can verify Equation \eqref{E:five1} as follows.

\begin{figure}[ht]
  \mbox{\includegraphics[width=\picSingleWidth]{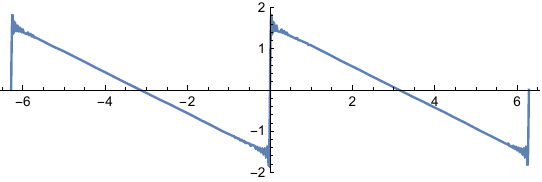}}
  \caption{Sum of the first 100 terms of the series in Equation \eqref{E:five1}}  
  \label{fig:f1Plot2Periods}
\end{figure}
The graph in Figure \ref{fig:f1Plot2Periods} shows \emph{two} periods of the sum in \eqref{E:five1}.

Consider the function $f(x) = (\pi - x)/2$, where one period is $0 < x < 2 \pi$.
Just for purposes of illustration, we will calculate its Fourier series in two ways.
First, we use Equations \eqref{E:anPeriodPab} and \eqref{E:bnPeriodPab}, with $a = 0$ and $b = 2 \pi$, so the period $p = 2 \pi$.
We integrate by parts to get:
\[
a_n = \frac{1}{\pi} \int_{0}^{2 \pi} \frac{\pi - x}{2} \cos(n x) \, dx = 0 \, ,
\]
and
\[
b_n = \frac{1}{\pi} \int_{0}^{2 \pi} \frac{\pi - x}{2} \sin(n x) \, dx = \frac{1}{n} \, .
\]
This verifies that $f(x)$ has the fourier series given in Equation \eqref{E:five1}.

Since Equation \eqref{E:five1} has only sine functions, a second approach is to consider $f(x)$ to be an \emph{odd} function over $-\pi < x < \pi$, whose period is still $2 \pi$.
In this case, we can use Equation \eqref{E:bnOdd} to compute only the sine coefficients, because the cosine coefficients of an odd function are 0.
Here, we integrate from $x = 0$ to $x = \pi$, then multiply by 2.
We get the same result as before:
\[
b_n = \frac{2}{\pi} \int_{0}^{\pi} \frac{\pi - x}{2} \sin(n x) \, dx = \frac{1}{n} \, .
\]
This second approach is simpler because we need to evaluate only one integral.

We should note that, for $-\pi < x < 0$, f(x) is \emph{not} $(\pi - x)/2$:
$f(x)$ is odd, so for $-\pi < x < 0$, we have $f(x) = -f(-x) = -(\pi - (-x))/2 = -(\pi + x)/2$.

\textbf{Parseval's equation.}

If
\[
F(x) \approx \frac{a_0}{2} + \sum_{n=1}^{\infty} a_n \cos(n x) + b_n \sin(n x)
\]
is the Fourier series for $F(x)$ over $[-\pi, \pi]$, then \cite[page 439]{Danese}, \cite[page 519]{Kaplan}, \cite[page 119]{Tolstov}
\begin{equation}\label{E:ParsevalsEquation}
\frac{1}{\pi} \int_{-\pi}^{\pi} F(x)^2 \, dx = \frac{a_0^2}{2} + \sum_{n=1}^{\infty} (a_n^2 + b_n^2) \, .
\end{equation}

More generally, if $F(x)$ has a Fourier series over the interval from $a$ to $b$, and, as above,
its period is $p = b - a$, then the Fourier series is given by \eqref{E:FxRepPeriodP}.
Then Parseval's equation \eqref{E:ParsevalsEquation} becomes
\begin{equation}\label{E:GenParsevalEqn}
\frac{2}{p} \int_{a}^{b} F(x)^2 \, dx = \frac{a_0^2}{2} + \sum_{n=1}^{\infty} (a_n^2 + b_n^2) \, .
\end{equation}

Parseval's equation often lets us calculate sums of interesting infinite series.
Sometimes, the sum of the series on the right side of Equation \eqref{E:GenParsevalEqn}, by itself, is difficult to compute,
while the integral on the left side is easy to evaluate.
In cases like this, Parseval's equation gives us the answer.
See, for example, Equation \eqref{E:Prob6241b}.

\textbf{The sinc function.}

The \emph{sinc} function is defined as sinc$(0) = 1$, and sinc$(x) = \sin(x)/x$ if $x \neq 0$.
This makes sinc a continuous function, because $\sin(x)/x$ approaches 1 as $x$ approaches 0.
See Figure \ref{fig:sincPlot} on page \pageref{fig:sincPlot} for a plot of $\text{sinc}(x)$ for $x \geq 0$.
Many of the equations in this paper could be written a little more concisely using $\text{sinc}(x)$ instead of $\sin(x)/x$.
However, the author made the arbitrary choice to use $\sin(x)/x$ instead.
Some properties of the sinc function are listed in \cite{MathWorldSincFunction}.

\section{Weird non-cancellation in a conditionally-convergent series} \label{S:Cancellation}  

The equations in this section aren't needed for the rest of the paper.
However, they present an amusing example of what can happen with a seemingly simple operation
like cancelling $+1/n$ and $-2/(2n)$ in a series that is conditionally-convergent.


We will prove that
\begin{equation} \label{E:HalfLn2}
\frac{\ln(2)}{2} = \frac{1}{1} - \frac{2}{2} + \frac{1}{3} + \frac{1}{5} - \frac{2}{6} + \frac{1}{7} + \frac{1}{9} - \frac{2}{10}
 + \frac{1}{11} + \frac{1}{13} - \frac{2}{14} + \frac{1}{15} + \frac{1}{17} - \frac{2}{18} + \ldots \, .
\end{equation}

Notice that, for every term of the form $1/n$, there is a corresponding term $-2/(2n)$ later on.
Thus, \emph{everything} on the right \emph{appears} to cancel out, producing a sum equal to 0.

What's \emph{really} going on is that the terms $+1/n$ and $-2/(2n)$ become arbitrarily far apart,
so the cancellation may produce an incorrect result (and in this case, it \emph{does} produce an incorrect result).
If the terms did \emph{not} become arbitrarily far apart, then Theorem 12-14 of \cite[p. 357]{Apostol} would apply;
this would allow us to legitimately combine each ($+1/n$ and $-2/(2n)$) pair.

Proof of Equation \eqref{E:HalfLn2}:

Begin with the Fourier series expansion (see Equation 14.1 on page 93 of \cite{Tolstov})
\begin{equation}
-\ln \left| 2 \sin{ \left( \frac{x}{2} \right) } \right|
 = \sum_{n=1}^{\infty} \frac{\cos{nx}}{n} \, ,
\end{equation}
which is valid if $x$ is not a multiple of $2 \pi$.
It can also be obtained by the standard methods for calculating Fourier series.

Replace $x$ with $\pi - 2x$ throughout this equation.
We get
\begin{equation*}
-\ln \left| 2 \sin{ \left( \frac{\pi - 2x}{2} \right) } \right|
 = \sum_{n=1}^{\infty} \frac{\cos{n(\pi - 2x)}}{n} \, .
\end{equation*}
This replacement is valid if $\pi - 2x$ is not a multiple of $2\pi$.
Simplify this using standard identities, including $\sin(\pi/2 - x) = \cos(x)$ and $\cos(x - y) = \cos(x) \cos(y) + \sin(x) \sin(y)$, we get
\begin{align*}
\cos(n\pi - 2nx) & = \cos(n\pi)\cos(2nx) + \sin(n\pi) \sin(2nx) \\ & = (-1)^n \cos(2nx) + 0 \cdot \sin(2nx) \\ & = (-1)^n (1 - 2 \sin^2(nx)) \, .
\end{align*}

We then have
\begin{equation} \label{E:HalfLn2B}
-\ln \left| 2 \cos(x) \right|
 = \sum_{n=1}^{\infty} \frac{(-1)^n(1 - 2\sin^2{nx})}{n}
 = \sum_{n=1}^{\infty} \frac{(-1)^n}{n} + 2 \sum_{n=1}^{\infty} (-1)^{n+1} \frac{\sin^2{nx}}{n} \, .
\end{equation}

Substitute $x = \pi/4$ into equation \eqref{E:HalfLn2B}.
$2 \cos(\pi/4) = \sqrt{2}$, so the left side is $-\ln(2)/2$.
The first sum on the right side is $-\ln(2)$.
Adding $\ln(2)$ to both sides, we get
\begin{equation} \label{E:HalfLn2C}
\frac{\ln(2)}{2} = 2 \sum_{n=1}^{\infty} (-1)^{n+1} \frac{\sin^2{nx}}{n} \, .
\end{equation}
For $n$ = 1, 2, 3, and 4, the values of $\sin^2(n \pi/4)$ are $1/2$, 1, $1/2$, and 0, and this cycle repeats for $n > 4$.
The first 12 fractions in the sum on the right side of \eqref{E:HalfLn2C} are
\[  
\frac{1}{2} \; , \; -\frac{1}{2} \; , \; \phantom{-}\frac{1}{6} \; , \; -\frac{0}{4} \; , \; \phantom{-}\frac{1}{10} \; , \; -\frac{1}{6} \; , \; \phantom{-}\frac{1}{14} \; , \; -\frac{0}{8} \; , \; \phantom{-}\frac{1}{18} \; , \; -\frac{1}{10} \; , \; \phantom{-}\frac{1}{22} \; , \; -\frac{0}{12} \; , \; \ldots \, .
\]
Multiplying by the 2 in front of the sum, the right side of \eqref{E:HalfLn2C} is
\[
\frac{1}{1} - \frac{2}{2} + \frac{1}{3} - \frac{0}{4} + \frac{1}{5} - \frac{2}{6} + \frac{1}{7} - \frac{0}{8} + \frac{1}{9} - \frac{2}{10} + \frac{1}{11} - \frac{0}{12} \ldots \, ,
\]
which proves Equation \eqref{E:HalfLn2}.

For many of the series we'll be working with, the absolute value of the $n$\textsuperscript{th} term will be at most $1/n^2$, so the series are absolutely convergent.
This means we can rearrange the order of the terms of the series any way we wish, and the sum will remain unchanged.

\ifthenelse {\boolean{BKMRK}}
  { \section{Some interesting series and formulas for \texorpdfstring{$\pi$}{pi}}\label{S:InterestingPiFormulas} }
  { \section{Some interesting series and formulas for pi}\label{S:InterestingPiFormulas} }

While doing calculations on a computer in the 1970's, the author noticed that
\[
\sum_{n=1}^{1000000} \left(\frac{\sin(n)}{n} \right)^{2} \approx \frac{\pi - 1}{2} \, .
\]
(This sum is about $(\pi - 1)/2 - 0.0000005$.)
This result was surprising, because the author knew from textbooks that
\[
\sum_{n=1}^{\infty} \frac{\sin(n)}{n} = \frac{\pi - 1}{2} \, .
\]
The above calculation eventually led to our first Theorem.

\begin{theorem} \label{T:Thm1}
\begin{equation} \label{E:Prob6241a}
\sum_{n=1}^{\infty} \frac{\sin(n)}{n} =
\sum_{n=1}^{\infty} \left(\frac{\sin(n)}{n} \right) ^{2} = \frac{\pi - 1}{2} \, .
\end{equation}
\end{theorem}

\textbf{Proof of Theorem \ref{T:Thm1}.}
First, here is what one might call a ``trick'':
Observe that the leftmost sum in Equation \eqref{E:Prob6241a} is just the Fourier series
\[
\sum_{n = 1}^{\infty} \frac{\sin(n x)}{n}
\]
evaluated at $x = 1$.
This sum is $(\pi - 1)/2$ because from Equation \eqref{E:five1}, we already know that
\[
\sum_{n = 1}^{\infty} \frac{\sin(n x)}{n} = \frac{\pi - x}{2}  \quad \text{ for $0 < x < \pi$} \, .
\]

\begin{figure}[ht] 
\centering
\begin{minipage}[t]{\figureMinipageWidth}
\centering
\includegraphics[width=\picDblWidth]{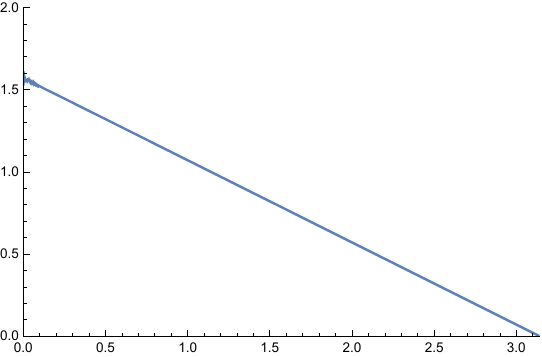}
\caption{$\frac{\pi - x}{2} = \sum_{n = 1}^{\infty} \frac{\sin(nx)}{n}$}
\label{fig:f1Plot}
\end{minipage}\hfill
\begin{minipage}[t]{\figureMinipageWidth}
\centering
\includegraphics[width=\picDblWidth]{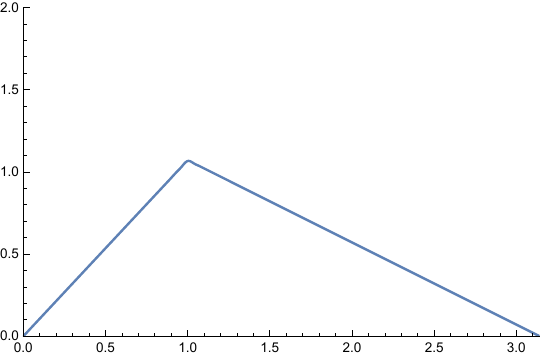}
\caption{$\sum_{n = 1}^{100} \frac{\sin(n)}{n^2} \cdot \sin(n x)$}
\label{fig:gxPlot}
\end{minipage}
\end{figure}

Figure \ref{fig:f1Plot} shows this function for $0 < x < \pi$.
We will now prove that
\begin{equation} \label{E:SumSincSquared}
\sum_{n=1}^{\infty} \left(\frac{\sin(n)}{n} \right) ^{2} = \frac{\pi - 1}{2} \, .
\end{equation}

We will use a ``trick'' similar to the one we used above:
We can express the sum
\[
\sum_{n=1}^{\infty} \left( \frac{\sin(n)}{n} \right) ^{2}
\]
as the Fourier series
\[
\sum_{n=1}^{\infty} \frac{\sin(n)}{n^2} \cdot \sin(n x)
\]
evaluated at $x = 1$.

So, let's consider the function
\begin{equation}\label{E:gxTentative}
g(x) = \sum_{n=1}^{\infty} \frac{\sin(n)}{n^2} \cdot \sin(nx).
\end{equation}
What does $g(x)$ look like?
The sum of the first 100 terms of this series is shown in Figure \ref{fig:gxPlot}.
Compare this graph to Figure \ref{fig:f1Plot}:
Notice that $g(x)$ in Figure \ref{fig:gxPlot} appears to match $f(x) = (\pi - x)/2$
over an interval from $x = 1$, to at least $x = \pi$.

For $0 \leq x \leq \pi$, $g(x)$ \emph{appears} to consist of two linear segments.
So, let's tentatively assume that $g(x)$ \emph{does} consist of two linear segments.
We will obtain expressions for these segments, then we'll compute the Fourier series
for the resulting function and see if it matches Equation \eqref{E:gxTentative}.

To get the equation(s) for $g(x)$, we will use Equations \eqref{E:slope} and \eqref{E:intercept}.

Segment 1:
From Equation \eqref{E:gxTentative}, $g(0) = 0$.
If we set $x = 1$ and add the first million terms, we get $g(1) \approx 1.0707958$, which is very close to $f(1) = (\pi - 1)/2 \approx 1.0707963$.
If $g(1)$ is, indeed, $(\pi - 1)/2$, then the first linear part of $g(x)$
goes from the point $(0, 0)$ to the point $(1, (\pi - 1)/2)$.
From Equations \eqref{E:slope} and \eqref{E:intercept}, $m = (\pi - 1)/2$ and $b = 0$.
So, the expression for the first segment of $g(x)$ would be $x(\pi - 1)/2$.

Segment 2:
Here, $g(x)$ starts at $g(1) = (\pi - 1)/2$ and, from Equation \eqref{E:gxTentative}, $g(\pi) = 0$.
Therefore, segment 2 of $g(x)$ goes from $(1, (\pi - 1)/2)$ to $(\pi, 0)$.
From Equations \eqref{E:slope} and \eqref{E:intercept}, $m = -1/2$ and $b = \pi/2$.
So, for $1 \leq x \leq \pi$, the expression for the second segment of $g(x)$ would be $(\pi - x)/2$.

Putting all this together, it \emph{appears} that for $0 \leq x \leq \pi$,
\begin{equation} \label{E:gxTentative2}
g(x) =
  \begin{cases}  
    x(\pi - 1)/2  &\text{for $0 \leq x < 1$,}\\
    (\pi - x)/2  &\text{for $1 \leq x \leq \pi$.}
  \end{cases}
\end{equation}
We hope that the Fourier series for this function is given by Equation \eqref{E:gxTentative}.
We will now determine whether this is the case.

We will consider $g(x)$ to be an \emph{odd} function with period $2 \pi$, so in the Fourier series, all cosine terms will be 0:
the values from 0 to $-\pi$ will be the negatives of the values from 0 to $\pi$.

The coefficient of $\sin(nx)$ in the Fourier series for this $g(x)$ is computed from Equation \eqref{E:bnOdd}, with $b = \pi$:
\begin{align}\label{E:TwoPieceIntegral}
b_n
 = & \frac{2}{\pi} \int_{0}^{\pi} g(x) \sin(nx) \, dx \\
 = & \frac{2}{\pi} \left( \int_{0}^{1} \frac{x(\pi - 1)}{2} \sin(nx) \, dx +
  \int_{1}^{\pi} \frac{\pi - x}{2} \sin(nx) \, dx 
  \right)  \notag \\
 = & \frac{2}{\pi} \left( -\frac{(\pi - 1) (n \cos(n)-\sin(n))}{2 n^2} +
    \frac{(\pi - 1) n \cos(n) + \sin(n)}{2 n^2}
      \right) \notag \\
=  & \frac{2}{\pi} \left(
  -\frac{(\pi - 1) n \cos(n) - (\pi - 1) \sin(n)}{2 n^2} +
  \frac{(\pi - 1) n \cos(n) + \sin(n)}{2n^2}
   \right) \notag \\
 = & \frac{2}{\pi} \left(  \frac{\pi \sin(n)}{2 n^2}  \right) \notag \\
 = & \frac{\sin(n)}{n^{2}} \notag \, .
\end{align}
Bingo!

So, the $g(x)$ given by Equation \eqref{E:gxTentative2} \emph{does} have the Fourier series given in Equation \eqref{E:gxTentative}.
In particular, we have proved that, for $1 \leq x \leq \pi$,
\begin{equation}
\sum_{n=1}^{\infty} \frac{\sin(nx)}{n} = \sum_{n=1}^{\infty} \frac{\sin(n)}{n^2} \cdot \sin(nx) = \frac{\pi - x}{2} \label{E:five2} \, .
\end{equation}
Special case: substituting $x = 1$, this shows that both sums in Equation \eqref{E:Prob6241a} equal $(\pi - 1)/2$.

This completes the proof of Theorem \ref{T:Thm1}.

We can use polylogarithms to directly \emph{derive} the sum $(\pi - 1)/2$ in Theorem \ref{T:Thm1}; see Section \ref{S:Dilogarithms}.

Here is \textit{Mathematica} code that uses Equation \eqref{E:bnOdd} to evaluate the two integrals in \eqref{E:TwoPieceIntegral}:
\begin{verbatim}
  i1 = Integrate[Sin[n x] * x * (Pi - 1)/2, {x, 0, 1}]
  i2 = Integrate[Sin[n x] * (Pi - x)/2, {x, 1, Pi}]
  total = Simplify[(2/Pi)*(i1 + i2), Assumptions -> Element[n, Integers]]
\end{verbatim}
(``\verb+Assumptions -> Element[n, Integers]+'' tells \textit{Mathematica} that $n$ is an integer, so \textit{Mathematica} can do things like replace $\sin(n \pi)$ with 0).
The result is that \verb+total+ equals $\sin(n)/n^2$.


\begin{theorem} \label{T:Thm1A}
\begin{equation} \label{E:Prob6241b}
\sum_{n=1}^{\infty} \frac{\sin^{2}(n)}{ n^{4} } = \frac{(\pi - 1)^{2} }{6} \, .
\end{equation}
\end{theorem}

\textbf{Proof.}
The $n^{\text{th}}$ Fourier sine coefficient in Equation \eqref{E:gxTentative} is $\sin(n)/n^2$.
The \emph{square} of this coefficient is just the $n^{\text{th}}$ term in Equation \eqref{E:Prob6241b}.
This suggests that Parseval's equation \eqref{E:GenParsevalEqn} might be useful.

We will prove Equation \eqref{E:Prob6241b} by applying Parseval's equation to the Fourier series for $g(x)$ in Equation \eqref{E:gxTentative}.
With $a = -\pi, b = \pi$, the period $p$ is $p = b - a = 2 \pi$, so $2/p = 1/\pi$.
Because $g(x)$ is an odd function, for any $x$, we have $g(-x)^2 = g(x)^2$, so we can integrate just from 0 to $\pi$ and multiply the integral in \eqref{E:GenParsevalEqn} by 2.
All $a_n$ are zero and $b_n = \sin(n)/n^2$, so
\begin{align}\label{E:ThreePartParseval}
 & \sum_{n=1}^{\infty} b_{n}^2 = \sum_{n=1}^{\infty} \frac{\sin^{2}(n)}{ n^{4} } =
 \frac{1}{\pi} \int_{-\pi}^{\pi} (g(x))^2 \, dx =
 \frac{2}{\pi} \int_{0}^{\pi} (g(x))^2 \, dx \\
 & = \frac{2}{\pi} \left(
   \int_{0}^{1} \left( \frac{x(\pi - 1)}{2} \right)^2 dx +
   \int_{1}^{\pi} \left( \frac{\pi - x}{2} \right)^2 dx
  \right) \notag \\
  & = \frac{2}{\pi} \left( \frac{(\pi - 1)^2}{12} + \frac{(\pi - 1)^3}{12} \right)
  = \frac{2}{\pi} \left(  \frac{\pi^3 - 2 \pi^2 + \pi}{12}  \right)
   \notag \\
  & = \frac{(\pi - 1)^2}{6} \notag \, .
\end{align}

Here is the \textit{Mathematica} code that carries out this somewhat tedious calculation:
\begin{verbatim}
  i1 = Integrate[( x*(Pi - 1)/2 )^2, {x, 0, 1}]
  i2 = Integrate[( (Pi - x)/2 )^2,   {x, 1, Pi}]
  total = Simplify[(2/Pi)*(i1 + i2)]
\end{verbatim}
The result is that \verb+total+ equals $(\pi - 1)^2/6$.

This completes the proof of Theorem \ref{T:Thm1A}.

\textbf{Discussion.} \textbf{1.}
Equations \eqref{E:Prob6241a} and \eqref{E:Prob6241b} originally appeared in \cite{Baillie6241}.
Equation \eqref{E:Prob6241a} is a nice example of a series with the counterintuitive property that $\sum a_{n} = \sum a_{n}^{2}$.
It is also the first of several examples we will encounter in which we can multiply the $n^{\text{th}}$ term of a series by $\sin(n)/n$ without changing the sum.
Equation \eqref{E:Prob6241b} is a pretty variation on the classic formula
\[
\sum_{n=1}^{\infty} \frac{1}{n^2} = \frac{\pi^{2}}{6} \, .
\]

\begin{figure}[ht] 
\centering
\begin{minipage}[t]{\figureMinipageWidth}
\centering
\includegraphics[width=\picDblWidth]{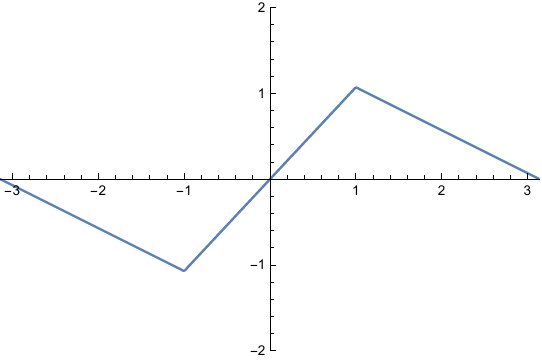}
\caption{$g(x), -\pi \leq x \leq \pi$}
\label{fig:gxLargeInterval1}
\end{minipage}\hfill
\begin{minipage}[t]{\figureMinipageWidth}
\centering
\includegraphics[width=\picDblWidth]{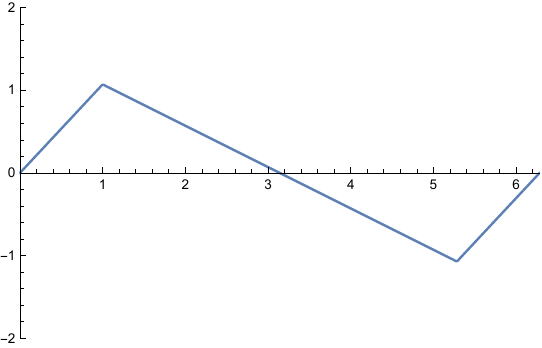}
\caption{$g(x), 0 \leq x \leq 2 \pi$}
\label{fig:gxLargeInterval2}
\end{minipage}
\end{figure}

\textbf{2.}
In Figures \ref{fig:gxLargeInterval1} and \ref{fig:gxLargeInterval2}, we plot the series
\[
g(x) = \sum_{n=1}^{\infty} \frac{\sin(n)}{n^2} \sin(nx)
\]
over larger intervals.
We used two intervals (0 to 1 and 1 to $\pi$) and the two expressions for $g(x)$ in Equation \eqref{E:gxTentative2} to get the Fourier series for $g(x)$.
But over these larger intervals, $g(x)$ consists of \emph{three} segments.
We could have used the three segments in these figures to derive the Fourier series.
For example, over the interval from $-\pi$ to $\pi$, we have
\begin{equation*} 
g(x) =
  \begin{cases}  
    -(\pi + x)/2  &\text{for $-\pi \leq x < -1$,} \\
    x(\pi - 1)/2  &\text{for $-1 \leq x \leq 1$,} \\
    (\pi - x)/2  &\text{for $1 < x \leq \pi$.}
  \end{cases}
\end{equation*}
Alternatively, over the interval from 0 to $2 \pi$, we have
\begin{equation*} 
g(x) =
  \begin{cases}  
    x(\pi - 1)/2  &\text{for $0 \leq x < 1$,} \\
    (\pi - x)/2  &\text{for $1 \leq x \leq 2 \pi - 1$,} \\
    (x - 2 \pi)(\pi - 1)/2  &\text{for $2 \pi - 1 < x \leq 2 \pi$.}
  \end{cases}
\end{equation*}
Instead of the two integrals in \eqref{E:TwoPieceIntegral}, we would have three, but the result would be the same.

\textbf{3.}  
It's unusual to see $\sin(n)$ as part of the \emph{coefficient} of $\sin(nx)$ in a Fourier series.
Where did it come from?
Most examples of Fourier series in textbooks involve functions over intervals with endpoints such as 0, $\pi/2$, and $\pi$.
The series for $g(x)$ in the above proof contains $\sin(n)$ because one of the ``corners'' of $g(x)$ occurred at $x = 1$.
To see this, let $b$ be some number between 0 and $\pi$.
Here's a function, $h(x)$, with the same overall shape as $g(x)$, but with a ``corner'' at $x = b$:
\[
h(x)=
  \begin{cases}
    x(\pi - b)/2  &\text{for $0 \leq x \leq b$,}\\
    b(\pi - x)/2  &\text{for $b < x \leq \pi$.}
  \end{cases}
\]

$h(x)$ is a generalization of $g(x)$: if $b = 1$, then $h(x) = g(x)$.

In \textit{Mathematica}, the coefficient of $\sin(n x)$ can be calculated with:
\begin{verbatim}
  i1 =  Integrate[Sin[n x] * x*(Pi - b)/2, {x, 0, b}]
  i2 = Integrate[Sin[n x] * b*(Pi - x)/2, {x, b, Pi}]
  total = Simplify[(2/Pi)*(i1 + i2), Assumptions -> Element[n, Integers]]
\end{verbatim}

The result is $\sin(bn)/n^2$, so the Fourier series for $h(x)$ is
\[
h(x) = \sum_{n=1}^{\infty} \frac{\sin(bn)}{n^2} \sin(nx) \, .
\]

If $b = 1$, we get $\sin(n)$ in the coefficient.
If $b$ is the $x$ value of an endpoint commonly found in textbooks, for example $b = \pi/2$, then $\sin(bn)$ reduces to 1, $-1$, or 0.
Later, we will again use this ``trick'' to put $\sin(n)$ into the $n$\textsuperscript{th} Fourier coefficient.

\textbf{4.}
Among the series in which we can multiply the $n$\textsuperscript{th} term by $\sin(n)/n$ without changing the sum is this classic Gregory-Liebniz series \cite{BGL}:
\[
\frac{\pi }{4} = 1 - \frac{1}{3} + \frac{1}{5} - \frac{1}{7} + \dots \, .
\]
Let's substitute $x = \pi /2$ into Equations \eqref{E:five1} and \eqref{E:five2}.
First, $(\pi - x)/2 = (\pi - (\pi/2))/2 = \pi/4$.
Next, when $n$ is even, $\sin(n\pi /2) = 0$ and when $n$ is odd, $\sin(n \pi /2) = \pm 1$.
For $x = \pi/2$, Equation \eqref{E:five1} gives
\[
  \frac{\pi }{4} = 1 - \frac{1}{3} + \frac{1}{5} - \frac{1}{7} + \dots \, .
\]
Likewise, with $x = \pi/2$, Equation \eqref{E:five2} becomes
\[
 \frac{\pi }{4} = 1\cdot \frac{\sin(1)}{1} - \frac{1}{3} \cdot \frac{\sin(3)}{3} + \frac{1}{5} \cdot \frac{\sin(5)}{5} - \frac{1}{7} \cdot \frac{\sin(7)}{7} + \dots \, .
\]

So, we get this amusing variant of the Gregory-Liebniz series:
\begin{equation}\label{E:GregoryGen1}
\frac{\pi }{4} = 1 - \frac{1}{3} + \frac{1}{5} - \frac{1}{7} + \dots =
1\cdot \frac{\sin(1)}{1} - \frac{1}{3} \cdot \frac{\sin(3)}{3} + \frac{1}{5} \cdot \frac{\sin(5)}{5} - \frac{1}{7} \cdot \frac{\sin(7)}{7} + \dots \, .
\end{equation}

We can also prove this more general result: for every non-zero $x$ in $-\pi /2 \leq x \leq \pi /2$,
\begin{equation}\label{E:GregoryGen2}
\frac{\pi }{4} = 1 - \frac{1}{3} + \frac{1}{5} - \frac{1}{7} + \dots =
1\cdot \frac{\sin(x)}{x} - \frac{1}{3} \cdot \frac{\sin(3x)}{3x} + \frac{1}{5} \cdot \frac{\sin(5x)}{5x} - \frac{1}{7} \cdot \frac{\sin(7x)}{7x} + \dots \, .
\end{equation}

To prove Equation \eqref{E:GregoryGen2}, consider the function defined by
\begin{equation}\label{E:GregoryGen3}
f(x)=
  \begin{cases}
    \pi x/4         &\text{for $0 \leq x \leq \pi/2$,}\\
    \pi(\pi - x)/4  &\text{for $\pi/2 < x \leq \pi$.}
  \end{cases}
\end{equation}
The Fourier sine coefficients for this function can be computed with this \textit{Mathematica} code:
\begin{verbatim}
  i1 = Integrate[ Sin[n x] * Pi x/4,        {x, 0, Pi/2}]
  i2 = Integrate[ Sin[n x] * Pi (Pi - x)/4, {x, Pi/2, Pi}]
  total = Simplify[(2/Pi) * (i1 + i2), Assumptions -> Element[n, Integers]]
\end{verbatim}
The result is $\sin(n\pi/2)/n^2$, which is 0 if $n$ is even, and the Fourier sine series of $f(x)$ is
\begin{equation}\label{E:GregoryGen4}
f(x) = \sum_{n=1}^{\infty} \frac{\sin(n\pi /2)}{n^{2} } \sin(nx) =
 \frac{\sin(x)}{1^{2} } - \frac{\sin(3x)}{3^{2} } + \frac{\sin(5x)}{5^{2} } - \cdots \, .
\end{equation}
This is actually valid over the larger interval $-\pi \leq x \leq \pi$ because $f(x)$ is an odd function.
Therefore, $f(x) = \pi x/4$ for $-\pi/2 \leq x \leq \pi/2$.
Substituting $x = 1$ into $f(x)$ gives the variant of the Gregory-Leibniz series stated above.

From the definition of $f(x)$ in Equation \eqref{E:GregoryGen3}, note that, for non-zero $x$ in $-\pi/2 \le x \le \pi/2$, $f(x)/x$ has the constant value $\pi /4$.
For $x \neq 0$, divide $f(x)$ and its Fourier series \eqref{E:GregoryGen4} through by $x$. For non-zero x in $x$ in $-\pi/2 \le x \le \pi/2$,
\[
  \frac{f(x)}{x} = \frac{\pi}{4} = 1\cdot \frac{\sin(x)}{x} -\frac{1}{3} \cdot \frac{\sin(3x)}{3x} +\frac{1}{5} \cdot \frac{\sin(5x)}{5x} - \cdots \, .
\]
We can also write this series in closed form, as follows:
\[
\frac{\pi}{4} = \sum_{n=1}^{\infty} \frac{(-1)^{n+1} }{2n-1} \cdot \frac{\sin((2n-1)x)}{(2n-1)x} \, .
\]

Recall from Section \ref{S:Preliminaries} that the \emph{sinc} function is defined as sinc$(0) = 1$, and sinc$(x) = \sin(x)/x$ if $x \neq 0$.
Using the sinc function, we can remove the condition in the previous equation that $x$ be non-zero.
For \emph{all} $x$ in $-\pi/2 \le x \le \pi/2$, we have
\[
\frac{\pi}{4} = 1\cdot \text{sinc}(x) -\frac{1}{3} \cdot \text{sinc}(3x) +\frac{1}{5} \cdot \text{sinc}(5x) - \cdots \, .
\]

Perhaps it should be pointed out that the fact that this holds for $x < 0$ is not really very interesting because $\text{sinc}(-x) = \text{sinc}(x)$.
In similar cases later on, we will sometimes not bother to say that an equation also holds for $x < 0$.

\textbf{5.}
Look again at the two series in Equation \eqref{E:GregoryGen1}.
In the series on the left, the terms are at most $1/(2 n - 1)$ in absolute value.
(This series is well-known for its \emph{slow} convergence to the sum, $\pi/4$.)
In the series on the right, the terms are at most $1/(2 n - 1)^2$ in absolute value.
One might suspect that the series on the right converges faster.
It does: if we add one million terms of each series, the partial sums on the left and right
differ from $\pi/4$ by about $2.5 \cdot 10^{-7}$ and $1.5 \cdot 10^{-13}$.
Nevertheless, because the series on the right requires the computation of many values of the $\sin$ function,
this series is not useful for computing $\pi$.

\textbf{6.}
Published proofs often obscure the method by which theorems, or proofs, were discovered.
Instead, things often appear to be pulled out of thin air for no apparent reason, and they magically seem to work.
In this author's view, this is an unfortunate practice (of which this author is also guilty).
So, it is worth pointing out that Equations \eqref{E:Prob6241a}, \eqref{E:Prob6241b}, and \eqref{E:five1} and \eqref{E:five2} were found and proved by the procedure in the above proof:
First, a numerical calculation suggested that the sum of the second series in \eqref{E:Prob6241a} was very close to $(\pi -1)/2$.
Later, the author noticed in a textbook that the sum of the first series in \eqref{E:Prob6241a} equals $(\pi -1)/2$.
This surprising observation meant that there was something interesting going on that merited further investigation.
Then the author realized that the second series in \eqref{E:Prob6241a} was just the Fourier series $\sum \frac{\sin(n)}{n^{2} }  \sin(nx)$, evaluated at $x = 1$.
Plotting this Fourier series suggested that the function it represented consisted of linear pieces that could be determined from the graph.
Finally, the pieces were integrated, as above, to compute the Fourier coefficients, thus verifying the Fourier series, and proving that the sums were equal.

\textbf{7.}
We've just found a series where $\sum a_{n} = \sum a_{n}^{2}$.
It's hard not to mention series with another curious property: $(\sum a_{n} )^{2} = \sum a_{n}^{2}$.
The square of this series
\[
\frac{\pi }{\sqrt{8} } = 1 + \frac{1}{3} - \frac{1}{5} - \frac{1}{7} + \frac{1}{9} + \frac{1}{11} \dots ++-- \dots
\]
is
\[
\frac{\pi ^{2} }{8} = 1 + \frac{1}{3^{2} } + \frac{1}{5^{2} } + \frac{1}{7^{2} } + \frac{1}{9^{2} } + \frac{1}{11^{2} } + \dots \, .
\]
That is, the first series can be squared by adding the squares of the individual terms!
These series are exercises 6(a) and 6(b) in \cite[page 503]{Kaplan}:
\[
 k_1(x) =
   \frac{1}{2} + \frac{2}{\pi} \left( \sum_{n=1}^{\infty} \frac{\sin((2 n - 1) x)}{2n - 1} \right) =
   \begin{cases}  
     0  &\text{for $-\pi \leq x < 0$,}\\
     1  &\text{for $0 \leq x \leq \pi$.}
   \end{cases}
\]
\[
 k_2(x) =
   \frac{\pi}{2} - \frac{4}{\pi} \left( \sum_{n=1}^{\infty} \frac{\cos((2 n - 1) x)}{(2n - 1)^2} \right) =
   \begin{cases}  
              -x  &\text{for $-\pi \leq x \leq 0$,}\\
    \phantom{-}x  &\text{for $0 \leq x \leq \pi$.}
   \end{cases}
\]
(Hint: Substitute $x = \pi/4$ in $k_1(x)$ and $x = 0$ in $k_2(x)$.)
Are there other examples like this?

\textbf{8.}
We have already seen that
\[
 \sum_{n=1}^{\infty} \frac{\sin(n)}{n} =
 \sum_{n=1}^{\infty} \left( \frac{\sin(n)}{n} \right) ^{2} =
 \frac{\pi - 1}{2}.
\]
It would be surprising if
\[
 \sum_{n=1}^{\infty} \left( \frac{\sin(n)}{n} \right) ^{3}
\]
had the same sum.
Unfortunately, this third sum is different.
This can easily be seen by numerically adding, say, 1 million terms of this series.
The sum is about
\[
  0.67809 \text{ } 72450 \text{ } 96172 \text{ } 46442 \, .
\]
\textit{Mathematica} claims this third sum is \emph{exactly} $(3 \pi - 4)/8$, which is about
\[
  0.67809 \text{ } 72450 \text{ } 96172 \text{ } 46376 \, .
\]
We can get the sum $(3 \pi - 4)/8$ by substituting $x = 1$ into Equation \eqref{E:unknownFxQuadratic};
we will prove that Equation in Section \ref{S:Recovering}.

We have seen that we can get nice series for $\pi$ by substituting an appropriate value of $x$ into an appropriate Fourier series.
But observe that we also did something else in the proof of Theorem \ref{T:Thm1}:
we used two functions that were distinct over $(0, \pi)$, but which were equal over the entire subinterval $[1, \pi]$.
The Fourier series were different, but the two series yielded the same values over that subinterval.
Of course, there are many such pairs of functions, but in what follows, we will use a few that lead to interesting results.
The reader is encouraged to experiment and look for other such functions.

\textbf{9.}
Serge Ballif \cite{Ballif} showed how to produce other series with the property that
\begin{equation} \label{E:BallifMK}
\left( \sum_{n=0}^{\infty} a_n \right)^m = \sum_{n=0}^{\infty} a^k_n \, ,
\end{equation}
where $m$ and $k$ are $\geq 1$ and $m \neq k$.
The series in Equation \eqref{E:Prob6241a} is of this form, with $m = 1$ and $k = 2$.
As a special case $(m = 1, k = 2)$, he shows that
\[
\sum_{n=0}^{\infty} (1 + r) r^n = \sum_{n=0}^{\infty} \left( (1 + r) r^n \right)^2 = \frac{1+r}{1-r} \, ,
\]
provided $-1 \leq r < 1$.
For example, setting $r = 1/2$, we have
\[
\sum_{n=0}^{\infty} \left( \frac{3}{2} \cdot \frac{1}{2^n} \right) = \sum_{n=0}^{\infty} \left( \frac{3}{2} \cdot \frac{1}{2^n} \right)^2 = 3 \, .
\]

Following Ballif's recipe with $m = 1$ and $k = 3$, another example is
\begin{equation} \label{E:Ballif3}
  \sum_{n = 0}^{\infty} \left( \frac{13}{8} \cdot \left( \frac{7}{8} \right)^n \right) =
  \sum_{n = 0}^{\infty} \left( \frac{13}{8} \cdot \left( \frac{7}{8} \right)^n \right)^3 = 13 \, .
\end{equation}
Equation \eqref{E:SumAndSumCubed} also has the property that $\sum a_n = \sum (a_n)^3$.

\ifthenelse {\boolean{BKMRK}}
  { \section{Sum of \texorpdfstring{$\sin(2 n x)/(2 n)$}{sin(2 n x)/(2 n)}}\label{S:sin2n} }
  { \section{Sum of $\sin(2 n x)/(k n)$}\label{S:sin2n} }

Look again at the sums in Equations \eqref{E:five1} and \eqref{E:five2}.
Is there some way we can we ``tweak'' them to get other interesting results?
One approach will be discussed later, in Theorem \ref{T:ThmFourEqualSums}.
For now, let's see what happens if we take the sums just over \emph{even} $n$.

Is the following true?
\[
  \sum_{n=1}^{\infty} \frac{\sin(2n)}{2n} \overset{?}{=} \sum_{n=1}^{\infty} \left( \frac{\sin(2n)}{2n} \right) ^2 \, .
\]

First, let's look at some numerical data.
The sums of the first million terms of these series are approximately 0.28539 78783 and 0.28539 80384, respectively.
\emph{Both} values are close to $(\pi - 2)/4 \approx$ 0.28539 81634.
Thus, it \emph{appears} that
\[
  \sum_{n=1}^{\infty} \frac{\sin(2n)}{2n} = \sum_{n=1}^{\infty} \left( \frac{\sin(2n)}{2n} \right) ^2 = \frac{\pi - 2}{4} \, .
\]
Let's see if we can prove this.

\begin{figure}[ht]
  \mbox{\includegraphics[width=\picSingleWidth]{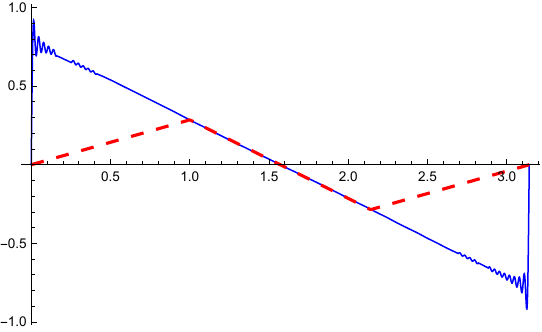}}
  \caption{$\sum_{n=1}^{100} \frac{\sin(2nx)}{2n}$ (solid) and $\sum_{n=1}^{100} \frac{\sin(2n)}{2n} \cdot \frac{\sin(2nx)}{2n}$ (dashed)}
  \label{fig:evenSincs}
\end{figure}

We will use the same trick we used earlier: consider the sums to be Fourier series evaluated at $x = 1$.
That is, we will consider the two series
\[
\sum_{n=1}^{\infty} \frac{\sin(2n x)}{2n} \quad \text{ and } \quad
\sum_{n=1}^{\infty} \frac{\sin(2n)}{2n} \frac{\sin(2nx)}{2n} \, .
\]
The sums of the first 100 terms are the solid (blue) and dashed (red) graphs in Figure \ref{fig:evenSincs}.

The two sums \emph{seem} to agree for every $x$ in the entire interval from $x = 1$ to $x = \pi - 1$ !

Reading the graphs carefully suggests that, for $1 \leq x \leq \pi - 1$, we appear to have
\[
\sum_{n=1}^{\infty} \frac{\sin(2nx)}{2n} = \sum_{n=1}^{\infty}  \frac{\sin(2n)}{2n} \frac{\sin(2nx)}{2n} = \frac{\pi - 2x}{4} \, .
\]
As far as we can tell from this plot, the solid graph looks like $(\pi - 2x)/4$.

Let's define $f(x) = (\pi - 2 x)/4$ for $0 < x < \pi$, with period $2 \pi$.
Let's further define $f(x)$ to be \emph{odd}, so it can be represented as a sum only of sines, without cosine terms.
Then we can use Equation \eqref{E:bnOdd} to compute the sine series of $f(x)$:
\[
\frac{2}{\pi} \int_{0}^{\pi} f(x) \sin(nx) \, dx = \frac{1 + (-1)^n}{2 n} \, .
\]
If $n$ is odd, then $1 + (-1)^n = 0$.
If $n$ is even, then $(1 + (-1)^n)/(2n)$ equals $1/2$, $1/4$, $1/6$, etc.
Therefore, for $0 < x < \pi$, the sine series for $f(x)$ is
\begin{equation} \label{E:sumSincEven}
\frac{\pi - 2x}{4} = \sum_{n=1}^{\infty} \frac{\sin(2nx)}{2n} \, .
\end{equation}

Similarly, the dashed lines in Figure \ref{fig:evenSincs} \emph{seem} to be these linear pieces:
\begin{equation}\label{E:fEven}
g(x)=
  \begin{cases}
    x(\pi - 2)/4          &\text{for $0 \leq x < 1$,}\\
    (\pi - 2x)/4          &\text{for $1 \leq x \leq \pi - 1$,}\\
    (x - \pi)(\pi - 2)/4  &\text{for $\pi - 1 < x \leq \pi$.}
  \end{cases}
\end{equation}
Let's also assume that $g(x)$ has period $2 \pi$ and that $g(x)$ is odd.
We use Equation \eqref{E:bnOdd} to compute the coefficient of $\sin(nx)$ in the sine series for $g(x)$.
We will need three separate integrals.
The result is:
\begin{align*}
& \frac{2}{\pi} \int_{0}^{\pi} g(x) \sin(nx) \, dx \\
& = \frac{2}{\pi} \int_{0}^{1} \frac{x (\pi - 2)}{4} \sin(nx) \, dx +
\frac{2}{\pi} \int_{1}^{\pi - 1} \frac{(\pi - 2 x)}{4} \sin(nx) \, dx +
\frac{2}{\pi} \int_{\pi - 1}^{\pi} \frac{(x - \pi) (\pi - 2)}{4} \sin(nx) \, dx \\
& = \frac{(1+(-1)^{n} )\sin(n)}{2n^{2} } \, .
\end{align*}
If $n$ is odd, this expression is 0.
If $n$ is even, then $1+(-1)^n = 2$, so for $n = 2$, $n = 4$, etc., the expression becomes $\sin(2)/2^2$, $\sin(4)/4^2$, etc.

Here is the \textit{Mathematica} code that gives this result:
\begin{verbatim}
  i1 = Integrate[Sin[n x] * x*(Pi - 2)/4, {x, 0, 1}]
  i2 = Integrate[Sin[n x] * (Pi - 2 x)/4, {x, 1, Pi - 1}]
  i3 = Integrate[Sin[n x] * (x - Pi)*(Pi - 2)/4, {x, Pi - 1, Pi}]
  total1 = (2/Pi)*(i1 + i2 + i3)
  total = Simplify[total1, Assumptions -> Element[n, Integers]]
\end{verbatim}

Therefore, the function $g(x)$ defined in Equation \eqref{E:fEven} does, indeed, have the sine series
\[
\sum_{n=1}^{\infty} \frac{(1+(-1)^{n} )\sin(n)}{2n^2} \sin(nx) =
\frac{\sin(2)}{2^2}\sin(2x) + \frac{\sin(4)}{4^2}\sin(4x) + \cdots
=\sum_{n=1}^{\infty} \frac{\sin(2n)}{2n}  \frac{\sin(2nx)}{2n} \, .
\]

We conclude that, for $0 \leq x \leq \pi$,
\[
g(x) = \sum_{n=1}^{\infty} \frac{\sin(2n)}{2n} \cdot \frac{\sin(2nx)}{2n} \, .
\]
From this, and from Equation \eqref{E:sumSincEven}, for $1 \leq x \leq \pi - 1$, we have
\begin{equation}\label{E:evenTerms}
\sum_{n=1}^{\infty} \frac{\sin(2nx)}{2n} = \sum_{n=1}^{\infty} \frac{\sin(2n)}{2n} \frac{\sin(2nx)}{2n} = \frac{\pi - 2x}{4} \, .
\end{equation}

This is just like Equations \eqref{E:five1} and \eqref{E:five2}, which are sums over \emph{all} (i.e., both even and odd) positive integers $n$:
\begin{equation}\label{E:allTerms}
\sum_{n=1}^{\infty} \frac{\sin(nx)}{n} = \sum_{n=1}^{\infty} \frac{\sin(n)}{n} \frac{\sin(nx)}{n} = \frac{\pi - x}{2} ,
\end{equation}
If we subtract Equation \eqref{E:evenTerms} from \eqref{E:allTerms}, we are left with just the terms for \emph{odd} $n$.
Also,
\[
\frac{\pi - x}{2} - \frac{\pi - 2x}{4} = \frac{\pi}{4} \, .
\]
We get
\begin{equation}\label{E:oddTerms}
\sum_{n=1}^{\infty} \frac{\sin((2n-1)x)}{2n-1} = \sum_{n=1}^{\infty} \frac{\sin(2n-1)}{2n-1} \frac{\sin((2n-1)x)}{2n-1} = \frac{\pi}{4} \, .
\end{equation}
This is true wherever \emph{both} Equations \eqref{E:evenTerms} and \eqref{E:allTerms} are true, namely, for $1 \leq x \leq \pi - 1$.

We have therefore proved

\begin{theorem}\label{T:Thm2}
For $1 \leq x \leq \pi - 1$, we have
\begin{equation} \label{E:Sin2nxOver2N}
\sum_{n=1}^{\infty} \frac{\sin(2nx)}{2n} = \sum_{n=1}^{\infty} \frac{\sin(2n)}{2n} \frac{\sin(2nx)}{2n} = \frac{\pi -2x}{4}
\end{equation}
and
\begin{equation} \label{E:SinOddNxOverOddN}
\sum_{n=1}^{\infty} \frac{\sin((2n-1)x)}{2n-1} = \sum_{n=1}^{\infty} \frac{\sin(2n-1)}{2n-1} \frac{\sin((2n-1)x)}{2n-1} = \frac{\pi }{4} \, .
\end{equation}
\end{theorem}

\textbf{Corollary.}
Evaluating the sums in Theorem \ref{T:Thm2} at $x = 1$, we get
\begin{equation} \label{E:SincEvenN}
\sum_{n=1}^{\infty} \frac{\sin(2n)}{2n} = \sum_{n=1}^{\infty} \left(\frac{\sin(2n)}{2n} \right) ^{2} = \frac{\pi -2}{4}
\end{equation}
and
\begin{equation} \label{E:SincOddN}
\sum_{n=1}^{\infty} \frac{\sin(2n-1)}{2n-1} = \sum_{n=1}^{\infty} \left(\frac{\sin(2n-1)}{2n-1} \right) ^{2} = \frac{\pi }{4} \, .
\end{equation}
So, the peculiar property of Equation \eqref{E:Prob6241a} that $\sum a_n = \sum (a_n)^2$ holds for the even and odd terms separately!

\begin{figure}[ht]
  \mbox{\includegraphics[width=\picSingleWidth]{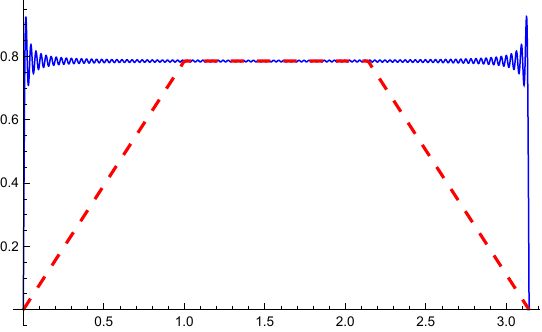}}
  \caption{$\sum_{n=1}^{100} \frac{\sin((2n-1)x)}{2n-1}$ (solid) and $\sum_{n=1}^{100}  \frac{\sin(2n-1)}{2n-1} \cdot \frac{\sin((2n-1)x)}{2n-1}$ (dashed)}
  \label{fig:oddSincs}
\end{figure}

\textbf{Discussion 1.}
We obtained the sum in \eqref{E:SinOddNxOverOddN} by subtracting Equation \eqref{E:Sin2nxOver2N} (even-numbered terms) from Equation \eqref{E:allTerms} (all terms).
We could also analyze the odd terms directly, just as we did above for the even terms.
The relevant graphs are shown in Figure \ref{fig:oddSincs}.
The dashed function is given by
\begin{equation}\label{E:gOdd}
g(x)=
  \begin{cases}
    \pi x/4         &\text{for $0 \leq x < 1$,}\\
    \pi/4           &\text{for $1 \leq x \leq \pi - 1$,}\\
    \pi(\pi - x)/4  &\text{for $\pi - 1 < x \leq \pi$.}
  \end{cases}
\end{equation}
For this $g(x)$, the following \textit{Mathematica} code calculates the coefficient of $\sin(nx)$:
\begin{verbatim}
  trigRules = { Sin[n*(c_ + Pi)] -> Cos[n Pi] * Sin[n c],
                Cos[n*(c_ + Pi)] -> Cos[n Pi] * Cos[n c] }
  i1  = Integrate[Sin[n x] * Pi x/4, {x, 0, 1}]
  i2Temp = Integrate[Sin[n x] * Pi/4,   {x, 1, Pi - 1}]
  i2Temp2 = i2Temp //. trigRules
  i2 = Simplify[i2Temp2, Assumptions -> Element[n, Integers]]
  i3Temp = Integrate[Sin[n x] * Pi (Pi - x)/4, {x, Pi - 1, Pi}]
  i3Temp2 = i3Temp //. trigRules
  i3 = Simplify[i3Temp2, Assumptions -> Element[n, Integers]]
  total = Simplify[(2/Pi)*(i1 + i2 + i3)]
\end{verbatim}
The $n$\textsuperscript{th} coefficient is
\[
\frac{\left(1 - (-1)^n\right) \sin(n)}{2 n^2} \, .
\]
For even $n$, the coefficients are 0, and for $n = 1$, $n = 3$, and $n = 5$, the coefficients are $\sin(1)$, $\sin(3)/3^2$, and $\sin(5)/5^2$.
We conclude that, for $0 < x < \pi$,
\[
g(x) = \sum_{n=1}^{\infty} \sin((2n-1)x) \cdot \frac{\sin(2n-1)}{(2n-1)^2} \, .
\]

\textbf{2.}

Let's apply Parseval's equation \eqref{E:ParsevalsEquation}, using the interval from $x = 0$ to $x = \pi$, to $f(x)$ and $g(x)$.
\[
f(x) = \sum_{n=1}^{\infty} \sin(2nx) \cdot \frac{\sin(2n)}{(2n)^2}
\]
so
\[
\sum_{n=1}^{\infty} \frac{\sin ^{2} (2n)}{(2n)^{4} } = \frac{2}{\pi}\int_0^{\pi} f(x)^2 \, dx \, .
\]

Just as we did in Equation \eqref{E:ThreePartParseval}, we must break up this integral into three pieces:

\begin{align*}
 \frac{2}{\pi}\int_0^{\pi} f(x)^2 \, dx & = \frac{2}{\pi} \left (
\int_0^1 f(x)^2 dx  + \int_1^{\pi-1} f(x)^2 \, dx + \int_{\pi-1}^{\pi} f(x)^2 \, dx
\right ) \\
 &
 = \frac{2}{\pi} \left( \frac{(\pi-2)^2}{48} + \frac{(\pi-2)^3}{48} + \frac{(\pi-2)^2}{48} \right)
  = \frac{(\pi - 2)^{2} }{24} \, .
\end{align*}

The result is
\[
\sum_{n=1}^{\infty} \frac{\sin ^{2} (2n)}{(2n)^{4} } = \frac{(\pi - 2)^{2} }{24} \, .
\]
Similarly, applying Parseval's equation to
\[
g(x) = \sum_{n=1}^{\infty} \sin((2n-1)x) \cdot \frac{\sin(2n-1)}{(2n-1)^2} \, ,
\]
we get
\[
\sum_{n=1}^{\infty} \frac{\sin ^{2} (2n-1)}{(2n-1)^{4} } = \frac{\pi ^{2} }{8} - \frac{\pi }{6} \, .
\]

\ifthenelse {\boolean{BKMRK}}
  { \section{Sum of \texorpdfstring{$\sin(3 n x)/(3 n)$}{sin(3 n x)/(3 n)}}\label{S:sin3n} }
  { \section{Sum of $\sin(3 n x)/(3 n)$}\label{S:sin3n} }

We have seen in Equations \eqref{E:Prob6241a} and \eqref{E:SincEvenN} that
\[
\sum_{n=1}^{\infty} \frac{\sin(n)}{n} = \sum_{n=1}^{\infty} \left( \frac{\sin(n)}{n} \right)^2 = \frac{\pi - 1}{2}
\]
and
\[
\sum_{n=1}^{\infty} \frac{\sin(2n)}{2n} = \sum_{n=1}^{\infty} \left( \frac{\sin(2n)}{2n} \right)^2 = \frac{\pi - 2}{4} \, .
\]
Is the following true?
\[
\sum_{n=1}^{\infty} \frac{\sin(3n)}{3n} = \sum_{n=1}^{\infty} \left( \frac{\sin(3n)}{3n} \right)^2 = \frac{\pi - 3}{6} \, .
\]

First, let's check numerically whether the sums are at least approximately equal. If we add one million terms of each series, we get
\begin{align*}
  \sum_{n=1}^{1000000} \frac{\sin(3n)}{3n} & \approx 0.02359 \text{ } 86235 \, , \\
  \sum_{n=1}^{1000000} \left( \frac{\sin(3n)}{3n} \right)^2 & \approx 0.02359 \text{ } 87200 \, . \\
  \text{Also,} \quad \frac{\pi - 3}{6} & \approx 0.02359 \text{ } 87756 \, .
\end{align*}

This is reasonably good agreement, given that we added only one million terms.
So, let's tentatively call this a Theorem and try to prove it.

\begin{theorem}\label{T:ThmSum3nx}
\begin{equation}\label{E:sin3n}
\sum_{n=1}^{\infty} \frac{\sin(3n)}{3n} = \sum_{n=1}^{\infty} \left( \frac{\sin(3n)}{3n} \right)^2 = \frac{\pi - 3}{6} \, .
\end{equation}
\end{theorem}

\begin{figure}[ht]
\centering
\begin{minipage}[t]{\figureMinipageWidth}
\centering
\includegraphics[width=\picDblWidth]{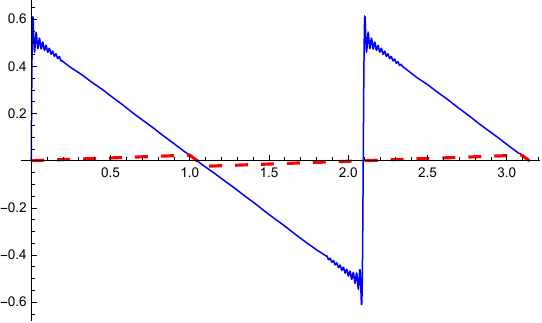}
\caption{$f(x)$ (solid) and $g(x)$ (dashed)}
\label{fig:f3g3x}
\end{minipage}\hfill
\begin{minipage}[t]{\figureMinipageWidth}
\centering
\includegraphics[width=\picDblWidth]{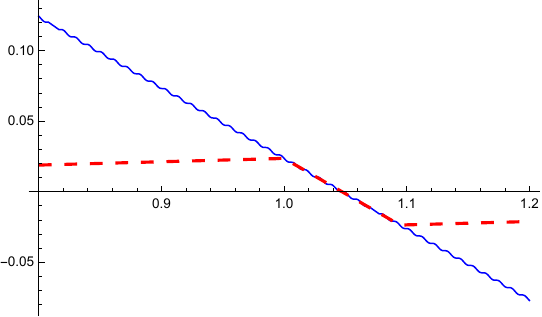}
\caption{Figure \ref{fig:f3g3x}, near $x = 1$}
\label{fig:f3g3xCloseup}
\end{minipage}
\end{figure}

\textbf{Proof.}
As above, let's convert the sums into functions of $x$ by assuming that one factor of $\sin(3 n)$ comes from $\sin(3 n x)$ evaluated at $x = 1$.
Then the above sums become special cases of these functions of $x$:
\begin{equation}\label{E:f3x}
f(x) = \sum_{n=1}^{\infty} \frac{\sin(3nx)}{3n}
\end{equation}
\begin{equation}\label{E:g3x}
g(x) = \sum_{n=1}^{\infty} \frac{\sin(3n)}{3n} \cdot \frac{\sin(3nx)}{3n}
\end{equation}
Their graphs are shown in Figure \ref{fig:f3g3x}.
Figure \ref{fig:f3g3xCloseup} is a closeup of the graphs near $x = 1$, this time, summing the first 200 terms of each series.
The graphs seem to match from about $x = 1$ to about $x = 1.1$.

\begin{figure}[ht]
  \mbox{\includegraphics[width=\picSingleWidth]{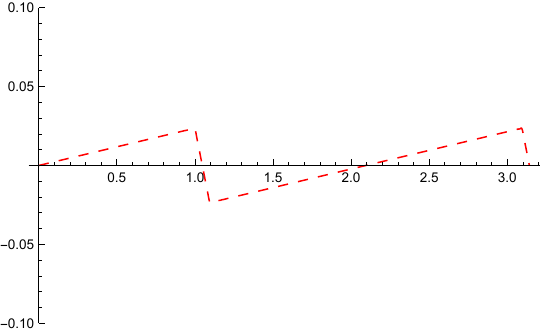}}
  \caption{Four linear pieces of $g(x)$ for $0 \leq x \leq \pi$}
  \label{fig:g3Plot}
\end{figure}

First, let's look at $f(x)$. $\sin(k \pi) = 0$ when $k$ is an integer.
Because $f(x)$ involves $\sin(3nx)$, we see that $f(x)$ must equal 0 when $3nx$ is a multiple of $\pi$.
This happens when $x$ is a multiple of $\pi/3$.
Adding 1 million terms of the series for $f(x)$ shows that $f(.001)$ is be near $\pi/6$.
Also, $f(x)$ seems to decrease linearly to $-\pi/6$ at $x = 2\pi/3$.
So, assuming $f(x)$ consists of linear pieces, we conjecture that
\begin{equation} \label{E:fx3Guess}
f(x) =
  \begin{cases}
    (\pi/3 - x)/2  &\text{for $0 < x < 2\pi/3$,}\\
    (\pi - x)/2    &\text{for $2\pi/3 < x \leq \pi$.}
  \end{cases}
\end{equation}
It is straightforward to verify that the Fourier sine series of this function has coefficients $1/(3n)$, that is, $f(x)$ is, indeed, represented by the series in \eqref{E:f3x}. Also, $f(1) = (\pi/3 - 1)/2 = (\pi - 3)/6$.

Next, we'll examine $g(x)$. This is a little trickier.
All of the following statements are conjectures based on figures \ref{fig:f3g3x}, \ref{fig:f3g3xCloseup}, and \ref{fig:g3Plot}.
Assume that from $x = 0$ to $x = \pi$, $g(x)$ consists of four linear segments, as suggested by Figure \ref{fig:g3Plot}.
After we obtain equations for these four segments, we will prove that the Fourier sine series is the one claimed in Equation \eqref{E:g3x}.

Segment 1:
$g(x)$ increases from $x = 0$, $y = 0$ to $x = 1$, $y = (\pi - 3)/6$.
From Equations \eqref{E:slope} and \eqref{E:intercept}, we get $m = (\pi - 3)/6$ and $b = 0$.
So, the equation for this segment is $y = x(\pi - 3)/6$.

Segment 2:
Over this segment, $g(x) = f(x)$. $g(1) = f(1) = (\pi - 3)/6$. The segment ends at $x = c_1$, where we have to find $c_1$.
Recall that $f(\pi/3) = g(\pi/3) = 0$.
By symmetry (see Figure \ref{fig:f3g3xCloseup}), $c_1$ is as far beyond $\pi/3$ as $\pi/3$ is beyond 1.
That is, $c_1 - \pi/3 = \pi/3 - 1$, so $c_1 = 2\pi/3 - 1 \approx 1.094$.
Also, $g(c_1) = -g(1) = -(\pi - 3)/6$.
So, the equation for this segment is $y = -x/2 + \pi/6$.

Segment 3:
Here, $g(x)$ increases from $x = c_1$ to $x = c_2$.
$g(c_1) = -(\pi - 3)/6$.
Again, by symmetry, $\pi - c_2 = c_1 - \pi/3$, so $c_2 = \pi +  \pi/3 - c_1 = 2\pi/3 + 1$.
$g(c_2)$ is again $(\pi - 3)/6$.
So, the equation for this segment is $y = x (\pi - 3)/6 - \pi(\pi - 3)/9$.

Segment 4:
$g(c_2) = (\pi - 3)/6$ and $g(\pi) = 0$, so the equation for this segment is $y = (\pi - x)/2$.
Summarizing, we conjecture that:
\begin{equation} \label{E:gx3Guess}
g(x) =
  \begin{cases}
    x(\pi - 3)/6     &\text{for $0 \leq x < 1$,}\\
    -x/2 + \pi/6     &\text{for $1 \leq x \leq 2\pi/3 - 1$,}\\
    x (\pi - 3)/6 - (\pi - 3)\pi/9  &\text{for $2\pi/3 - 1 < x \leq 2\pi/3 + 1$,}\\
    (-x + \pi)/2     &\text{for $2\pi/3 + 1 < x \leq \pi$.}
  \end{cases}
\end{equation}
The Fourier sine coefficients for this $g(x)$ can be calculated in \textit{Mathematica} with this code:
\begin{verbatim}
i1 = Integrate[ Sin[n x] * x * (Pi - 3)/6, {x, 0, 1} ]
i2 = Integrate[ Sin[n x] * (-x/2 + Pi/6) , {x, 1, 2 Pi/3 - 1} ]
i3 = Integrate[ Sin[n x] * (x * (Pi - 3)/6 - Pi (Pi - 3)/9), 
                               {x, 2 Pi/3 - 1, 2*Pi/3 + 1} ]
i4 = Integrate[ Sin[n x] * (Pi - x)/2, {x, 2*Pi/3 + 1, Pi} ]
totalgx = Simplify[ (2/Pi) * (i1 + i2 + i3 + i4) ]
\end{verbatim}

Below, we show some of the work:
\begin{align}\label{E:fourIntegrals}
 \frac{2}{\pi}
 \Bigl(
&
  \int_{0}^{1}           \frac{x(\pi - 3)}{6}  \sin(nx) \, dx +
  \int_{1}^{2\pi/3 - 1}  \left( \frac{-x}{2} + \frac{\pi}{6}  \right)  \sin(nx) \, dx + \\
&
  \int_{2\pi/3 - 1}^{2\pi/3 + 1}  \left(  \frac{x (\pi - 3)}{6} - \frac{\pi(\pi - 3)}{9}  \right)  \sin(nx) \, dx +
  \int_{2\pi/3 + 1}^{\pi}         \frac{\pi - x}{2}  \sin(nx) \, dx
  \Bigr)
\notag \\
  &
 =  \frac{2}{\pi} \cdot
 \frac{
     \left(1 + 2 \cos \left(\frac{2 n \pi}{3}\right) \right)
     \left (\pi  \sin(n)-3 \sin \left(\frac{n \pi}{3}\right) \right) }{6 n^2}
\notag \\
  &
 = \frac{
     \left(1 + 2 \cos \left(\frac{2 n \pi}{3}\right) \right)
     \left (\pi  \sin(n)-3 \sin \left(\frac{n \pi}{3}\right) \right) }{3 \pi  n^2}
\notag \\
  &
 = \frac{\sin(3n)}{(3n)^{2}} \, . \notag
\end{align}
Let's examine the fourth line.
If $n$ is \emph{not} a multiple of 3, then $\cos(2 n \pi/3) = -1/2$, so $1 + 2 \cos(2 n \pi/3)$ is 0, amd the entire expression is 0.

If $n$ \emph{is} a a multiple of 3, then $1 + 2 \cos(2 n \pi/3) = 3$.
And if $n$ is a multiple of 3, then $\sin(n \pi/3) = 0$, so the expression becomes $3 \pi \sin(n)/(3 \pi n^2) = \sin(n)/n^2$.
But remember that this is for the case when $n$ is a multiple of 3, that is, $n = 3 k$.
Then we can write $\sin(n)/n^2 = \sin(3 k)/(3 k)^2$.
It is merely a switch of variables to write this as $\sin(3 n)/(3 n)^2$.
This explains how the fourth line in \eqref{E:fourIntegrals} simplifies to the fifth line.

\textbf{QED.}

This completes the proof of Theorem \ref{T:ThmSum3nx}.
We have also proved that, $1 \leq x \leq 2\pi/3 - 1 \approx 1.094$,
\[
\sum_{n=1}^{\infty} \frac{\sin(3nx)}{3n} = \sum_{n=1}^{\infty}  \frac{\sin(3n)}{3n} \cdot \frac{\sin(3nx)}{3n} \, .
\]

Applying Parseval's equation \eqref{E:ParsevalsEquation} to $g(x)$ (Equation \eqref{E:g3x}) requires evaluating four integrals.
When we have \textit{Mathematica} do the work, we get
\[
\sum_{n=1}^{\infty} \left( \frac{\sin(3n)}{(3n)^2} \right)^2 = \frac{(\pi - 3)^2}{54} \, .
\]

\textbf{A Simplification.}
The above proof was the one originally worked out by the author.
However, we can simplify the calculation of the Fourier coefficients of $g(x)$.
Instead of integrating from $x = 0$ to $x = \pi$, we need to integrate only up to $x = \pi/3$.
Here's why.

$g(x)$ as defined by Equation \eqref{E:g3x} has period $p = 2\pi/3$ .
We can see this if we graph two periods of $g(x)$, from $x = 0$ to $x = 2 \cdot 2\pi/3$ .
\begin{figure}[ht]
  \mbox{\includegraphics[width=\picSingleWidth]{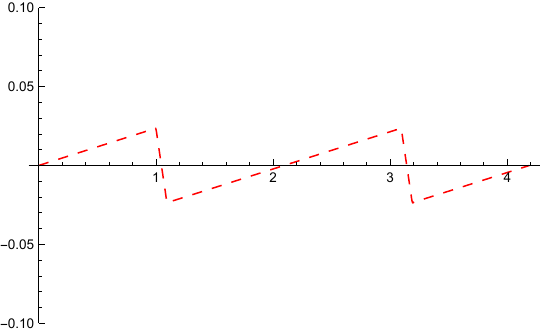}}
  \caption{Two periods of $g(x)$ (one period = $2 \pi/3$)}
  \label{fig:g3TwoPeriods}
\end{figure}
See Figure \ref{fig:g3TwoPeriods}.
Also, $g(x)$ is an odd function, because only the $\sin$ terms are present in the Fourier series (i.e., the $\cos$ terms are all 0).

Therefore, we can use Equation \eqref{E:bnOddPeriodP} to integrate over the shorter interval
$x = 0$ to $x = b = p/2 = \pi/3$, instead of from $x = 0$ to $x = \pi$ as we did above.
This means only the first \emph{two} pieces of $g(x)$ from Equation \eqref{E:gx3Guess} will be needed in the integrals.
Only two pieces are needed because $\pi/3 \approx 1.047$ is less than the righthand endpoint of the second piece, $2\pi/3 - 1 \approx 1.094$.
We have $p = 2 b = 2\pi/3$, so $b = p/2 = \pi/3$, $2/b = 6/\pi$, and $n \cdot 2\pi/p \cdot x = 3nx$.
So, applying Equation \eqref{E:bnOddPeriodP}, we get
\begin{align*}
b_n
  & = \frac{2}{b} \int_{0}^{b} F(x) \sin(n \cdot \frac{2 \pi}{p} \cdot x) \, dx \\
  & = \frac{6}{\pi} \int_{0}^{\pi/3} g(x) \cdot \sin(3nx) \, dx \\
  & = \frac{6}{\pi}
      \left(
        \int_{0}^{1}           \frac{x(\pi - 3)}{6}  \sin(3nx) \, dx +
        \int_{1}^{\pi/3} \left( \frac{-x}{2} + \frac{\pi}{6}  \right)  \sin(3nx) \, dx
      \right) \\
  & = \frac{6}{\pi}
      \left(
        -\frac{(\pi -3) (3 n \cos (3 n)-\sin (3 n))}{54 n^2} +
        \frac{\sin (3 n)+(\pi -3) n \cos (3 n)}{18 n^2}
      \right) \\
  & = \frac{\sin(3 n)}{(3 n)^2}
\end{align*}

Here is the \textit{Mathematica} code to do this simplified calculation:
\begin{verbatim}
  i1 = Integrate[ Sin[3 n x] * x * (Pi - 3)/6, {x, 0, 1} ]
  i2 = Integrate[ Sin[3 n x] * (-x/2 + Pi/6) , {x, 1, Pi/3} ]
  total = Simplify[ (6/Pi) * (i1 + i2), Assumptions -> Element[n, Integers] ]
\end{verbatim}

\ifthenelse {\boolean{BKMRK}}
  { \section{Sum of \texorpdfstring{$\sin(k n x)/(k n)$}{sin(k n x)/(k n)}}\label{S:sinkn} }
  { \section{Sum of $\sin(k n x)/(k n)$}\label{S:sinkn} }

We have seen that, for $k = 1$, $k = 2$, and $k = 3$, we have
\begin{equation}\label{E:sin123x}
\sum_{n=1}^{\infty} \frac{\sin(kn)}{kn} = \sum_{n=1}^{\infty} \left( \frac{\sin(kn)}{kn} \right)^2 = \frac{\pi - k}{2k} \, .
\end{equation}

Are these just three special cases of a more general pattern?
Does this identity hold for non-integer values of $k$ between 1 and 3?
Does it hold for $k < 1$ or $k > 3$?

If we numerically add, say, one million terms of these series for different values of $k$, it \emph{appears} that Equation \eqref{E:sin123x} is true, provided $0 < k < \pi$.
These sums are shown in Table \ref{Ta:Sum1Sum2Table}.

\begin{table}[ht]
 \begin{center}
  \begin{tabular}{ r r r r }
   $k$  &  $ \sum \sin(kn)/(kn)$  & $\sum (\sin(kn)/(kn))^2$  &  $(\pi-k)/(2k)$   \\ \hline
   1/2    & 2.64159 66853  &  2.64159 06536  &  2.64159 26536  \\
   3/4    & 1.59439 52964  &  1.59439 42135  &  1.59439 51024  \\
   1      & 1.07079 52944  &  1.07079 58268  &  1.07079 63268  \\
	 3/2    & 0.54719 80291  &  0.54719 73290  &  0.54719 75512  \\
	 5/3    & 0.44247 80756  &  0.44247 76161  &  0.44247 77961  \\
	 2      & 0.28539 78783  &  0.28539 80384  &  0.28539 81634  \\
	 5/2    & 0.12831 87282  &  0.12831 84507  &  0.12831 85307  \\
	 3      & 0.02359 86235  &  0.02359 87200  &  0.02359 87756  \\
	 31/10  & 0.00670 86384  &  0.00670 84405  &  0.00670 84925  \\
	 $\pi$  & 0              &  0              &  0   \\
	 32/10  & $-0$.00912 62882 &  0.00879 29527  &  $-0$.00912 61479  \\
	 7/2    & $-0$.05120 09220 &  0.04071 48188  &  $-0$.05120 10495  \\
  \end{tabular}
  \caption{Sums of $\sin(kn)/(kn)$ and $(\sin(kn)/(kn))^2$}
  \label{Ta:Sum1Sum2Table}
 \end{center}
\end{table}

Note that Equation \eqref{E:sin123x} can't possibly hold if $k > \pi$ because in this case, $(\pi-k)/(2k) < 0$,
whereas the second sum in Equation \eqref{E:sin123x} must be $\geq 0$.

So, let's compare some graphs of
\begin{equation}\label{E:fkDef1}
f_k(x) = \sum_{n=1}^{\infty} \frac{\sin(knx)}{kn}
\end{equation}
and
\begin{equation}\label{E:gkDef1}
g_k(x) = \sum_{n=1}^{\infty} \frac{\sin(kn)}{kn} \cdot \frac{\sin(knx)}{kn}
\end{equation}
for various values of $k$.

\begin{figure}[ht] 
\centering
\begin{minipage}[t]{\figureMinipageWidth}
\centering
\includegraphics[width=\picDblWidth]{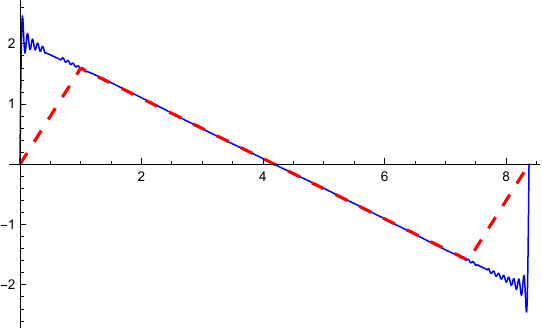}
\caption{Sums with $k = 3/4$, up to $x = 2\pi/k = 8\pi/3$}
\label{fig:plotKa1p}
\end{minipage}\hfill
\begin{minipage}[t]{\figureMinipageWidth}
\centering
\includegraphics[width=\picDblWidth]{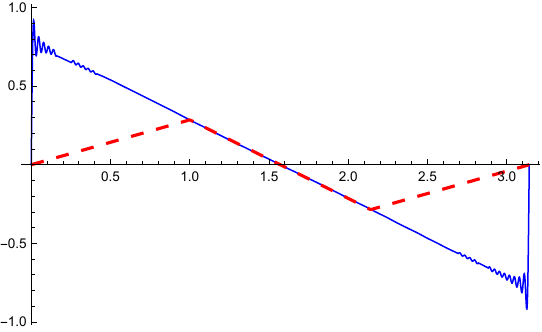}
\caption{Sums with $k = 2$, up to $x = 2\pi/k = \pi$}
\label{fig:plotKb1p}
\end{minipage}
\end{figure}

\begin{figure}[ht]
\centering
\begin{minipage}[t]{\figureMinipageWidth}
\centering
\includegraphics[width=\picDblWidth]{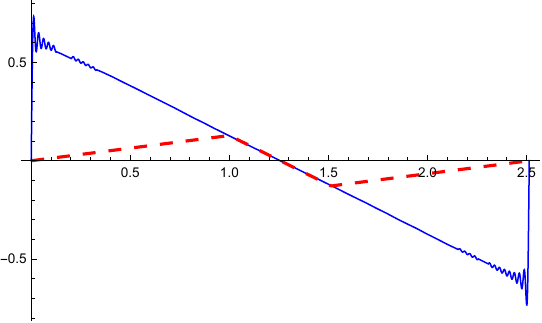}
\caption{Sums with $k = 5/2$, up to $x = 2\pi/k = 4\pi/5$}
\label{fig:plotKc1p}
\end{minipage}\hfill
\begin{minipage}[t]{\figureMinipageWidth}
\centering
\includegraphics[width=\picDblWidth]{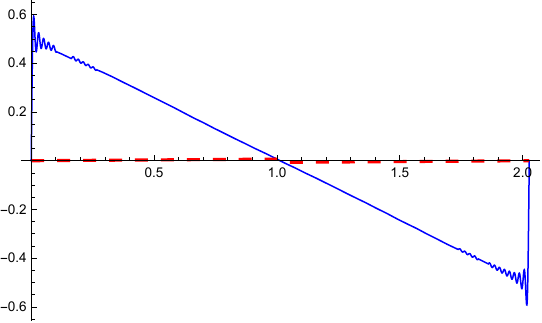}
\caption{Sums with $k = 31/10$, up to $x = 2\pi/k = 20\pi/31$}  
\label{fig:plotKd1p}
\end{minipage}
\end{figure}

\begin{figure}[ht]
\centering
\begin{minipage}[t]{\figureMinipageWidth}
\centering
\includegraphics[width=\picDblWidth]{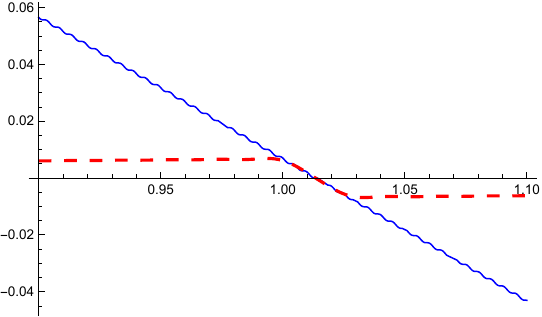}
\caption{Closeup of Figure \ref{fig:plotKd1p}}

\label{fig:plotKd1p2}
\end{minipage}\hfill
\begin{minipage}[t]{\figureMinipageWidth}
\centering
\includegraphics[width=\picDblWidth]{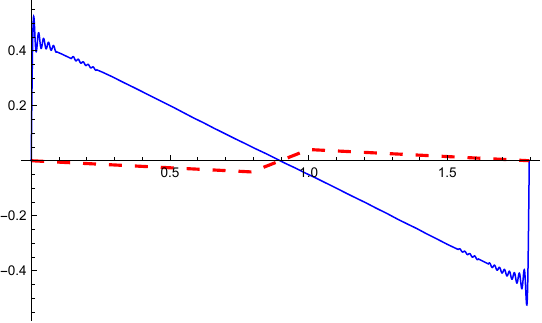}
\caption{Sums with $k = 7/2$, up to $x = 2\pi/k = 4\pi/7$}
\label{fig:plotKe1p}
\end{minipage}
\end{figure}

Note that both $f_k(x)$ and $g_k(x)$ have period $2\pi/k$.
With this in mind, let's plot exactly one full period of $f_k(x)$ and $g_k(x)$ with different values of $k$.
In Figures \ref{fig:plotKa1p} - \ref{fig:plotKe1p}, we plot the graphs for $k = 3/4$, $k = 2$, $k = 5/2$, $k = 31/10$, and $k = 7/2$.
In each graph, the solid curve is $f_k(x)$, and the dashed curve is $g_k(x)$.
Each graph shows the sum of 100 terms of each series.
Figure \ref{fig:plotKd1p2} is a closeup view of the graph for $k = 31/10$, using 400 terms of the series.

For $k = 7/2 = 3.5 > \pi$, (Figure \ref{fig:plotKe1p}), the graphs seem \emph{not} to agree except for those $x$ for which $f_k(x) = g_k(x) = 0$.
Otherwise, for each value of $k < \pi$, it appears that there is an \emph{interval} from $x = 1$ to $x = 2\pi/k - 1$ over which $f_k(x) = g_k(x)$.
Let's see if we can prove this.

First, each $f_k(x)$ appears to consist of one segment, from $f_k(0) = \pi/(2k)$ to $f_k(2\pi/k) = -\pi/(2k)$.
That is, it \emph{seems} that, for $0 < x < 2\pi/k$, we have
\begin{equation}\label{E:fkDef2}
f_k(x) = -\frac{x}{2} + \frac{\pi}{2k} = \frac{\pi - kx}{2 k} \, .
\end{equation}
The coefficient of $\sin(nx)$ in the Fourier sine series for this $f_k(x)$ can be computed using
Equation \eqref{E:bnOddPeriodP}, with $p = 2 b = 2 \pi/k$.
In this case, $2/b = 2k/\pi$ and $2 \pi/p = k$.
So, Equation \eqref{E:bnOddPeriodP} becomes
\[
\frac{2}{b} \int_{0}^{b} F(x) \sin \left( n \cdot \frac{2 \pi}{p} \cdot x \right) \, dx \, =
\frac{2 k}{\pi} \int_{0}^{\pi/k} \left( \frac{\pi - k x}{2k} \right) \sin(n \cdot kx) \, dx = \frac{1}{kn} \, .
\]
This shows that the $f_k(x)$ given by Equation \eqref{E:fkDef2} \emph{does} have the sine series given by Equation \eqref{E:fkDef1}.

$g_k(x)$ given in Equation \eqref{E:gkDef1} is more complicated.
From the graphs, it looks like $g_k(x)$ consists of three linear pieces.
By closely examining the graphs, or by calculating $g_k(x)$ for a variety of values of $k$ and $x$, the reader can verify that the three pieces \emph{appear} to be:
\begin{equation}\label{E:gkDef2}
g_k(x)=
  \begin{dcases}
    \frac{x(\pi - k)}{2k}  &\text{for $0 \leq x < 1$,}\\
    \frac{\pi - k x}{2k}  &\text{for $1 \leq x \leq \frac{2\pi}{k} - 1$,}\\
    \frac{x(\pi - k)}{2k} - \frac{\pi(\pi - k)}{k^2}  &\text{for $\frac{2\pi}{k} - 1 < x \leq \frac{2\pi}{k}$.}
  \end{dcases}
\end{equation}

Given the value of $k$, $g_k(x)$ is an odd function having period $p = 2\pi/k$.
Again, we will use Equation \eqref{E:bnOddPeriodP}, with $p = 2 b = 2 \pi/k$.
As above, $2/b = 2k/\pi$ and $2 \pi/p = k$.
Equation \eqref{E:bnOddPeriodP} becomes
\[
\frac{2}{b} \int_{0}^{b} F(x) \sin \left( n \cdot \frac{2 \pi}{p} \cdot x \right) \, dx \, =
\frac{2 k}{\pi} \int_{0}^{\pi/k} g_k(x) \sin(n \cdot kx) \, dx \, .
\]

Because we integrate only up to $x = p/2 = \pi/k < 2 \pi/k - 1$, we don't need to integrate the third linear expression that occurs in equation \eqref{E:gkDef2}.

The $n$\textsuperscript{th} coefficient in the sine series for the $g_k(x)$ in \eqref{E:gkDef2} is therefore given by
\begin{align}\label{E:gkTwoIntegrals}
&
 \frac{2 k}{\pi}
 \left(
  \int_{0}^{1}           \frac{x(\pi - k)}{2k}  \sin(knx) \, dx +
  \int_{1}^{\pi/k}  \left( \frac{\pi - k x}{2 k} \right)  \sin(knx) \, dx   \right)
\notag \\
 = & \frac{ 2k}{\pi}
  \left(
  \frac{(k-\pi ) (k n \cos (k n)-\sin(k n))}{2 k^3 n^2} + 
  \frac{ n (-k + \pi) \cos (k n)+ \sin (k n) - \sin (n \pi)}{2 k^2 n^2}
   \notag
  \right)
\notag \\
 = & \frac{ 2k}{\pi}
  \left(
  \frac{(k - \pi ) (k n \cos (k n) - \sin(k n))}{2 k^3 n^2} + 
  \frac{(k - \pi) (-kn) \cos (k n) +  k \sin (k n) }{2 k^3 n^2}
   \notag
  \right)
\notag \\
 = & \frac{ 2k}{\pi}
  \left(
  \frac{(k - \pi )(- \sin(k n))}{2 k^3 n^2} + 
  \frac{k \sin (k n) }{2 k^3 n^2}
   \notag
  \right)
\notag \\
 = & \frac{\sin(kn)}{(kn)^{2}} \, . \notag
\end{align}
Here is the \textit{Mathematica} code to do this integration:
\begin{verbatim}
  j1 = Integrate[Sin[k n x]*x (Pi - k)/(2 k), {x, 0, 1}]
  j2 = Integrate[Sin[k n x]*((Pi - k x)/(2 k)), {x, 1, Pi/k}]
  totalGk = Simplify[(2 k/Pi) (j1 + j2), Assumptions -> Element[n, Integers]]
\end{verbatim}


This proves that $g_k(x)$ defined by Equation \eqref{E:gkDef2} does have the sine series given by Equation \eqref{E:gkDef1}.

(One would get the same result by integrating from 0 to $2\pi/k$ using all three parts of \eqref{E:gkDef2}.
In this case, the constant in front of the integral sign would be $k/\pi$, not $2k/\pi$).

So, using Equations \eqref{E:fkDef1}, \eqref{E:gkDef1}, \eqref{E:fkDef2}, and \eqref{E:gkDef2},
we have proved that, for all k in $0 < k < \pi$ and all $x$ in $1 \leq x \leq 2\pi/k - 1$,
\begin{equation}\label{E:twoGeneralSumsWithSinkx}
\sum_{n=1}^{\infty} \frac{\sin(knx)}{kn} = \sum_{n=1}^{\infty} \frac{\sin(kn)}{kn} \cdot \frac{\sin(knx)}{kn} = \frac{\pi - kx}{2k} \, .
\end{equation}
(This holds trivially for $k = \pi$ because then, with $x = 1$, all terms are 0).

Substituting $x = 1$ into Equation \eqref{E:twoGeneralSumsWithSinkx} proves that, for $0 < k \leq \pi$,
\begin{equation}\label{E:twoSumsWithSinkx}
\sum_{n=1}^{\infty} \frac{\sin(kn)}{kn} = \sum_{n=1}^{\infty} \left( \frac{\sin(kn)}{kn} \right)^2 = \frac{\pi - k}{2k} \, .
\end{equation}

Applying Parseval's Equation \eqref{E:ParsevalsEquation} to $g_k(x)$ in \eqref{E:gkDef2}, we get (with \textit{Mathematica}'s help):
\begin{align*}
& \sum_{n=1}^{\infty} \left( \frac{\sin(kn)}{(kn)^2} \right)^2
 = \frac{k}{\pi} \int_{0}^{2\pi/k} \left( g_k(x) \right) ^2 \, dx
  \\
 = \frac{k}{\pi}
 \Bigl(
&
  \int_{0}^{1}           \left( \frac{x(\pi - k)}{2k} \right) ^2  \, dx +
  \int_{1}^{2\pi/k - 1}  \left( \frac{\pi - k x}{2k}  \right) ^2  \, dx +
   \\
&
  \int_{2\pi/k - 1}^{2\pi/k}  \left(  \frac{x (\pi - k)}{2k} - \frac{\pi(\pi - k)}{k^2}  \right) ^2 \, dx
  \Bigr)
 \\
 = \frac{k}{\pi}
&
  \left(
  \left( \frac{1}{12} - \frac{\pi}{6k} + \frac{\pi^2}{12k^2} \right) - \frac{(\pi-k)^3}{6k^3} +
  \left( \frac{1}{12} - \frac{\pi}{6k} + \frac{\pi^2}{12k^2} \right)
  \right) \, .
\notag
\end{align*}
Simplifying, we get the following, which holds for $0 < k \leq \pi$:
\begin{equation}\label{E:gkParseval}
\sum_{n=1}^{\infty} \left( \frac{\sin(kn)}{(kn)^2} \right)^2 = \frac{(\pi-k)^2}{6k^2} \, .
\end{equation}

Here is the \textit{Mathematica} code for these details:
\begin{verbatim}
  i1 = Integrate[ (x (Pi - k)/(2k))^2, {x, 0, 1}]
  i2 = Integrate[ ((Pi - k x)/(2 k))^2, {x, 1, 2 Pi/k - 1}]
  i3 = Integrate[ (x(Pi - k)/(2k) - Pi(Pi - k)/k^2)^2, {x, 2 Pi/k - 1, 2 Pi/k}]
  parsevalGk = Simplify[ (k/Pi)(i1 + i2 + i3) ]
\end{verbatim}

\textbf{Discussion.}
The results in this section generalize some of the results in Sections \ref{S:InterestingPiFormulas}, \ref{S:sin2n}, and \ref{S:sin3n}.
The reader may therefore wonder why those sections are present in this paper.

In many math papers, the theorems and proofs are presented without any hint as to how they were discovered.
One purpose of this paper is to show \emph{how} the discovery process gave rise to the results presented here.
The results in the last four sections were presented in roughly the order in which the author discovered them.
Moreover, the above discussions attempt to illustrate \emph{how} the results were discovered and proved.

\section{Series with alternating signs} \label{S:AlternatingSigns}

Now let's obtain a few series whose terms have alternating signs.
We will use Equations \eqref{E:five1}, \eqref{E:five2}, and \eqref{E:evenTerms}, which we restate here for convenience:
\begin{align*}
  \sum_{n=1}^{\infty} \frac{\sin(nx)}{n} & = \frac{\pi - x}{2}  \quad \text{ for $0 < x < 2 \pi$} & \text{Equation \eqref{E:five1},} \\
  \sum_{n=1}^{\infty} \frac{\sin(n)}{n} \cdot \frac{\sin(nx)}{n} & = \frac{\pi - x}{2}  \quad \text{ for $1 \leq x \leq 2\pi - 1$} & \text{Equation \eqref{E:five2},} \\
  \sum_{n=1}^{\infty} \frac{\sin(2nx)}{2n} = \sum_{n=1}^{\infty} \frac{\sin(2n)}{2n} \frac{\sin(2nx)}{2n} & = \frac{\pi - 2x}{4} \quad \text{ for $1 \leq x \leq \pi - 1$} & \text{Equation \eqref{E:evenTerms}.}
\end{align*}
Let's subtract Equation \eqref{E:five1} minus two times Equation \eqref{E:evenTerms}.
First, subtracting the right-hand sides, we have
\[
\frac{\pi - x}{2} - 2 \cdot \frac{\pi - 2x}{4} = \frac{x}{2} \, .
\]
Next, subtract the left-hand sides.
The sum in Equation \eqref{E:five1}, minus two times the \emph{first} sum in Equation \eqref{E:evenTerms} is
\begin{align*}
\frac{x}{2} = & \frac{\sin(x)}{1} + \frac{\sin(2x)}{2} + \frac{\sin(3x)}{3} + \frac{\sin(4x)}{4} \dots \\
 - 2 & \left( \frac{\sin(2x)}{2} + \frac{\sin(4x)}{4} + \frac{\sin(6x)}{6} + \frac{\sin(8x)}{8} \dots \right) \\
= & \frac{\sin(x)}{1} - \frac{\sin(2x)}{2} + \frac{\sin(3x)}{3} - \frac{\sin(4x)}{4} \dots \\
= & \sum_{n=1}^{\infty} (-1)^{n+1} \frac{\sin(nx)}{n} \, .
\end{align*}
Similarly, if we subtract the sum in Equation \eqref{E:five2}, minus two times the \emph{second} sum in Equation \eqref{E:evenTerms}, we get
\[
\frac{x}{2} = \sum_{n=1}^{\infty} (-1)^{n+1} \frac{\sin(n)}{n} \frac{\sin(nx)}{n} \, .
\]
These equations are true for those $x$ for which all of the original sums are true, that is, for $1 \leq x \leq \pi - 1$.
We can also prove that these sums hold for a larger interval, namely, for $0 \leq x \leq \pi - 1$:

\begin{theorem}\label{T:TwoAlternatingSums}
For $0 \leq x \leq \pi - 1$,
\begin{equation} \label{E:xOver2ab}
\sum_{n=1}^{\infty} (-1)^{n+1} \frac{\sin(nx)}{n} = 
\sum_{n=1}^{\infty} (-1)^{n+1} \frac{\sin(n)}{n} \cdot \frac{\sin(nx)}{n} = \frac{x}{2} \, .
\end{equation}
The first sum equals $x/2$ over the larger interval, $0 \leq x < \pi$.
\end{theorem}

\textbf{Proof.}

\begin{figure}[ht]
  \mbox{\includegraphics[width=\picSingleWidth]{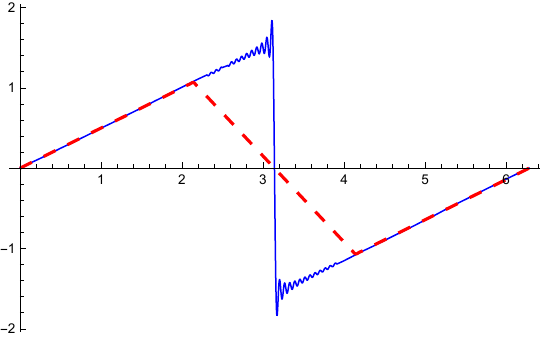}}
  \caption{$\sum_{n=1}^{100} (-1)^{n+1} \frac{\sin(nx)}{n}$ (solid) and $\sum_{n=1}^{100} (-1)^{n+1} \frac{\sin(n)}{n} \cdot \frac{\sin(nx)}{n}$ (dashed)}
  \label{fig:xOver2ab}
\end{figure}

Figure \ref{fig:xOver2ab} shows the graphs of
\begin{equation}  \label{E:altSignF1}
f(x) = \sum_{n=1}^{\infty} (-1)^{n+1} \frac{\sin(nx)}{n}
\end{equation}
and
\begin{equation}  \label{E:altSignG1}
g(x) = \sum_{n=1}^{\infty} (-1)^{n+1} \frac{\sin(n)}{n} \cdot \frac{\sin(nx)}{n} \, .
\end{equation}

First, the fact that the function $f(x) = x/2$ (for $-\pi < x < \pi$) has the Fourier series given in Equation \eqref{E:altSignF1} is a standard result.
See, for example, \cite[Equation (a) on p.\ 506]{Kaplan} or \cite[Equation (13.9) on p.\ 32]{Tolstov}.
It's easy enough to do by hand, but here is the \textit{Mathematica} code for this calculation:
\begin{verbatim}
  i1 = Integrate[ Sin[n x] * (x/2), {x, 0, Pi}]
  Simplify[ (2/Pi) * i1, Assumptions -> Element[n, Integers] ]
\end{verbatim}

Figure \ref{fig:xOver2ab} suggests that, for $\pi - 1 < x \leq \pi + 1$, $g(x)$ is
\[
  \frac{\pi - 1}{2}(\pi - x) = \frac{x}{2} + \frac{\pi}{2}(\pi - 1 - x) \, .
\]
We will now prove that over the interval $0 \leq x \leq \pi + 1$, $g(x)$ is given by
\begin{equation} \label{E:altSignG2}
g(x) = \sum_{n=1}^{\infty} (-1)^{n+1} \frac{\sin(n)}{n} \cdot \frac{\sin(nx)}{n} =
  \begin{dcases}  
    \frac{x}{2}  &\text{for $0 \leq x \leq \pi-1$,}\\
    \frac{x}{2} + \frac{\pi}{2}(\pi - 1 - x) &\text{for $\pi - 1 \leq x \leq \pi + 1$.}
  \end{dcases}
\end{equation}


Frst, this \textit{Mathematica} code calculates the Fourier coefficients for the $g(x)$ in Equation \eqref{E:altSignG2},
over the more convenient interval $0 \leq x \leq \pi$:
\begin{verbatim}
  i1 = Integrate[ Sin[n x] * (x/2), {x, 0, Pi - 1}]
  i2 = Integrate[ Sin[n x] * (x/2 + (Pi/2)*(Pi - 1 - x)), {x, Pi - 1, Pi}]
  trigRules = { Sin[n*(c_ + Pi)] -> Cos[n Pi] * Sin[n c], Cos[n Pi] -> (-1)^n }
  total = Simplify[ (2/Pi) * (i1 + i2),
                    Assumptions -> Element[n, Integers] ] //. trigRules
\end{verbatim}
The output of this calculation is
\[
-\frac{ (-1)^n \sin(n)}{n^2} = \frac{(-1)^{n+1} \sin(n)}{n^2} \, ,
\]
so that Equation \eqref{E:altSignG2} holds at least up to $x = \pi$.

We will now show that the second expression in Equation \eqref{E:altSignG2} actually holds up to $x = \pi + 1$.
If we replace $x$ by $2 \pi - x$ in Equation \eqref{E:altSignG2}, the left side merely changes sign, because $\sin(x) = -\sin(2\pi-x)$.
The right side also changes sign:
\[
\frac{2 \pi - x}{2} + \frac{\pi}{2}(\pi - 1 - (2 \pi - x))
= - \left( \frac{x}{2} + \frac{\pi}{2}(\pi - 1 - x) \right) \, .
\]

Therefore, Equation \eqref{E:altSignG2} is actually valid for $0 \leq x \leq \pi + 1$.

\textbf{QED.}

Note: In Section \ref{S:Dilogarithms}, we will show how to \emph{derive} the polynomials in \eqref{E:altSignG2} from the series for $g(x)$.

If we substitute $x = 1$ into Equation \eqref{E:xOver2ab}, we get
\begin{equation} \label{E:alternatingSigns}
\sum_{n=1}^{\infty} (-1)^{n+1} \frac{\sin(n)}{n} = 
\sum_{n=1}^{\infty} (-1)^{n+1} \left( \frac{\sin(n)}{n} \right) ^2 = \frac{1}{2} \, ,
\end{equation}
which is a nice variation of \eqref{E:Prob6241a}, but with alternating signs.

One is naturally led to wonder whether Equation \eqref{E:xOver2ab} or \eqref{E:alternatingSigns} holds with higher powers of $\sin(n)/n$.
If we perform a small numeric experiment, we make the interesting discovery that the sum of the first million terms of this series
\begin{equation*}
  \sum_{n=1}^{\infty} (-1)^{n+1} \left( \frac{\sin(n)}{n} \right) ^3
\end{equation*}
is about $0.50000 \text{ } 00000 \text{ } 00000 \text{ } 00067 \, .$
\textit{Mathematica} claims that the sum of the infinite series is \emph{exactly} $1/2$:
\begin{verbatim}
  Sum[(-1)^(n + 1) * (Sin[n]/n)^3, {n, 1, Infinity}] // FullSimplify
\end{verbatim}
In Sections \ref{S:PolylogApplicationX} through \ref{S:AltSum4xResult2}, we will prove this and much more.
We will see that Equation \eqref{E:xOver2ab} holds with one or two more powers of $\sin(n)/n$, although for shorter ranges of $x$ values.
See Equation \eqref{E:ThreeAlternatingSincSums} and Theorem \ref{T:AltSumkxTheorem} for a collection of related results.


\ifthenelse {\boolean{BKMRK}}
  { \section{Multiplying by powers of \texorpdfstring { $\sin(n)/n$}{sin(n)/n} } \label{S:sinTokth} }
  { \section{Multiplying by powers of $(\sin(n)/n)$}\label{S:sinTokth} }

We now discuss another way to extend Equations \eqref{E:five1} and \eqref{E:five2}, which are restated here for convenience.
Those equations say that, for $1 \leq x \leq 2\pi - 1$, we have
\begin{equation}\label{E:TwoSumsEqual}
  \sum_{n=1}^{\infty} \frac{\sin(nx)}{n} =
  \sum_{n=1}^{\infty} \frac{\sin(n)}{n} \frac{\sin(nx)}{n} =
  \frac{\pi - x}{2}  \, .
\end{equation}
We were able to multiply the $n$\textsuperscript{th} term in the series on the left by $\sin(n)/n$ without changing the sum.
Can we introduce even more factors of $\sin(n)/n$ into the $n$\textsuperscript{th} term without changing the sum?
Numerical experiments with Equation \eqref{E:TwoSumsEqual} suggested that sometimes we can.

With $x = 2$, Equation \eqref{E:TwoSumsEqual} guarantees that both of these sums are equal to $(\pi-2)/2$, which is about
0.57079 63267 94896 61923:
\[
\sum_{n=1}^{\infty} \frac{\sin(2n)}{n} = \sum_{n=1}^{\infty} \frac{\sin(n)}{n} \frac{\sin(2n)}{n} = \frac{\pi - 2}{2} \, .
\]

Let's put in higher powers of $\sin(n)/n$, add a million terms, and see what we get. First,
\[
  \sum_{n=1}^{1000000} \left(\frac{\sin(n)}{n} \right) ^{2} \frac{\sin(2n)}{n}
    \approx 0.57079 \text{ } 63267 \text{ } 94896 \text{ } 61906 \, .
\]
That sum is pretty close to $(\pi - 2)/2$.
However, the following sum is \emph{not} very close to $(\pi - 2)/2$:
\[
  \sum_{n=1}^{1000000} \left(\frac{\sin(n)}{n} \right) ^{3} \frac{\sin(2n)}{n} \approx 0.50534 \text{ } 64798 \, .
\]

So, we conjecture (and will later prove) that all three of these sums equal $(\pi - 2)/2$:
\begin{equation}\label{E:ThreeEqualSumsk2}
\sum_{n=1}^{\infty} \frac{\sin(2n)}{n} =
\sum_{n=1}^{\infty} \frac{\sin(n)}{n} \frac{\sin(2n)}{n} =
\sum_{n=1}^{\infty} \left(\frac{\sin(n)}{n} \right) ^{2} \frac{\sin(2n)}{n} =
\frac{\pi - 2}{2} \, .
\end{equation}

With $x = 4$, Equation \eqref{E:TwoSumsEqual}, tells us that
\[
  \sum_{n=1}^{\infty} \frac{\sin(4n)}{n} =
  \sum_{n=1}^{\infty} \frac{\sin(n)}{n} \frac{\sin(4n)}{n} = \frac{\pi - 4}{2} \, .
\]
Numerical calculation with $x = 4$ suggests that the following sum also equals
 $(\pi - 4)/2 \approx$ -0.42920 36732 05103 38076 :
\[
\sum_{n=1}^{1000000} \left(\frac{\sin(n)}{n} \right) ^{2} \frac{\sin(4n)}{n}
   \approx -0.42920 \text{ } 36732 \text{ } 05103 \text{ } 38027 \, .
\]
However, a numerical calculation shows that when we introduce $(\sin(n)/n)^3$, the result is different:
\[
\sum_{n=1}^{1000000} \left(\frac{\sin(n)}{n} \right) ^{3} \frac{\sin(4n)}{n} \approx -0.40509 \text{ } 74412 \, .
\]

So, we conjecture (and will later prove) that
\begin{equation}\label{E:ThreeEqualSumsk4}
\sum_{n=1}^{\infty} \frac{\sin(4n)}{n} =
\sum_{n=1}^{\infty} \frac{\sin(n)}{n} \frac{\sin(4n)}{n} =
\sum_{n=1}^{\infty} \left(\frac{\sin(n)}{n} \right) ^{2} \frac{\sin(4n)}{n} =
\frac{\pi - 4}{2} \, .
\end{equation}

For $x = 3$, the results are even more interesting.
Equation \eqref{E:TwoSumsEqual} shows that
\[
\sum_{n=1}^{\infty} \frac{\sin(3n)}{n} =
\sum_{n=1}^{\infty} \frac{\sin(n)}{n} \frac{\sin(3n)}{n} =
\frac{\pi - 3}{2} \, .
\]
Numeric calculations suggest that the following two sums \emph{also} equal $(\pi - 3)/2$,
which is about $ \approx 0.07079 \text{ } 63267 \text{ } 94896 \text{ } 61923 \text{ } 132 $:
\[
\sum_{n=1}^{1000000} \left(\frac{\sin(n)}{n} \right) ^{2} \frac{\sin(3n)}{n}
  \approx 0.07079 \text{ } 63267 \text{ } 94896 \text{ } 61941 \text{ } 941
\]
\[
\sum_{n=1}^{1000000} \left(\frac{\sin(n)}{n} \right) ^{3} \frac{\sin(3n)}{n}
  \approx 0.07079 \text{ } 63267 \text{ } 94896 \text{ } 61927 \text{ } 298  \, .
\]

So, we conjecture (and will later prove) that \emph{all four} of the following sums are \emph{exactly} equal to $(\pi - 3)/2$:
\begin{equation}\label{E:FourSumsWithN}
\sum_{n=1}^{\infty} \frac{\sin(3n)}{n} =
\sum_{n=1}^{\infty} \frac{\sin(n)}{n} \frac{\sin(3n)}{n} =
\sum_{n=1}^{\infty} \left(\frac{\sin(n)}{n} \right) ^{2} \frac{\sin(3n)}{n} =
\sum_{n=1}^{\infty} \left(\frac{\sin(n)}{n} \right) ^{3} \frac{\sin(3n)}{n} =
\frac{\pi - 3}{2} \, .
\end{equation}

A numeric calculation shows that, with $(\sin(n)/n)^4$, we have
\[
\sum_{n=1}^{1000000} \left( \frac{\sin(n)}{n} \right) ^{4} \frac{\sin(3n)}{n} \approx 0.06477 \text{ } 50586 \, ,
\]
which is significantly different from $(\pi - 3)/2 \approx 0.07079 \text{ } 63268 $.
So, Equation \eqref{E:FourSumsWithN} can \emph{not} be extended to a sum involving $(\sin(n)/n)^4$.

In fact, graphs and numeric calculations also suggest that, in Equation \eqref{E:TwoSumsEqual}, we can introduce up to three factors of $\sin(n)/n$, not just for $x = 3$, but for some $x$ values a little more than 3, and that, for each of these $x$, all four sums appear to equal $(\pi - x)/2$.

These experiments suggest the following theorem, which extends Equations \eqref{E:ThreeEqualSumsk2}, \eqref{E:ThreeEqualSumsk4}, and \eqref{E:FourSumsWithN} to an \emph{entire interval} of $x$ values.

\begin{theorem}\label{T:ThmFourEqualSums}
For all $x$ in the interval $3 \leq x \leq 2\pi - 3$, we have
\begin{equation} \label{E:FourEqualSums}
\sum_{n=1}^{\infty} \frac{\sin(nx)}{n} =
\sum_{n=1}^{\infty} \frac{\sin(n)}{n} \frac{\sin(nx)}{n} =
\sum_{n=1}^{\infty} \left(\frac{\sin(n)}{n} \right) ^{2} \frac{\sin(nx)}{n} =
\sum_{n=1}^{\infty} \left(\frac{\sin(n)}{n} \right) ^{3} \frac{\sin(nx)}{n} =
\frac{\pi - x}{2} \, .
\end{equation} 
Moreover: \\
The first sum equals $(\pi - x)/2$ for $0 < x < 2\pi \approx 6.28 \, $. \\
The second sum equals $(\pi - x)/2$ for $1 \leq x \leq 2\pi - 1 \approx 5.28 \, $. \\
The third sum equals $(\pi - x)/2$ for $2 \leq x \leq 2\pi - 2 \approx 4.28 \, $. \\
The fourth sum equals $(\pi - x)/2$ for $3 \leq x \leq 2\pi - 3 \approx 3.28 \, $.
\end{theorem}

\begin{figure}[ht]
  \mbox{\includegraphics[width=\picSingleWidth]{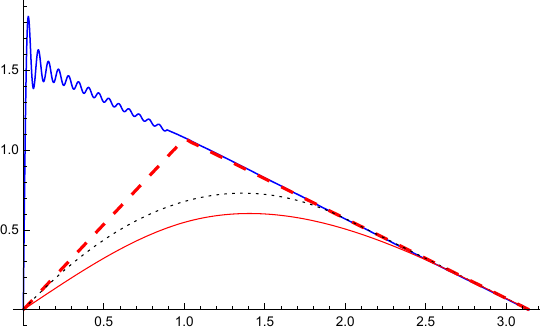}}
  \caption{The four sums in Theorem \ref{T:ThmFourEqualSums}}
  \label{fig:sincPowers}
\end{figure}
The graphs of 100 terms of these four sums are shown in Figure \ref{fig:sincPowers}.
The four sums, from left to right, have decreasing values at $x = .5$.
The four sums, from left to right, are shown in: solid blue, dashed red, dotted black, solid red.

\textbf{Proof.}
We'll prove the result for $x \leq \pi$.
The results for the larger intervals will follow by replacing $x$ with $2\pi - x$.
The first two sums are just restatements of Equations \eqref{E:five1} and \eqref{E:five2}.
We now prove that the third and fourth sums also equal $(\pi - x)/2$ over the stated intervals.

The third sum is the Fourier series of some unknown function $f(x)$ over $0 \leq x \leq \pi$:
\begin{equation}\label{E:unknownFxSeries}
f(x) = \sum_{n=1}^{\infty} \frac{\sin^2(n)}{n^3} \sin(nx) \, .
\end{equation}
But \emph{what} function would that be?
We need to find another expression for $f(x)$ which makes it clear that $f(x) = (\pi - x)/2$ for $2 \leq x \leq \pi$.

It turns out that
\begin{equation}\label{E:unknownFxQuadratic}
f(x) =
  \begin{dcases}
    \frac{x(\pi - 1)}{2} - \frac{\pi}{8} x^2  &\text{for $0 \leq x < 2$ ,} \\
    \frac{\pi - x}{2}                         &\text{for $2 \leq x \leq \pi$ ,}
  \end{dcases}
\end{equation}
provided we extend $f(x)$ for $x < 0$ by taking its odd periodic extension.
That is, if $x < 0$, we define $f(x) = -f(-x)$.
So, for $x < 0$, $f(x)$ would be:
\begin{equation*}
f(x) =
  \begin{dcases}
    - \left( \frac{\pi - (-x)}{2} \right)                            &\text{for $-\pi \leq x < -2$ ,} \\
    - \left( \frac{(-x)(\pi - 1)}{2} - \frac{\pi}{8} (-x)^2 \right)  &\text{for $-2 \leq x < 0$ .}
  \end{dcases}
\end{equation*}

In Section \ref{S:RecoveringFx}, we will see how the polynomials in \eqref{E:unknownFxQuadratic} were found.

This \textit{Mathematica} code
\begin{verbatim}
  i1 = Integrate[ Sin[n x] * ( x (Pi - 1)/2 - (Pi/8) x^2 ), {x, 0, 2}]
  i2 = Integrate[ Sin[n x] * (Pi - x)/2, {x, 2, Pi}]
  total = Simplify[ (2/Pi)(i1 + i2), Assumptions -> Element[n, Integers] ]
\end{verbatim}
gives $\sin^2(n)/n^3$ as the  $n$\textsuperscript{th} coefficient of the Fourier sine series for the $f(x)$ that is defined in Equation \eqref{E:unknownFxQuadratic},
and shows that $f(x)$ does, indeed, have the desired Fourier sine series \eqref{E:unknownFxSeries} over $0 \leq x \leq \pi$.
In particular, for $2 \leq x \leq \pi$, we do have $f(x) = (\pi - x)/2$, as claimed.

Now look at the fourth sum in Equation \eqref{E:FourEqualSums}.
Consider the function $g(x)$ such that
\begin{equation}\label{E:unknownGxSeries}
  g(x) = \sum_{n=1}^{\infty} \frac{\sin^3(n)}{n^4} \sin(nx) \, .
\end{equation}
We want an alternative expression for $g(x)$ which makes it clear that $g(x) = (\pi - x)/2$ for $3 \leq x \leq \pi$.
That expression turns out to be
\begin{equation}\label{E:unknownGxCubic}
g(x) =
  \begin{dcases}
    \left( \frac{3\pi}{8} - \frac{1}{2} \right) x - \frac{\pi}{24} x^3  &\text{for $0 \leq x \leq 1$,} \\
    -\frac{\pi}{16} + \left( \frac{9\pi}{16} - \frac{1}{2} \right) x - \frac{3\pi}{16} x^2 + \frac{\pi}{48} x^3  &\text{for $1 < x < 3$,} \\
    \frac{\pi - x}{2}                         &\text{for $3 \leq x \leq \pi$.}
  \end{dcases}
\end{equation}
In Section \ref{S:RecoveringGx}, we will see how the polynomials in \eqref{E:unknownGxCubic} were found.
This \textit{Mathematica} calculation
\begin{verbatim}
  i1 = Integrate[ Sin[n x] * ( (3 Pi/8 - 1/2) x - (Pi/24) x^3 ), {x, 0, 1}]
  i2 = Integrate[ Sin[n x] * 
       ( -Pi/16 + (9 Pi/16 - 1/2) x - (3 Pi/16) x^2 + (Pi/48) x^3 ), {x, 1, 3}]
  i3 = Integrate[ Sin[n x] * ( (Pi - x)/2 ), {x, 3, Pi}]
  total = Simplify[(2/Pi)(i1 + i2 + i3), Assumptions -> Element[n, Integers]]
\end{verbatim}
gives $\sin^3(n)/n^4$ as the  $n$\textsuperscript{th} coefficient of the Fourier sine series for the $g(x)$ in Equation \eqref{E:unknownGxCubic}.
This shows that $g(x)$ does, indeed, have the desired Fourier sine series \eqref{E:unknownGxSeries} over $0 \leq x \leq \pi$.
In particular, $g(x) = (\pi - x)/2$ for $3 \leq x \leq \pi$, as claimed.

\textbf{QED.}

Substituting $x = 2$ into Equation \eqref{E:FourEqualSums} shows that the three sums in Equation \eqref{E:ThreeEqualSumsk2} all equal $(\pi - 2)/2$.

Substituting $x = 4$ shows that the three sums in Equation \eqref{E:ThreeEqualSumsk4} all equal $(\pi - 4)/2$.

Substituting $x = 3$ in Equation \eqref{E:FourEqualSums} shows that the \emph{four} sums in Equation \eqref{E:FourSumsWithN} all equal $(\pi - 3)/2$.

Equation \eqref{E:FourSumsWithN} can also be proven using the much more general Theorem 1 of \cite{BBB}.
Suppose we are given $N > 0$ and $N + 1$ positive numbers $a_0$, $a_1$, $a_2$, ... $a_N \, .$
Among other things, Theorem 1 of \cite{BBB} says: If both of these inequalities are true
\begin{equation}\label{E:BBBTheoremPart1}
  \sum_{k = 0}^N a_k \leq 2 \pi \text{  and  } \sum_{k = 0}^N a_k \leq 2 a_0 \, ,
\end{equation}
then
\begin{equation}\label{E:BBBTheoremPart2}
  \sum_{n=1}^{\infty} \left( \prod_{k=0}^N \frac{\sin(a_k \cdot n)}{a_k \cdot n} \right)
    = \frac{\pi}{2a_0} - \frac{1}{2} \, .
\end{equation}
(Note that each term in the sum is the product of $N+1$ factors).
Using this Theorem with $N = 0$, $N = 1$, $N = 2$, and $N = 3$, we can evaluate the four sums in Equation \eqref{E:FourSumsWithN}.

For example, the second sum in Equation \eqref{E:FourSumsWithN} can be evaluated by taking $N = 1$, $a_0 = 3$, and $a_1 = 1$.
The conditions in \eqref{E:BBBTheoremPart1} are satisfied, and, Equation \eqref{E:BBBTheoremPart2} becomes
\[
  \sum_{n=1}^{\infty}
    \frac{\sin(3 n)}{3 n} \cdot
    \frac{\sin(n)}{n}
    = \frac{\pi}{6} - \frac{1}{2} \, .
\]
This is equivalent to the second sum in Equation \eqref{E:FourSumsWithN}.

Likewise, the fourth sum in Equation \eqref{E:FourSumsWithN} can be evaluated by taking
$N = 3$, $a_0 = 3$, and $a_1 = a_2 = a_3 = 1$. The conditions in \eqref{E:BBBTheoremPart1} hold, and Equation \eqref{E:BBBTheoremPart2} becomes
\[
  \sum_{n=1}^{\infty}
    \frac{\sin(3 n)}{3 n} \cdot
    \left( \frac{\sin(n)}{n} \right)^3
    = \frac{\pi}{6}
      - \frac{1}{2} \, .
\]

It is interesting that as long as the conditions in \eqref{E:BBBTheoremPart1} are satisfied, the value on the right of Equation \eqref{E:BBBTheoremPart2} depends \emph{only} on $a_0$ and \emph{not} on the values of the other $a_k$, \emph{nor} on how many $a_k$ there are.
So, for example, let $a_0 = 3$, and let $a_1 = a_2 = \dots = a_6 = 1/2$. For every $N \leq 6$, all of the sums
\[
  \sum_{k = 0}^N a_k
\]
are less than or equal to both $2\pi$ and $2 a_0$, so \eqref{E:BBBTheoremPart1} holds.
Therefore, for \emph{every} $N \leq 6$, we have
\[
  \sum_{n=1}^{\infty}
    \frac{\sin(3 n)}{3 n} \cdot \left( \frac{\sin(n/2)}{n/2} \right) ^N
    = \frac{\pi}{6} - \frac{1}{2} \, .
\]
(For $N = 7$, the result is a polynomial in $\pi$ of degree 8.)

\textbf{Exact Calculations.}\\
Above, we displayed the values of a few sums to about 20 decimal places.
Besides being able to do such numerical calculations to extremely high precision, computer algebra systems like \textit{Mathematica} can also compute \emph{exact} values of some of these expressions.

For example, \textit{Mathematica} claims that all three of the following sums are \emph{exactly} equal to $(\pi - 2)/2$:
\[
\sum_{n=1}^{\infty} \frac{\sin(2n)}{n} =
\sum_{n=1}^{\infty} \frac{\sin(n)}{n} \frac{\sin(2n)}{n} =
\sum_{n=1}^{\infty} \left(\frac{\sin(n)}{n} \right) ^{2} \frac{\sin(2n)}{n} =
\frac{\pi - 2}{2} \, .
\]
We proved this above; this is Equation \eqref{E:ThreeEqualSumsk2}.
\textit{Mathematica} also says that this sum is exactly equal to $23 \pi/48 - 1$:
\[
\sum_{n=1}^{\infty} \left(\frac{\sin(n)}{n} \right) ^{3} \frac{\sin(2n)}{n} = \frac{23\pi}{48} - 1 \approx 0.50534 \text{ } 64798 \, .
\]

Here's the \textit{Mathematica} code to do these calculations, and a few others.
We will calculate six values, corresponding to sums involving $(\sin(n)/n)^0 = 1$ through $(\sin(n)/n)^5$.
First, we define a function, \verb+sumk2[k]+, whose parameter, \verb+k+, is the power of $\sin(n)/n$ in the sum.
Then, we calculate a table of six expressions, \verb+sumk2[0]+ through \verb+sumk2[5]+.
(The trailing semicolon in the code below prevents the display of the lengthy output of these expressions).
Finally, we simplify the results as much as possible:
\begin{verbatim}
  sumk2[k_] := Sum[(Sin[n]/n)^k * Sin[2 n]/n, {n, 1, Infinity}]
  tab = Table[sumk2[k], {k, 0, 5}] ;
  FullSimplify[tab]
\end{verbatim}
This output is a table with six expressions.
The first three expressions are $(\pi - 2)/2$. We proved these results above.
The fourth expression, corresponding to the sum with $(\sin(n)/n)^3$, is $23\pi/48 - 1$.
In the fifth expression, \textit{Mathematica} is claiming that
\[
  \sum_{n=1}^{\infty} \left(\frac{\sin(n)}{n} \right) ^{4} \frac{\sin(2n)}{n} = \frac{11\pi}{24} - 1 \, .
\]
Finally, the sixth expression in \textit{Mathematica}'s output claims that
\begin{align*}
  \sum_{n=1}^{\infty} \left(\frac{\sin(n)}{n} \right) ^{5} \frac{\sin(2n)}{n} & =
  \frac{-3840 + \pi  (18489-2 \pi  (12005+4 \pi  (-1715 + 2 \pi  (245+\pi  (-35 + 2 \pi)))))}{3840} \\
  & = 
  -1+\frac{6163 \pi }{1280}-\frac{2401 \pi ^2}{384}+\frac{343 \pi ^3}{96}-\frac{49
   \pi ^4}{48}+\frac{7 \pi ^5}{48}-\frac{\pi ^6}{120}
  \, ,
\end{align*}

These results look intriguing, but we shall not prove them here; see \cite{BBB} for details.

\textbf{What Constitutes a Proof?}\\
It is important to note the following: because we are unable to examine or verify the internal, proprietary algorithms that \textit{Mathematica} uses to do these exact calculations, such output does not, in any sense, constitute a proof.
Still, such symbolic output is useful because it gives us the likely answer, so it tells us what it is that we should attempt to prove.

That said, we are cheating a little bit when we allow a computer algebra system like \textit{Mathematica} to evaluate an integral like
\begin{verbatim}
  i1 = Integrate[ Sin[n x] * ( (3Pi/8 - 1/2) x - (Pi/24) x^3 ), {x, 0, 1}]
\end{verbatim}
whose output is
\[
\frac{\left(  n^3 (6-4 \pi )-3 n \pi \right) \cos (n) + 3 \left(  n^2 (\pi -2) + \pi \right) \sin (n)}{12 n^4}
\]
if we don't also check the result by some other method.
Although this calculation is somewhat tedious, it can, in theory, be done by hand.
The reader could also gain confidence in the result by doing the same calculation in a \emph{different} computer algebra system.
Of course, doing such a calculation by hand instead of in \textit{Mathematica} does \emph{not} guarantee that it is correct!
However, if you show all the steps, this at least enables others to check your work.

\textbf{Replacing $\sin$ With $\cos$.}\\
Equation \eqref{E:ThreeEqualSumsk2} (restated here for convenience) says:
\[
\sum_{n=1}^{\infty} \frac{\sin(2n)}{n} =
\sum_{n=1}^{\infty} \frac{\sin(n)}{n} \frac{\sin(2n)}{n} =
\sum_{n=1}^{\infty} \left(\frac{\sin(n)}{n} \right) ^{2} \frac{\sin(2n)}{n} =
\frac{\pi - 2}{2} \, .
\]
If we use the identity $\sin(2n) = 2 \sin(n) \cos(n)$, we get pretty sums involving cosines:
\begin{equation} \label{E:PrettyCosSums}
\sum_{n=1}^{\infty} \frac{\sin(n)}{n} \cos (n) =
\sum_{n=1}^{\infty} \left(\frac{\sin(n)}{n} \right)^{2} \cos (n) =
\sum_{n=1}^{\infty} \left(\frac{\sin(n)}{n} \right)^{3} \cos (n) =
\frac{\pi -2}{4} \, .
\end{equation}

\begin{figure}[ht]
  \mbox{\includegraphics[width=\picSingleWidth]{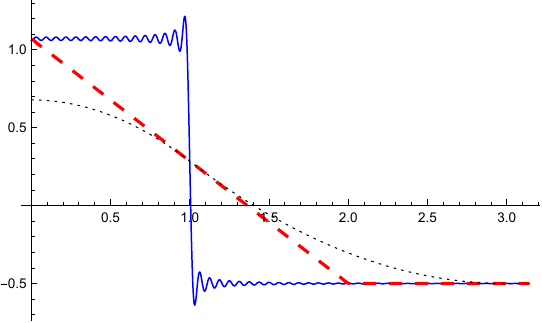}}
  \caption{Sum of $\left( \frac{\sin(n)}{n} \right)^k \cos(nx)$ for $k = 1, 2, 3$}
  \label{fig:sincToKcosPlot}
\end{figure}

The graphs of
\[
\sum_{n=1}^{\infty} \left(\frac{\sin(n)}{n} \right) ^{k} \cos (nx)
\]
for $k = 1$, 2, and 3 illustrate why \eqref{E:PrettyCosSums} holds for $x = 1$.
In Figure \ref{fig:sincToKcosPlot}, $k = 1$ is the solid graph, $k = 2$ is the dashed graph, and $k = 3$ is the dotted graph.
All three curves cross at $x = 1$, $y = (\pi - 2)/4 \approx 0.285 $.

\begin{figure}[ht]
  \mbox{\includegraphics[width=\picSingleWidth]{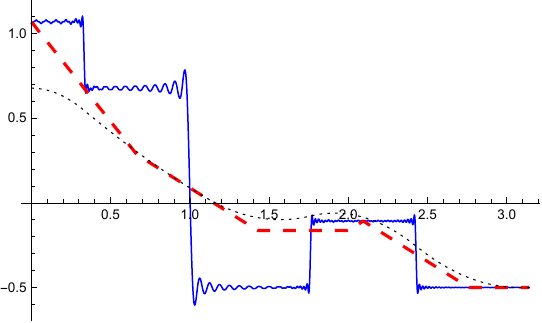}}
  \caption{Sum of $\left( \frac{\sin(n)}{n} \right)^k \cos^3(nx)$ for $k = 1, 2, 3$}
  \label{fig:sincToKcos3Plot}
\end{figure}

The reader might also be intrigued by the graphs of the following sums for $k = 1$, $2$, and $3$:
\[
\sum_{n=1}^{\infty} \left(\frac{\sin(n)}{n} \right) ^{k} \cos ^{3} (nx) \, .
\]
The graphs are in Figure \ref{fig:sincToKcos3Plot}.
The graphs for $k = 1$, $k = 2$, and $k = 3$ are solid, dashed, and dotted, respectively.
The graphs all appear to meet at the point $x = 1$, $y \approx \, .09$.
In fact, \textit{Mathematica} claims that, for $k = 1$, 2, and 3, the above sums are exactly $3\pi/16 - 1/2 \approx \, .089$, but we will not prove this.
For $k = 4$, the sum is about 0.075699, so this sum is different from the sums for $k = 1$, 2, and 3.

\begin{figure}[ht]
  \mbox{\includegraphics[width=\picSingleWidth]{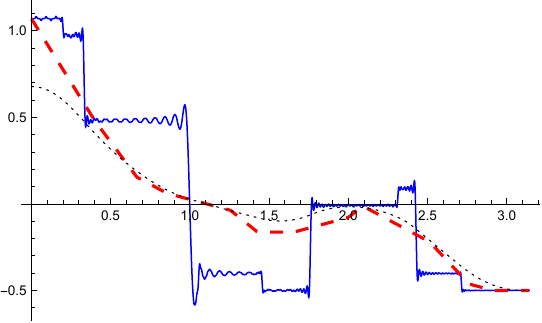}}
  \caption{Sum of $\left( \frac{\sin(n)}{n} \right)^k \cos^5(nx)$ for $k = 1, 2, 3$}
  \label{fig:sincToKcos5Plot}
\end{figure}

Question: for $k = 1$, 2, and 3, are any of the following sums equal at $x = 1$?
\[
  \sum_{n=1}^{\infty} \left(\frac{\sin(n)}{n} \right) ^{k} \cos ^{5} (nx) \, .
\]
(Hint: seeing is \emph{not} believing!)
The graphs are in Figure \ref{fig:sincToKcos5Plot}.
The graphs for $k = 1$, 2, and 3 are solid, dashed, and dotted, respectively.

Answer: at $x = 1$, the sums are all different.
For $k = 1$, 2, and 3, the sums at $x = 1$ are about $-0.00912$, 0.02606, and 0.02704.

\section{Recovering a function from its fourier series}\label{S:Recovering}

In the proof of Theorem \ref{T:ThmFourEqualSums} in the last section, we started with functions $f(x)$ and $g(x)$ defined by their Fourier series \eqref{E:unknownFxSeries} and \eqref{E:unknownGxSeries}.
The proof required expressing these functions using polynomials in \eqref{E:unknownFxQuadratic} and \eqref{E:unknownGxCubic}.
In this section, we will show how those polynomials were found.

It's worth stating that Section \ref{S:UsingPolylogs} shows how polylogarithms let us \emph{derive} the polynomials directly from the Fourier series.
That method eliminates the guesswork that we employ here.



For $0 \leq x \leq 2 \pi$, we know that
\begin{equation} \label{E:pmxOver2}
\frac{\pi - x}{2} = \sum_{n=1}^{\infty} \frac{\sin(n x)}{n} \, .
\end{equation}
Also, Equation \eqref{E:Zeta2} states that, for $0 \leq x \leq 2 \pi$,
\[
\frac{3 x^2 - 6 \pi x + 2 \pi^2}{12} =
  \sum_{n=1}^{\infty} \frac{\cos(n x)}{n^2} \, .
\]

\end{comment}

But suppose we didn't know, but \emph{guessed}, that, over $0 \leq x \leq 2 \pi$, \emph{some} quadratic has the Fourier series
\begin{equation} \label{E:QuadraticCBA}
C x^2 + B x + A = \sum_{n=1}^{\infty} \frac{\cos(n x)}{n^2} \, .
\end{equation}

How can we find these coefficients?
First, we can compute the Fourier series for $A$, $B x$, and $C x^2$.
The Fourier ``series'' for the constant $A$ is just $A$, with no sums.
From Equation \eqref{E:pmxOver2},
\[
B x = B \pi - 2 B \sum_{n=1}^{\infty} \frac{\sin(n x)}{n} \, .
\]
If we do the integrations by parts to compute the Fourier series of $C x^2$, we get
\[
C x^2 = 4C \frac{\pi^2}{3} +
4 C     \sum_{n=1}^{\infty} \frac{\cos(n x)}{n^2} -
4 C \pi \sum_{n=1}^{\infty} \frac{\sin(n x)}{n} \, .
\]
So, a general quadratic has the Fourier series
\[
C x^2 + B x + A =
\left( A + B \pi + 4 C \frac{\pi^2}{3} \right)
+ 4 C \sum_{n=1}^{\infty} \frac{\cos(n x)}{n^2}
- (2B + 4 C \pi) \sum_{n=1}^{\infty} \frac{\sin(n x)}{n} \, .
\]
To match the right-hand side of Equation \eqref{E:QuadraticCBA}, we must have
\[
4 C = 1 \quad \text{ so } \quad C = \frac{1}{4} \, ,
\]
\[
A + B \pi + 4 C \frac{\pi^2}{3} = A + B \pi + \frac{\pi^2}{3} = 0 \, ,
\]
and
\[
2 B + 4 C \pi = 2 B + \pi = 0 \, .
\]
From this last equation, $B = -\pi/2$, so
\[
0 = A + B \pi + \frac{\pi^2}{3} = A -\frac{\pi^2}{2} + \frac{\pi^2}{3} = A - \frac{\pi^2}{6} \, ,
\]
so $A = \pi^2/6$.
Then we have
\begin{equation} \label{E:FS-SpecialQuadratic}
\frac{x^2}{4} - \frac{\pi x}{2} + \frac{\pi^2}{6}
= \frac{3 x^2 - 6 \pi x + 2 \pi^2}{12}
= \sum_{n=1}^{\infty} \frac{\cos(n x)}{n^2} \, .
\end{equation}
In a similar fashion, we can calculate that, for $0 \leq x \leq 2 \pi$,
\[
\frac{x^3 - 3 \pi x^2 + 2 \pi^2 x}{12}
= \sum_{n=1}^{\infty} \frac{\sin(n x)}{n^3} \, ,
\]
and
\[
-\frac{15 x^4 - 60 \pi x^3 + 60 \pi^2 x^2 - 8 \pi^4 }{720}
= \sum_{n=1}^{\infty} \frac{\cos(n x)}{n^4} \, .
\]

In general, for $0 < x < 2 \pi$ (so that $0 < x/(2 \pi) < 1$), and for any integer $k \geq 1$, the sums
\[
\sum_{n=1}^{\infty} \frac{\sin(n x)}{n^{2 k - 1}} \quad\ \text{ and } \quad
\sum_{n=1}^{\infty} \frac{\cos(n x)}{n^{2 k}}
\]
can be expressed in terms of Bernoulli polynomials $B_j(\cdot)$; see \cite[Eqs. 24.8.1 and 24.8.2]{DLMF}:
\begin{equation}  \label{E:BernoulliEven}
 B_{2k}\left( \frac{x}{2 \pi} \right) = (-1)^{k+1} \; \frac{ 2 \cdot (2k)! }{ (2 \pi)^{2k} } \; \sum_{n=1}^{\infty} \frac{\cos(n x)}{n^{2k}} \, ,
\end{equation}
and
\begin{equation}  \label{E:BernoulliOdd}
 B_{2k-1}\left( \frac{x}{2 \pi} \right) = (-1)^{k} \; \frac{2 \cdot (2k-1)! }{ (2 \pi)^{2k-1} } \; \sum_{n=1}^{\infty} \frac{\sin(n x)}{n^{2k-1}} \, .
\end{equation}
With the code below, we can use the \verb+fSeries[ ]+ function in the \verb+FS.m+ package to check a few instances of Equations \eqref{E:BernoulliOdd} and \eqref{E:BernoulliEven}.
See Section \ref{S:Package} for details on this software package.
\begin{verbatim}
  bCoeff[k_Integer?Positive] :=
  Module[
  (* DLMF, eq 24.8.1, 24.8.2; Abramowitz and Stegun, eq 23.1.17 and 23.1.18. *)
    { expo },
    If[OddQ[k] , expo = (k + 1)/2 , expo = (k/2) + 1 ];
    Return[ (-1)^expo * 2 * k! / (2 Pi)^k ]
  ];  (* end of Module bCoeff *)

  Table[ fSeries[x, BernoulliB[k, x/(2 Pi)]/bCoeff[k], 0, 2 Pi, n] , {k, 1, 6} ]
\end{verbatim}
Items 1, 3, and 5 below confirm \eqref{E:BernoulliOdd};
items 2, 4, and 6 confirm \eqref{E:BernoulliEven}:
\[
\left \{ 0 , 0 , \frac{1}{n} \right \} ,
\left\{ 0 , \frac{1}{n^2} , 0 \right \} ,
\left \{ 0 , 0 , \frac{1}{n^3} \right \} ,
\left \{ 0 , \frac{1}{n^4} , 0 \right \} ,
\left \{ 0 , 0 , \frac{1}{n^5} \right \} ,
\left \{ 0 , \frac{1}{n^6} , 0 \right \} \, .
\]

\subsection{Identifying f(x)}\label{S:RecoveringFx}

\mbox{}  

Equations \eqref{E:unknownFxSeries} and \eqref{E:unknownFxQuadratic} claimed that $f(x)$, defined as a series, has a polynomial representation:
\begin{equation*}
f(x) = \sum_{n=1}^{\infty} \frac{\sin^2(n)}{n^3} \sin(nx) =
  \begin{dcases}
     \frac{\pi - 1}{2} x - \frac{\pi}{8} x^2  &\text{for $0 \leq x < 2$,} \\
    \frac{\pi - x}{2}                         &\text{for $2 \leq x \leq \pi$.}
  \end{dcases}
\end{equation*}
How did we obtain these polynomials?

\begin{figure}[ht]
  \mbox{\includegraphics[width=\picSingleWidth]{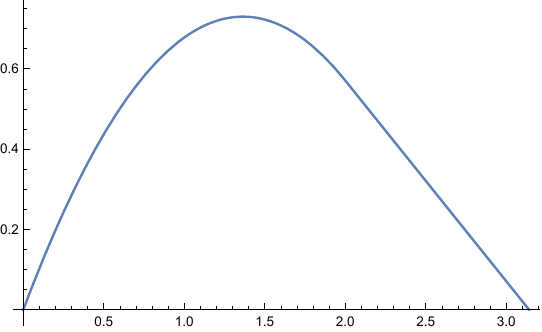}}
  \caption{$f(x) = \sum_{n=1}^{100} \frac{\sin(n)^2}{n^3} \cdot \sin(nx)$}
  \label{fig:fxPlot23}
\end{figure}

The graph of the series in \eqref{E:unknownFxSeries} is shown in Figure \ref{fig:fxPlot23}.
Unlike the graphs we encountered in earlier sections, such as Figures \ref{fig:f1Plot} and \ref{fig:gxPlot},
here we cannot easily read off from the graph any pieces that may comprise $f(x)$.

As noted in the discussion point \textbf(3) after the proof of Theorem \ref{T:Thm1A}, the presence of $\sin(n)$ in the Fourier coefficients probably means that $f(x)$ is composed of two or more functions over intervals with integer endpoints.
The graph in Figure \ref{fig:fxPlot23} \emph{appears} to consist of two (or more) pieces: first, quadratic or higher-degree polynomials over $0 < x < 2$, and a first-degree (linear) polynomial for $2 < x < \pi$. (We'll worry about the endpoints later.)
Of course, there's no guarantee that $f(x)$ \emph{does} consist of polynomials.
If it does not, then we may be out of luck.
Note that, in order to match the right side of Equation \eqref{E:FourEqualSums}, the linear part would have to be $(\pi - x)/2$, but we will show how to find this expression without assuming it.

Figure \ref{fig:fxPlot23} gives us a couple of hints about how to proceed.
First, let's assume that $f(x)$ is represented by just two polynomials: one over (0, 2) and another over (2, $\pi $).
We can then compute the sum in \eqref{E:unknownFxSeries} to reasonable precision to get pairs of data points $(x, f(x))$.
Then, we can do least-squares curve fitting over the two intervals to try to find polynomials that match the data.
Finally, we will try to recognize the exact values of the numeric coefficients, and write the coefficients as expressions.

Let's start by fitting $f(x)$ from $x = 0$ to $x = 2$ to a cubic polynomial.
We'll use 10000 terms of the series for $f(x)$ to compute a set of $(x, y)$ data points from $x = 0$ to $x = 2$, then use \textit{Mathematica}'s \verb+Fit+ function to fit the data to a cubic.
Why a cubic? This is just a guess.
$f(x)$ might be a polynomial of higher degree, or it might not be a polynomial at all.
The only way to find out is trial and error: fit $f(x)$ to a cubic and see what happens.

This \textit{Mathematica} code (the trailing semicolon suppresses lengthy output that we don't need)
\begin{verbatim}
  fx[x_] := Sum[Sin[n]^2/n^3 * Sin[n x], {n, 1, 10000}]
  data02 = Table[{x, fx[x]}, {x, 0, 2, .01}] ;  (* x from 0 to 2 *)
  Fit[data02, {1, x, x^2, x^3}, x]
\end{verbatim}
gives the polynomial
\[
-2.13691153825 \cdot 10^{-12} + 1.0707963268 x - 0.392699081705 x^2 + 1.36298199138 \cdot 10^{-12} x^3 \, .
\]
The last few digits of these coefficients may be meaningless because the data came from only the first 10000 terms of the series for $f(x)$.
The relative smallness of the constant term and the coefficient of $x^3$ suggest that these terms really should be 0.

So, let's try again, this time fitting the data to a polynomial of the form $a x + b x^2$.
This time, \verb+Fit[data02, {x, x^2}, x]+ gives $1.0707963268 x - 0.392699081699 x^2 \, $.

In the previous sections, coefficients often turned out to be simple fractions involving $\pi$.
With that in mind, it doesn't take long to see that the coefficient of $x$ is very close to $(\pi - 1)/2 \approx 1.07079632679$.
Also, $\pi/0.392699081699 \approx 7.99999999999$, so 0.392699081699 is probably an approximation to $\pi/8 = 0.3926990816987 \dots$.
So, we conjecture that, between $x = 0$ and $x = 2$, $f(x)$ has the polynomial representation
\[
f(x) = \frac{\pi - 1}{2} x - \frac{\pi}{8} x^2 \, .
\]

Figure \ref{fig:fxPlot23} suggests that $f(x)$ is linear from $x = 2$ to $x = \pi$.
But just in case there is also a quadratic term, let's fit a quadratic over this interval.
\begin{verbatim}
  data2Pi = Table[{x, fx[x]}, {x, 2, Pi, .01}] ;  (* x from 2 to Pi *)
  Fit[data2Pi, {1, x, x^2}, x]
\end{verbatim}
gives the polynomial
\[
1.57079632682 - 0.500000000016 x + 2.92866805496 \cdot 10^{-12} x^2 \, .
\]
The quadratic term seems to be 0, so let's fit a linear function.
\verb+Fit[data2Pi, {1, x}, x]+ gives $1.5707963268 - 0.500000000001 x$.
This suggests that from $x = 2$ to $x = \pi$,
\[
f(x) = \frac{\pi}{2} - \frac{x}{2} = \frac{\pi - x}{2} \, .
\]
This is how the polynomial expressions in Equation \eqref{E:unknownFxQuadratic} were obtained.

\textbf{Alternative Method.}
Here's another approach that does not involve trying to guess the exact values of the imprecise decimal values above.
Again, we do not assume that $f(x) = (\pi - x)/2$ for $x$ between 2 and $\pi$, although we hope that turns out to be the case.

Let's assume that $f(x)$ is the quadratic $f_1(x) = ax^2 + bx$ over $0 \leq x < 2$.
(From the graph, $f(0)$ appears to be 0, so the constant term in the quadratic is 0).
Let's also assume $f(x)$ is the linear function $f_2(x) = cx + d$ over $2 < x \leq \pi$.
Let's also assume that both the function values and the slopes agree at $x = 2$, and that $f_2(\pi) = 0$.
Then the coefficients must be related to each other in the following ways (here, the primes $'$ denote derivatives):
\begin{align*}
 f_1(2) & = f_2(2),   & \text{ so } 4a + 2b = 2c + d, \\
 f_1'(2) & = f_2'(2), & \text{ so } 4a + b = c, \\
 f_2(\pi) & = 0,      & \text{ so } \pi c + d = 0 \, .
\end{align*}

If we solve these for $a$, $b$, and $d$ in terms of $c$, we get $a = \pi c /4$, $b = c(1-\pi)$, and $d = -c \pi$.
Next, we compute the integrals \eqref{E:bnOdd} to get the coefficients of the Fourier sine series:
\begin{verbatim}
  i1 = (2/Pi) Integrate[Sin[n x] (a x^2 + b x), {x, 0, 2}]
  i2 = (2/Pi) Integrate[Sin[n x] (c x + d),     {x, 2, Pi}]
  expr1 = FullSimplify[i1 + i2, Assumptions -> Element[n, Integers]]
\end{verbatim}
(As before, \verb+Assumptions -> Element[n, Integers]+ sets terms like $\sin(n\pi)$ to 0).
This gives us another expression involving the parameters $a$, $b$, $c$, and $d$:
\[
\frac{2}{\pi } \int_{0}^{2}(ax^{2} +bx)\sin(nx)dx +\frac{2}{\pi } \int_{2}^{\pi }(cx+d)\sin(nx)dx =
\]
\[
\frac{-4 a-2 (-1)^n n^2 (d + c \pi) +  \left(4 a + 2 (-4 a-2 b+2 c+d) n^2\right) \cos (2 n)  + 2 (4 a+b-c) n \sin(2 n)}{n^3 \pi} \, .
\]

This looks pretty bad!
However, we know $a$, $b$, and $d$ in terms of $c$.
When we make those substitutions, we get an expression involving only $c$.
\begin{verbatim}
  expr2 = expr1 /. {a -> Pi c/4, b -> c(1 - Pi), d -> -c Pi}
  Simplify[expr2]
\end{verbatim}
The result of this \textit{Mathematica} calculation is $-2 c \sin^2(n)/n^3 \, .$

In order for this to match the coefficient of $\sin(nx)$ in \eqref{E:unknownFxSeries}, we must have $c = -1/2$.
Then, $a = -\pi/8$, $b = (\pi - 1)/2$, and $d = \pi/2$, and the two functions are just what we claimed above in Equation \eqref{E:unknownFxQuadratic}:
\[
f_{1} (x) = \frac{\pi -1}{2} x - \frac{\pi}{8} x^{2}
\]
and
\[
f_{2} (x) = \frac{\pi}{2} - \frac{1}{2} x \, .
\]

\subsection{Identifying g(x)}\label{S:RecoveringGx}

\mbox{}  

Now, let's look at $g(x)$ in Equations \eqref{E:unknownGxSeries} and \eqref{E:unknownGxCubic}.
These equations claimed that
\begin{equation*}
g(x) = \sum_{n=1}^{\infty} \frac{\sin^3(n)}{n^4} \sin(nx) =
  \begin{dcases}
    \left( \frac{3\pi}{8} - \frac{1}{2} \right) x - \frac{\pi}{24} x^3  &\text{for $0 \leq x \leq 1$,} \\
    -\frac{\pi}{16} + \left( \frac{9\pi}{16} - \frac{1}{2} \right) x - \frac{3\pi}{16} x^2 + \frac{\pi}{48} x^3  &\text{for $1 < x < 3$,} \\
    \frac{\pi - x}{2}                         &\text{for $3 \leq x \leq \pi$.}
  \end{dcases}\end{equation*}

\begin{figure}[ht]
  \mbox{\includegraphics[width=\picSingleWidth]{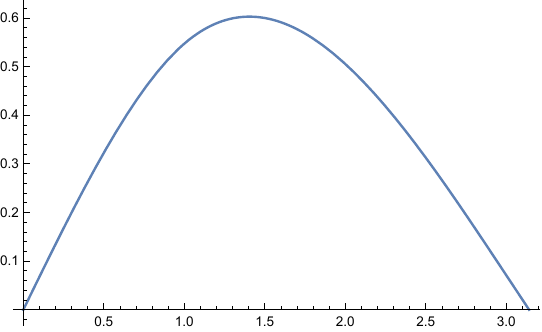}}
  \caption{$g(x) = \sum_{n=1}^{100} \frac{\sin(n)^3}{n^4} \cdot \sin(nx)$}
  \label{fig:gxPlot34}
\end{figure}
Figure \ref{fig:gxPlot34} shows the graph of $g(x)$ from $x = 0$ to $x = \pi$.
This graph is smoother than the graph of $f(x)$ in Figure \ref{fig:fxPlot23}.
Unlike the graph of $f(x)$, this graph doesn't give any visible hints about its pieces, or about the endpoints of intervals on the $x$ axis.
All we expect (based on the statement of Theorem \ref{T:ThmFourEqualSums}) is that, from $x = 3$ to $x = \pi$, $g(x) = (\pi - x)/2$.

Again, the presence of $\sin(n)$ in the coefficients of the Fourier series suggest that $g(x)$ might consist of pieces defined over intervals with integer endpoints.
Curve-fitting with even a 5\textsuperscript{th} degree polynomial over [0, 3] did not produce a good match to the data.
This suggested that, over [0, 3], either $g(x)$ was a polynomial of degree 6 or more, or that $g(x)$ was composed of two or more polynomials over several subintervals of [0, 3].
Trial and error curve-fitting over a number of intervals finally gave good matches over the three intervals [0, 1], [1, 2], and [2, 3].
The polynomials over [1, 2] and [2, 3] were the same, suggesting that just one polynomial over [1, 3] will work.
From $x = 3$ to $x = \pi$, the first-degree polynomial $(\pi - x)/2$ fit the data pretty well.
All of this implied that $g(x)$ was probably composed of polynomials over thee three intervals [0, 1], [1, 3], and [3, $\pi $].

First, using the series definition of $g(x)$ in Equation \eqref{E:unknownGxSeries}, let's fit $g(x)$ over the interval from $x = 0$ to $x = 1$.
These denominators contain $n^4$, so the terms in the series get small very quickly.
So, we will use only 1000 terms of this series to generate the data.
\begin{verbatim}
  gx[x_] := Sum[Sin[n]^3/n^4 * Sin[n x], {n, 1, 1000}]
  data01 = Table[{x, gx[x]}, {x, 0, 1, .01}] ;  (* x from 0 to 1 *)
  Fit[data01, {1, x, x^2, x^3, x^4}, x]
\end{verbatim}
gives a polynomial with tiny constant, $x^2$, and $x^4$ coefficients (on the order of $10^{-9}$ or smaller).
This suggests that, from $x = 0$ to $x = 1$, $g(x)$ may be a polynomial with only $x$ and $x^3$ terms.

In fact, \verb+Fit[data01, {x, x^3}, x]+ gives $0.678097245105 x - 0.13089969392 x^3 \, .$

What are these numbers, 0.678097245105 and $-0.13089969392$ ? 
For example, are they related to $\pi$ in some simple way?
We can try looking them up in the on-line Inverse Symbolic Calculator \cite{ISC}, but that did not help in this case.
Keep in mind that the last few digits might be meaningless.
One possibility is that these two numbers are fractions that involve $\pi$.
We find that $\pi/-0.13089969392 = -23.9999999963$, so $-0.13089969392$ might be an approximation to $-\pi/24 \approx -0.13089 96938 9957 \dots$.

If we wish to search more systematically, we can use \textit{Mathematica}'s \verb+LatticeReduce+ function to help identify this number.
\verb+LatticeReduce+ looks for non-trivial linear relations among its input values.
You can learn more about the technique of ``lattice reduction'' at \cite{Lattice}.
This code shows how we can find a linear expression that involves the three numbers $c = -0.13089969392$, $\pi$, and 1:
\begin{verbatim}
  c = -0.13089969392
  v = {c, Pi, 1}
  {a0, a1, a2} = Round[10^10 v]
  a = { {1, 0, -a0}, {0, 1, -a1} }
  b = LatticeReduce[a]
\end{verbatim}
The result is $b = \{ \{24, 1, 0\}, \{1, 0, 1308996939\} \}$.
There are two rows in the vector $b$, which means that two linear relationships were found.
The first row has small coefficients $b_{11} = 24$, $b_{12} = 1$, $b_{13} = 0$, so we'll use them.
The linear relation implied by this first row is
\[
b_{11} \cdot v_1 + b_{12} \cdot v_2 + b_{13} \cdot v_3 \, .
\]
Here, $v_1$, $v_2$, and $v_3$ are the components of the vector $v$, namely, $c$, $\pi$, and 1.
We can write this dot product in \textit{Mathematica} as \verb+b[[1]].v+, which \textit{Mathematica} says is about $-4.9 \cdot 10^{-10}$.
Written out, this dot product is
\[
24 c + 1 \cdot \pi + 0 \cdot 1 = 24 c + \pi \approx -4.9 \cdot 10^{-10} \, .
\]

Since $24 c + \pi \approx 0$, we have $c \approx -\pi/24$.

Note: the second row of the vector $b$ gives another linear relationship.
However, this row contains a 10-digit number.
This is about the number of digits as in our input value $c = -0.13089969392$, so we should not expect this to be useful.
The second row gives rise to the relation
\[
b_{21} \cdot v_1 + b_{22} \cdot v_2 + b_{23} \cdot v_3 \approx 1.3 \cdot 10^9 \, .
\]

The reader may wonder, ``How did we know that $-0.13089969392$ did not involve, say, $\pi^2$?''
We didn't know this at first, but we can check for this possibility, too:
\begin{verbatim}
  v = {c, Pi, Pi^2, 1} ;
  {a0, a1, a2, a3} = Round[10^10 v] ;
  a = { {1, 0, 0, -a0}, {0, 1, 0, -a1}, {0, 0, 1, -a2} } ;
  b = LatticeReduce[a]
  b[[1]].v
\end{verbatim}
The first row of $b$ is $\{24, 1, 0, 0\}$, which implies that
\begin{align*}
&b_{11} \cdot v_1 + b_{12} \cdot v_2 + b_{13} \cdot v_3 + b_{14} \cdot v_4 \\
 = & \, 24 \cdot c + 1 \cdot \pi + 0 \cdot \pi^2 + 0 \cdot 1 \\
  \approx & -4.9 \cdot 10^{-10} \, .
\end{align*}
Note that the coefficient of $\pi^2$ is 0, which means that $\pi^2$ does not enter into this linear relation.

The second and third rows of $b$ are 
\[
  \{-1997, 47669, -15200, 35433\} \text{ and } \{-1273, 30367, -9683, -63546\}
\]
which together have more digits than our input value $-0.13089969392 \, .$
If you allow coefficients to have enough digits, you can always write down linear relations between any set of numbers, so we'll ignore these two rows of $b$.

Now let's identify $c = 0.678097245105$, the coefficient of $x$.
It turns out that this code
\begin{verbatim}
  c = 0.678097245105 ;
  v = {c, Pi, 1} ;
  {a0, a1, a2} = Round[10^10 v] ;
  a = { {1, 0, -a0}, {0, 1, -a1} } ;
  b = LatticeReduce[a]
\end{verbatim}
does not work well because all rows of $b$ have large coefficients.
So, try again, with an additional 1 in the input vector:
\begin{verbatim}
  v = {c, 1, Pi, 1} ;
  {a0, a1, a2, a3} = Round[10^10 v] ;
  a = { {1, 0, 0, -a0}, {0, 1, 0, -a1}, {0, 0, 1, -a2} } ;
  b = LatticeReduce[a]
  b[[1]].v
\end{verbatim}
This gives $\{8, 4, -3, 0\}$ as the first row of $b$.
Also, \verb+b[[1]].v+ $\approx 7.1*10^{-11}$.
So, we have
\begin{align*}
&b_{11} \cdot v_1 + b_{12} \cdot v_2 + b_{13} \cdot v_3 + b_{14} \cdot v_4 \\
  = & 8 \cdot c + 4 \cdot 1 - 3 \cdot \pi + 0 \cdot 1 \\
  = & 8c + 4 - 3\pi \\
  \approx & \, 7.1*10^{-11} \approx 0 \, .
\end{align*}
Now, if $8c + 4 - 3\pi = 0$, then $c = (3\pi - 4)/8 \approx 0.678097245096 \, .$
Finally, putting all this together, we now suspect that, from $x = 0$ to $x = 1$,
\[
g(x) = \frac{3 \pi - 4}{8} x - \frac{\pi}{24}x^3 \, .
\]
That takes care of the polynomial that comprises $g(x)$ from $x = 0$ to $x = 1$.

Now, we'll fit $g(x)$ to a polynomial from $x = 1$ to $x = 3$.
\begin{verbatim}
  data13 = Table[{x, gx[x]}, {x, 1, 3, .01}] ;  (* x from 1 to 3 *)
  Fit[data13, {1, x, x^2, x^3, x^4}, x]
\end{verbatim}
gives $1.5 \cdot 10^{-11}$ as the coefficient of $x^4$, so we'll assume that the polynomial is a cubic.

\verb+Fit[data13, {1, x, x^2, x^3}, x]+ gives
\[
-0.19634954096 + 1.26714586782 x - 0.589048622634 x^2 + 0.0654498469637 x^3
\]

We'll start with the simple idea of dividing the coefficients by $\pi$:
\begin{align*}
-0.19634954096/\pi & \approx -.0625 = -1/16 \\
-0.589048622634/\pi & \approx -.1875 = -3/16 \\
 0.0654498469637/\pi & \approx 0.0208333333378 \approx 1/48
\end{align*}

Therefore, it appears that the polynomial expression for $g(x)$ between $x = 1$ and $x = 3$ is
\[
-\frac{\pi}{16} + 1.26714586782 x - \frac{3 \pi}{16} x^2 + \frac{\pi}{48} x^3 \, .
\]
However, $1.26714586782/\pi \approx 0.403345056964$ is not an easily-recognizable fraction, so we'll use \verb+LatticeReduce+ again.
\begin{verbatim}
  c = 1.26714586782 ;
  v = {c, 1, Pi, 1} ;
  {a0, a1, a2, a3} = Round[10^10 v] ;
  a = { {1, 0, 0, -a0}, {0, 1, 0, -a1}, {0, 0, 1, -a2} } ;
  b = LatticeReduce[a]
  b[[1]].v
\end{verbatim}
The first row of $b$ has coefficients $\{-16, -8, 9, 24\}$.
\textit{Mathematica} says that the dot product of this row with $v$ is \verb+b[[1]].v+ $\approx 23.9999999972 \,$.
The dot product of this row with $v$ can be written out as
\begin{align*}
&b_{11} \cdot v_1 + b_{12} \cdot v_2 + b_{13} \cdot v_3 + b_{14} \cdot v_4 \\
  = & -16 \cdot c - 8 \cdot 1 + 9 \cdot \pi + 24 \cdot 1 \\
  = & -16c - 8 + 9 \pi + 24 \\
  \approx & \, 23.9999999972 \approx 24 \, .
\end{align*}

If $-16c - 8 + 9\pi + 24 \approx 24$, then $c \approx (9\pi - 8)/16$.

Therefore, it appears that from $x = 1$ to $x = 3$,
\[
g(x) = -\frac{\pi}{16} + \frac{(9\pi - 8)}{16} x - \frac{3\pi}{16} x^2 + \frac{\pi}{48} x^3 \, .
\]
This is the second polynomial in Equation \eqref{E:unknownGxCubic}.

Finally, between $x = 3$ and $x = \pi$, the expression $(\pi - x)/2$ is a good fit to the data from the series in Equation \eqref{E:unknownGxSeries}.
This gives us the last part of Equation \eqref{E:unknownGxCubic}.

To summarize, what we have done is this:
\begin{itemize}
  \item (1) use the first thousand terms of the series for $g(x)$ to compute data points from $x = 0$ to $x = 1$, and from $x = 1$ to $x = 3$
  \item (2) used least-squares fitting to get polynomials with numeric coefficients
  \item (3) used a combination of guessing and \verb+LatticeReduce+ to determine analytic expressions for the numeric coefficients.
\end{itemize}
In Section \ref{S:sinTokth}, we showed that this set of polynomials has the Fourier coefficients $\sin(n)^3/n^4$.

\ifthenelse {\boolean{BKMRK}}
  { \section{Multiplying by higher powers of \texorpdfstring { $\sin(nx)$}{sin(nx)} } \label{S:HigherPowersOfSinNx} }
  { \section{Multiplying by higher powers of $(\sin(nx))$} \label{S:HigherPowersOfSinNx} }

\begin{figure}[ht]
  \mbox{\includegraphics[width=\picSingleWidth]{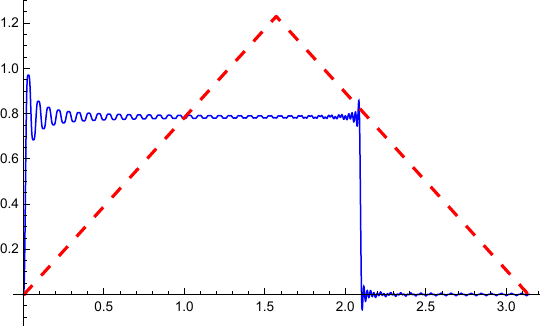}}
  \caption{$\sum_{n=1}^{100} \frac{\sin ^{3} (nx)}{n}  $ (solid) and $\sum_{n=1}^{100} \frac{\sin ^{4} (nx)}{n^{2} }$ (dashed)}
  \label{fig:hpPlot}
\end{figure}

The terms in a Fourier series involve the first power of $\sin(nx)$ and $\cos(nx)$.
But after a little experimentation, we discover that the graphs of functions like
\[
  \sum_{n=1}^{\infty} \frac{\sin ^{3} (nx)}{n}
\]
and
\[
  \sum_{n=1}^{\infty} \frac{ \sin ^{4} (nx) }{ n^{2} }
    = \sum_{n=1}^{\infty} \frac{\sin ^{3} (nx)}{n} \cdot \frac{\sin(nx)}{n}
\]
look interesting, too. See Figure \ref{fig:hpPlot}.

Notice that the graphs appear to cross at $x = 1$, suggesting that
\[
\sum_{n=1}^{\infty} \frac{\sin ^{3} (n)}{n}
  = \sum_{n=1}^{\infty} \frac{\sin ^{3} (n)}{n} \cdot \frac{\sin(n)}{n}
  = \sum_{n=1}^{\infty} \frac{\sin ^{4} (n)}{n^{2} } \, .
\]
We will prove this, and we will evaluate these and similar sums.
We will also show that there are analogues of Equation \eqref{E:Prob6241a} that involve higher powers of $\sin(n)/n$.

The identity $\cos (2nx) = 1 - 2\sin ^{2} (nx)$ suggests that if we simply replace $x$ with $2x$ in a Fourier series, we can obtain sums that involve $\sin{}^{2}(nx)$ and $\cos{}^{2}(nx)$.
In our examples, we will usually replace $x$ not with $kx$, but with $\pi - kx$, although the idea is the same.

We will use these five standard identities from trigonometry: for any $a$ and $b$,
\begin{equation} \label{E:sineSum}
  \sin(a+b) = \sin(a) \cos(b) + \cos(a) \sin(b)
\end{equation}
\begin{equation} \label{E:cosineSum}
  \cos(a+b) = \cos(a) \cos(b) - \sin(a) \sin(b)
\end{equation}
\begin{equation} \label{E:sine3a}
  \sin(3a) = 3\sin(a) - 4\sin^{3}(a)
\end{equation}
\begin{equation} \label{E:cosine2a}
  \cos(2a) = 1 - 2 \sin^{2}(a)
\end{equation}
\begin{equation} \label{E:cosine4a}
  \cos (4a) = 1 - 8\sin^2(a) + 8\sin^4(a)
\end{equation}

By performing elementary manipulations on Fourier series, we can obtain higher-power analogues of Equation \eqref{E:Prob6241a}.
The next Theorem provides more examples of series in which we can multiply the $n$\textsuperscript{th} term by $\sin(n)/n$ without changing the sum.

\begin{theorem}\label{T:thmHP}
\begin{equation}\label{E:thmHP1}
 \sum_{n=1}^{\infty} \frac{\sin ^{3} (nx)}{n} = \frac{\pi }{4}        \quad \text{ and } \quad
 \sum_{n=1}^{\infty} \frac{\sin ^{4} (nx)}{n^{2} } = \frac{\pi x}{4}
\end{equation}
for $0 < x < 2\pi /3$ and $0 \leq x \leq \pi /2$, respectively.

Also,
\begin{equation}\label{E:thmHP2}
 \sum_{n=1}^{\infty} \frac{\sin ^{5} (nx)}{n} = \frac{3\pi }{16}        \quad \text{ and } \quad
 \sum_{n=1}^{\infty} \frac{\sin ^{6} (nx)}{n^{2} } = \frac{3\pi x}{16}
\end{equation}
for $0 < x < 2\pi /5$ and $0 \leq x \leq \pi /3$, respectively.
\end{theorem}

\textbf{Proof.}
We start with this example from a textbook \cite[Equation 13.9, p.\ 32]{Tolstov}:
\begin{equation} \label{E:Tolstov32}
\sum_{n=1}^{\infty} \frac{(-1)^{n+1} }{n} \sin(nx) = \frac{x}{2} ,
\end{equation}
which is valid for $-\pi < x < \pi$.
We will replace $x$ by $\pi - 3x$ throughout \eqref{E:Tolstov32}.
The result holds provided $-\pi < \pi - 3x < \pi$, that is, for $0 < x < 2\pi /3$.

First, with this substitution, $\sin(n x)$ becomes $\sin(n(\pi - 3x)) \, .$
Using Equation \eqref{E:sineSum},
\begin{align*}
  \sin(n(\pi - 3x)) & = \sin(n \pi - 3 n x) = \sin(n \pi) \cos(3 n x) - \cos(n \pi) \sin(3 n x) \\
                    & = 0 - (-1)^n \sin(3 n x) \\
                    & = (-1)^{n+1} \sin(3 n x) \, .
\end{align*}

Using this and Equation \eqref{E:sine3a}, we have
\[
  \sin(n(\pi - 3x)) = (-1)^{n+1} (3\sin(nx) - 4\sin ^{3} (nx)) \, .
\]

Making this replacement, we get
\[
  \sum_{n=1}^{\infty} \frac{(-1)^{n+1} }{n} \sin(n(\pi - 3x))
  = 3 \sum_{n=1}^{\infty} \frac{\sin(nx)}{n} - 4\sum_{n=1}^{\infty} \frac{\sin ^{3} (nx)}{n}
  = \frac{\pi - 3x}{2} \, .
\]
Rearranging, this becomes
\begin{equation} \label{E:thmHP10}
  4 \sum_{n=1}^{\infty} \frac{\sin ^{3} (nx)}{n} = -\frac{\pi -3x}{2} + 3 \sum_{n=1}^{\infty} \frac{\sin(nx)}{n} \, .
\end{equation}

We know from Equation \eqref{E:five1} that, for $0 < x < 2\pi$,
\[
  \sum_{n=1}^{\infty} \frac{\sin(nx)}{n} = \frac{\pi - x}{2} \, .
\]
Substituting this into \eqref{E:thmHP10}, we get the first Equation in \eqref{E:thmHP1}:
\begin{equation} \label{E:thmHP20}
\sum_{n=1}^{\infty} \frac{\sin ^{3} (nx)}{n}  = \frac{\pi }{4} \, .
\end{equation}

The original substitution in Equation \eqref{E:Tolstov32} was valid if $-\pi < \pi - 3x < \pi$ (equivalently, $0 < x < 2\pi/3$).
We also used Equation \eqref{E:five1}, which is valid for $0 < x < 2\pi$.
Therefore, Equation \eqref{E:thmHP10} is valid whenever both of these conditions hold (that is, the intersection of the two intervals), which is $0 < x < 2\pi /3$.

To get the second Equation in \eqref{E:thmHP1}, we use a similar technique, starting with this easily-verified Fourier series \cite[p.\ 25]{Tolstov}, valid for $-\pi \leq x \leq \pi$ :
\begin{equation} \label{E:thmHP30}
  \sum_{n=1}^{\infty} \frac{(-1)^{n+1} }{n^{2} } \cos (nx) = \frac{ \pi^2 - 3x^2 }{12} \, .
\end{equation}

We will replace $x$ by $\pi - 2x$ in Equation \eqref{E:thmHP30}.
From Equations \eqref{E:cosineSum} and \eqref{E:cosine2a},
\[
  \cos (n\pi - 2nx) = (-1)^{n} \cos (2nx) = (-1)^{n} (1 - 2\sin ^{2} (nx)) \, .
\]

When we replace $x$ by $\pi - 2x$ throughout Equation \eqref{E:thmHP30}, we get
\[
\sum_{n=1}^{\infty} \frac{(-1)^{2n+1} }{n^{2} } (1-2\sin ^{2} (nx)) =
\sum_{n=1}^{\infty} \frac{2\sin ^{2} (nx)-1}{n^{2} } =
\frac{\pi ^{2} -3(\pi -2x)^{2} }{12} \, .
\]

Rearranging this, and using the fact that
\[
\sum_{n=1}^{\infty} \frac{1}{n^{2}} = \frac{\pi ^{2}}{6} ,
\]
we get
\begin{equation} \label{E:thmHP40} 
\sum_{n=1}^{\infty} \frac{\sin ^{2} (nx)}{n^{2} } = \frac{x(\pi -x)}{2} \, .
\end{equation}

We began with Equation \eqref{E:thmHP30}, which is valid for $-\pi \le x\le \pi $.
Therefore, \eqref{E:thmHP40} is valid for $-\pi \leq \pi -2x \leq \pi$, that is, for $0 \leq x \leq \pi$.
We have thus obtained a series involving the second power of $\sin(nx)$.
Note that Equation \eqref{E:thmHP40} generalizes the second sum in Equation \eqref{E:Prob6241a}.

We now carry this one step further.
If we replace $x$ by $\pi - 4x$ in \eqref{E:thmHP30}, we can obtain a sum involving the \emph{fourth} power of $\sin(nx)$.
From Equations \eqref{E:cosineSum} and \eqref{E:cosine4a},
\[
  \cos (n(\pi -4x)) = (-1)^{n} \cos(4nx) = (-1)^{n} (1 - 8\sin^{2}(nx) + 8\sin^{4}(nx)) \, .
\]
Replacing $x$ by $\pi - 4x$ in \eqref{E:thmHP30} and combining the powers of $(-1)$, we get
\[
  \sum_{n=1}^{\infty} \frac{-1}{n^{2} } + 8 \sum_{n=1}^{\infty} \frac{\sin ^{2} (nx)}{n^{2} } - 8 \sum_{n=1}^{\infty} \frac{\sin ^{4} (nx)}{n^{2} }
    = \frac{\pi ^{2} -3(\pi -4x)^{2} }{12} \, .
\]
This equation holds whenever $-\pi \leq \pi -4x \leq \pi$, that is, $0 \leq x \leq \pi /2$.
Solving for the third sum on the left, we get
\[
\sum_{n=1}^{\infty} \frac{\sin ^{4} (nx)}{n^{2} } = \frac{-1}{8}
  \left( \frac{\pi ^{2} -3(\pi - 4x)^{2} }{12} + \sum_{n=1}^{\infty} \frac{1}{n^{2} } - 8 \sum_{n=1}^{\infty} \frac{\sin ^{2} (nx)}{n^{2} } \right) \, .
\]

The second sum is $\pi^2/6$.
From Equation \eqref{E:thmHP40}, we know that the third sum on the right is $x(\pi - x)/2$ provided $0 \leq x \leq \pi$.
So, making these substitutions, we get half of Equation \eqref{E:thmHP1}:
\[
  \sum_{n=1}^{\infty} \frac{\sin ^{4} (nx)}{n^{2} } = \frac{\pi x}{4} \, .
\]

This result is valid where \eqref{E:thmHP30} and all intermediate steps are valid, that is, for $0 \leq x \leq \pi /2$.

Equation \eqref{E:thmHP2} can be derived in a similar way.
We replace $x$ by $\pi - 5x$ in \eqref{E:Tolstov32} and use a standard identity for $\sin(5x)$:
\[
\sin(5x) = 5 \sin(x) - 20 \sin(x)^3 + 16 \sin(x)^5 \, .
\]
Likewise, we can replace $x$ by $\pi - 6x$ in \eqref{E:thmHP30}.
After carrying out the algebra, we obtain both sums in \eqref{E:thmHP2}.

\textbf{QED.}

\textbf{Discussion. 1.}
All sums in Equations \eqref{E:thmHP1} and \eqref{E:thmHP2} are valid at $x = 1$.
Substituting $x = 1$ into \eqref{E:thmHP1} and \eqref{E:thmHP2} gives us:
\begin{equation} \label{E:SumSinCubedOverN}
 \sum_{n=1}^{\infty} \frac{\sin ^{3} (n)}{n} =
\sum_{n=1}^{\infty} \frac{\sin ^{4} (n)}{n^{2} } =
 \frac{\pi }{4}
\end{equation}
and
\begin{equation} \label{E:SumSinFifthOverN}
 \sum_{n=1}^{\infty} \frac{\sin ^{5} (n)}{n} =
 \sum_{n=1}^{\infty} \frac{\sin ^{6} (n)}{n^{2} } =
 \frac{3\pi }{16} \, .
\end{equation}

\textbf{2}.
The reader may wish to verify that, had we substituted $2x$ for $x$ in \eqref{E:thmHP30}, that instead of Equation \eqref{E:thmHP40}, we would have obtained the following, valid over $-\pi/2 \leq x \leq \pi/2$:
\begin{equation} \label{E:SumAltSinSquaredNx}
\sum_{n=1}^{\infty} \frac{(-1)^{n+1} }{n^{2} } \sin ^{2} (nx) = \frac{x^{2} }{2} \, .
\end{equation}

Equations \eqref{E:thmHP40} and \eqref{E:SumAltSinSquaredNx} are both valid over $0 \leq x \leq \pi /2$.
Therefore, if we first add \eqref{E:thmHP40} to \eqref{E:SumAltSinSquaredNx}, and then subtract \eqref{E:SumAltSinSquaredNx} from \eqref{E:thmHP40}, we get the following results, valid for $0 \leq x \leq \pi /2$:
\[\sum_{n=1}^{\infty} \frac{\sin ^{2} (2nx)}{(2n)^{2} } = \frac{x(\pi -2x)}{4} \]
\[\sum_{n=1}^{\infty} \frac{\sin ^{2} ((2n-1)x)}{(2n-1)^{2} } = \frac{\pi x}{4} \, .\]

Notice the similarity to Equations \eqref{E:Sin2nxOver2N} and \eqref{E:SinOddNxOverOddN}.

\textbf{3.}
Equation \eqref{E:thmHP1} shows that
\[
\sum_{n=1}^{\infty} \frac{\sin ^{3} (nx)}{n} = \frac{\pi}{4}
\]
for $0 < x < 2\pi/3$.
The graph in Figure \ref{fig:hpPlot} suggests that
\[
\sum_{n=1}^{\infty} \frac{\sin ^{3} (nx)}{n} = 0
\]
for $2\pi/3 < x \leq \pi$.
This can be proved the same way we proved the first half of \eqref{E:thmHP1}.
But instead of starting with \eqref{E:Tolstov32}, we start with the fact that, for $-3\pi < x < -\pi$, \eqref{E:Tolstov32} has this slightly different form:
\[
\sum_{n=1}^{\infty} \frac{(-1)^{n+1} }{n} \sin(nx) = \frac{x}{2} + \pi \, .
\]
Then, replace $x$ with $\pi - 3x$ and proceed as before.

\textbf{4.}
Figure \ref{fig:hpPlot} shows the series
\[
\sum_{n=1}^{\infty} \frac{\sin ^{3} (nx)}{n} \, .
\]
This series represents the function
\[
f(x)=
  \begin{cases}
    \pi/4  &\text{for $0 < x < 2\pi/3$,}\\
    0      &\text{for $2\pi/3 < x < \pi$.}
  \end{cases}
\]
provided we extend $f(x)$ to $x < 0$ by taking its odd periodic extension.
The above series, as written, is not a Fourier series, but this $f(x)$ does have this sine series over $0 \leq x \leq \pi $:
\[
\sum_{n=1}^{\infty} \frac{\sin ^{2} (n\pi /3)}{n} \sin(nx) = 
\frac{3}{4} \left(\frac{\sin x}{1} + \frac{\sin(2x)}{2} + \frac{\sin(4x)}{4} + \frac{\sin(5x)}{5} + \frac{\sin(7x)}{7} + \frac{\sin(8x)}{8} + \dots \right) \, ,
\]
which skips the terms whose denominators are multiples of 3.

\begin{figure}[ht]
  \mbox{\includegraphics[width=\picSingleWidth]{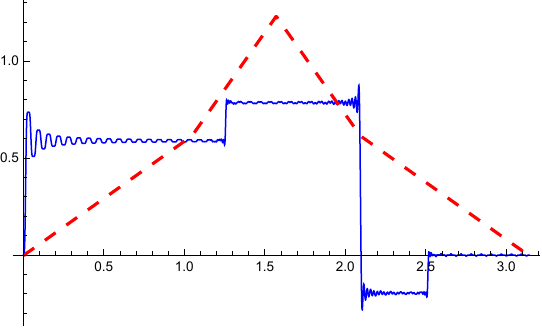}}
  \caption{$\sum_{n=1}^{100} \frac{\sin ^{5} (nx)}{n}  $ (solid) and $\sum_{n=1}^{100} \frac{\sin ^{6} (nx)}{n^{2} }$ (dashed)}
  \label{fig:hpPlot56}
\end{figure}

\textbf{5.}
Figure \ref{fig:hpPlot56} shows the graphs of the sums
\[
\sum_{n=1}^{\infty} \frac{\sin ^{5} (nx)}{n}
\]
and
\[
\sum_{n=1}^{\infty} \frac{\sin ^{6} (nx)}{n^2}
\]
for $0 \leq x \leq \pi$.
Equation \eqref{E:thmHP2} tells us that the functions are equal at $x = 1$.
The first sum takes the constant $y$ values $3\pi/16$, $\pi/4$, $-\pi/16$, and 0, with jumps occurring at $x = 0$, $2\pi/5$, $2\pi/3$, $4\pi/5$, and $\pi$.
Just as in discussion point 3, these $x$ and $y$ values can be calculated by considering what \eqref{E:Tolstov32} looks like outside of $-\pi < x < \pi$.

The reader may wish to verify (with computer assistance) that, for $-\pi < x < \pi$, the sine series for $\sum \sin ^{5} (nx)/n $ is
\[
\frac{1}{8} \sum_{n=1}^{\infty} \frac{3 + \cos (2n\pi /5) + \cos (4n\pi /5) - 5\cos (2n\pi /3)}{n} \sin(nx) \, .
\]
(For $n = 1$ through $n = 15$, the expression $3 + \cos (2n\pi /5) + \cos (4n\pi /5) - 5\cos (2n\pi /3)$ takes the values
$5, 5, -5/2, 5, 15/2, -5/2, 5, 5, -5/2, 15/2, 5, -5/2, 5, 5, 0$. For larger $n$, this sequence repeats.)
  
\begin{figure}[ht]
  \mbox{\includegraphics[width=\picSingleWidth]{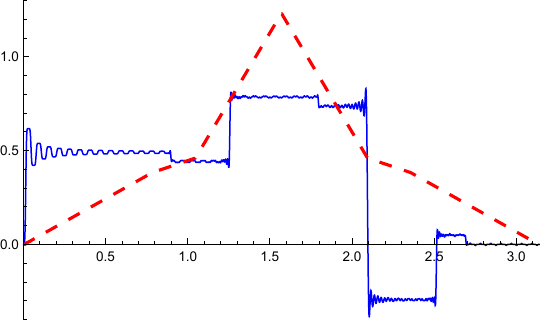}}
  \caption{$\sum_{n=1}^{100} \frac{\sin ^{7} (nx)}{n}  $ (solid) and $\sum_{n=1}^{100} \frac{\sin ^{8} (nx)}{n^{2} }$ (dashed)}
  \label{fig:plot78}
\end{figure}

\textbf{6.}
Alas, a computer calculation can verify that equations corresponding to \eqref{E:SumSinCubedOverN} and \eqref{E:SumAltSinSquaredNx}, but with 7\textsuperscript{th} and 8\textsuperscript{th} powers, does \emph{not} hold.
In fact, \textit{Mathematica} claims the sums are
\[
\sum_{n=1}^{\infty} \frac{\sin^7(n)}{n} = \frac{9\pi}{64}
\]
and
\[
\sum_{n=1}^{\infty} \frac{\sin^8(n)}{n^2} = \frac{\pi(6 + \pi)}{64} \, .
\]
Figure \ref{fig:plot78} shows the graphs of the sums of $\sin^7(nx)/n $ and $\sin^8(nx)/n^2$ \, .
One must examine them very closely to see that they cross, not at $x = 1$, but at $x \approx 0.976$.
This illustrates why one must use graphs only as guidelines, not as proofs.

\section{Using polylogarithms} \label{S:UsingPolylogs}

We can express many of the sums we've previously considered in terms of polylogarithms, which were mentioned above in Section \ref{S:sinTokth}.
The polylogarithm $\text{Li}_k(z)$ can be defined by the series
\begin{equation} \label{E:PolylogSeriesK}
\text{Li}_k(z) = \sum_{n=1}^{\infty} \frac{z^n}{n^k} \, .
\end{equation}
If $k \geq 2$ is an integer, this converges for $\lvert z \rvert \leq 1$.
With analytic continuation, polylogarithms can be extended beyond this region.
For more details, see these books on polylogarithms: \cite{Lewin} or \cite{Srivastava}.

\subsection{Dilogarithms} \label{S:Dilogarithms}
$\text{Li}_2(z)$ is called the \textit{dilogarithm}.

Here, we will show how to use dilogarithms to compute the sum in Equation \eqref{E:SumSincSquared}:
\[
  \sum_{n=1}^{\infty} \left( \frac{\sin(n)}{n} \right)^2 \, .
\]

An equation of Euler says that, for all $n$ (here, $n$ can be \textit{any} real or complex number),
\[
  e^{in} = \cos(n) + i \sin(n) \, .
\]
Likewise, $e^{-in} = \cos(n) - i \sin(n)$.
Subtracting these equations and dividing by $2i$, we get
\begin{equation} \label{E:EulerFormForSin}
  \sin(n) = \frac{e^{in} - e^{-in}}{2i} =
  \frac{i \left( e^{in} - e^{-in} \right) }{2i^2} =
  \frac{i}{2} (e^{-in} - e^{in} ) \, .
\end{equation}
Likewise, for any $n$ and any $x$,
\begin{equation} \label{E:EulerFormForSinx}
  \sin(n x) = \frac{i}{2} (e^{-i n x} - e^{i n x} ) \, .
\end{equation}

From Equation \eqref{E:EulerFormForSin}, we have
\begin{equation} \label{E:EulerFormForSinSquared}
\sin^2(n) = \left( \frac{i}{2} \right)^2 (e^{-in} - e^{in} )^2
= \frac{-1}{4} \left(  e^{-2in} + e^{2in} - 2 \right) \, .
\end{equation}
Dividing by $n^2$ and summing over $n$,
\[
\sum_{n=1}^{\infty} \left( \frac{\sin(n)}{n} \right)^2
= \frac{-1}{4} \left(
 \sum_{n=1}^{\infty} \frac{(e^{-2i})^n}{n^2}
+ \sum_{n=1}^{\infty} \frac{(e^{2i})^n}{n^2}
- 2 \cdot \sum_{n=1}^{\infty} \frac{1}{n^2}
\right) \, .
\]
Expressing this in terms of dilogarithms, this is
\begin{equation} \label{E:SumWithDiLog}
\sum_{n=1}^{\infty} \left( \frac{\sin(n)}{n} \right)^2
= \frac{-1}{4} \left(  \text{Li}_{2} (e^{-2i}) + \text{Li}_{2} (e^{2i}) - \frac{\pi^2}{3} \right) \, .
\end{equation}

This combination of dilogarithms can be evaluated with the ``inversion'' identity:
\begin{equation*}
\text{Li}_{2} (z) + \text{Li}_{2} (1/z) = -\frac{\pi^2}{6} - \frac{1}{2} (\ln(-z))^2  \, ,
\end{equation*}
According to \cite[Eq. 25.12.4, Chapt. 25]{DLMF}, this holds unless $z$ is a real number $\geq 0$.
(This seems to hold unless $z$ is real and $0 \leq z < 1$.)

It is slightly more convenient to use this identity in the form
\begin{equation} \label{E:Li2Inversion}
\text{Li}_{2} (-z) + \text{Li}_{2} (-1/z) = -\frac{\pi^2}{6} - \frac{1}{2} (\ln(z))^2  \, ,
\end{equation}
which holds unless $z$ is real and $-1 < z \leq 0$; however, we need this only for \emph{complex} $z$.

In Equation \eqref{E:SumWithDiLog}, we want to evaluate $\text{Li}_{2} (e^{-2i}) + \text{Li}_{2} (e^{2i})$.
So, we will use Equation \eqref{E:Li2Inversion} with $-z = e^{-2i}$.
Then $z = (-1)e^{-2i} = e^{\pi i} e^{-2i} = e^{i(\pi - 2)}$.

If $-\pi < x \leq \pi$, then $\ln(e^{i x}) = i x$.
Here, $-\pi < \pi - 2 \leq \pi$, so if $z = e^{i(\pi - 2)}$, then $\ln(z) = i(\pi - 2)$.
Then the right side of \eqref{E:Li2Inversion} is
\begin{equation}  \label{E:dilogRightSide}
-\frac{\pi^2}{6} - \frac{1}{2} \, i^2 \, (\pi - 2)^2 =
-\frac{\pi^2}{6} + \frac{1}{2}(\pi^2 - 4 \pi + 4) = \frac{\pi^2}{3} - 2 \pi + 2 \, .
\end{equation}
Equation \eqref{E:SumWithDiLog} then becomes
\[
\sum_{n=1}^{\infty} \left( \frac{\sin(n)}{n} \right)^2 =
\frac{-1}{4} \left( \frac{\pi^2}{3} - 2 \pi + 2 - \frac{\pi^2}{3} \right) = \frac{\pi - 1}{2} \, .
\]
This confirms Equation \eqref{E:SumSincSquared}.
Note that we computed the sum \emph{directly}, that is, without guesswork or having to analyze graphs.

\textbf{Logs of complex numbers.}

Shortly, we will be computing expressions like $\ln(e^{i x})$ where $x$ is a real variable.

$\ln(e^{i x})$ is not always equal to $i x$.
For example, $e^{3 \pi i} = -1$, but $\ln(e^{3 \pi i}) = \pi i$.
This happens because $e^{3 \pi i + 2 \pi i k} = -1$ for \emph{any} integer $k$.
So, the proper use of $\ln( \cdot )$ requires that we select the correct value for the logarithm.
Books such as \cite[Chapt. 6]{Danese} or \cite[Sect. 8.7]{Kaplan} that discuss complex variables, cover this in more detail.
In any case, it is helpful to keep in mind that, if $-\pi < x \leq \pi$, then $\ln(e^{i x}) = i x$.


\textbf{Sums that involve $x$.}

Now, let's work with a series that is not constant, but which is a function of $x$.
Here is Equation \eqref{E:altSignG2}, repeated here for convenience:
\begin{equation} \label{E:altSignG2-copy}
g(x) = \sum_{n=1}^{\infty} (-1)^{n+1} \frac{\sin(nx)}{n} \cdot \frac{\sin(n)}{n} =
  \begin{dcases}  
    \frac{x}{2}                              &\text{ for $0 \leq x \leq \pi-1$} \, , \\
    \frac{x}{2} + \frac{\pi}{2}(\pi - 1 - x) &\text{ for $\pi - 1 \leq x \leq \pi + 1$} \, .
  \end{dcases}
\end{equation}
The series in Equation \eqref{E:altSignG2-copy} is represented by two different polynomials over the intervals $0 \leq x \leq \pi - 1$ and $\pi - 1 \leq x \leq \pi + 1$
(note that at $x = \pi - 1$, both polynomials equal $(\pi - 1)/2$).
These polynomials were arrived at from guesswork based on the graph of $g(x)$.
Here, we will use polylogarithms to \emph{derive} these polynomials directly from the series: no guesswork needed!

Using Equations \eqref{E:EulerFormForSin} and \eqref{E:EulerFormForSinx},
we can express the sum in Equation \eqref{E:altSignG2-copy} using polylogarithms:
\begin{align*}
g(x) =
& \sum_{n = 1}^{\infty} (-1)^{n+1} \; \frac{\sin(n x)}{n} \cdot \frac{\sin(n)}{n}
 = \sum_{n = 1}^{\infty} \frac{(-1)^{n+1}}{n^2} \left( \frac{i}{2} \right)^2 (e^{-i n x} - e^{i n x} )(e^{-in} - e^{in} ) \\
& = \frac{1}{4} \sum_{n = 1}^{\infty} \frac{(-1)^{n}}{n^2}
  \left( e^{-i n - i n x} + e^{i n + i n x} - e^{i n - i n x} - e^{-i n + i n x} \right) \\
& = \frac{1}{4} \sum_{n = 1}^{\infty} \frac{(-1)^{n}}{n^2}
  \left( (e^{-i - i x})^n + (e^{i + i x})^n - (e^{i - i x})^n - (e^{-i + i x})^n \right) \\
& = \frac{1}{4} \sum_{n = 1}^{\infty} \frac{(-1)^{n}}{n^2}
  \left( (e^{-i(x + 1)})^n + (e^{i(x + 1)})^n - (e^{-i(x - 1)})^n - (e^{i(x - 1)})^n \right) \\
& = \frac{1}{4} \sum_{n = 1}^{\infty} \frac{1}{n^2}
  \left( (-e^{-i(x + 1)})^n + (-e^{i(x + 1)})^n - \left( (-e^{-i(x - 1)})^n + (-e^{i(x - 1)})^n \right) \right) \, .
\end{align*}
Writing this as four $\text{Li}_2( \cdot )$'s, this is:
\[
\frac{1}{4} \left( \text{Li}_2\left(-e^{-i (x+1)}\right) + \text{Li}_2\left(-e^{i (x+1)} \right)
          - \left( \text{Li}_2\left(-e^{-i (x-1)}\right) + \text{Li}_2\left(-e^{i (x-1)} \right) \right) 
           \right) \, .
\]
It will be a little easier if we put the ``positive'' exponents first:
\begin{equation} \label{E:gxDilogs}
g(x) = \frac{1}{4} \left( \text{Li}_2\left(-e^{i (x+1)}\right) + \text{Li}_2\left(-e^{-i (x+1)} \right)
          - \left( \text{Li}_2\left(-e^{i (x-1)}\right) + \text{Li}_2\left(-e^{-i (x-1)} \right) \right) 
           \right) \, .
\end{equation}
It should be emphasized that this expression for $g(x)$ agrees with the series in \eqref{E:altSignG2-copy} \emph{for all} $x$.
But each polynomial in \eqref{E:altSignG2-copy} agrees with the series only for a finite interval of $x$ values.

To convert \eqref{E:gxDilogs} to polynomials, we will use Equation \eqref{E:Li2Inversion} twice, with $z = e^{i(x + 1)}$ and $z = e^{i(x - 1)}$.

First, we need to compute the logs in \eqref{E:Li2Inversion}.
Let's start with an easy calculation.
If $x$ is real, what is $\ln(e^{ix})$?
The value depends on $x$:
\begin{align}
 \ln( e^{i x} ) & = i x + 2 \pi i             \quad \text{for $-3\pi  < x \leq  -\pi$}             \, , \label{E:SimpleLogA} \\
 \ln( e^{i x} ) & = i x  \phantom{+ 2 \pi i}  \quad \text{ for $ \phantom{3}-\pi  < x \leq   \pi$} \, , \label{E:SimpleLog} \\
 \ln( e^{i x} ) & = i x - 2 \pi i             \quad \text{for $  \phantom{-3}\pi  < x \leq  3\pi$} \, . \label{E:SimpleLogB}
\end{align}

Replace $x$ with $x - 1$ in Equation \eqref{E:SimpleLog}.
This gives $\ln( e^{i (x - 1)} ) = i (x - 1)$, provided $-\pi < x - 1 \leq \pi$, that is, if $-\pi + 1 < x \leq \pi + 1$.

Similarly, replacing $x$ with $x + 1$, $\ln( e^{i (x + 1)} ) = i (x + 1)$ if $-\pi < x + 1 \leq \pi$,
that is, if $-\pi - 1 < x \leq \pi - 1$.

The intersection of these two intervals
\begin{align*}
-\pi + 1 & < x \leq \pi + 1 \quad \text{and} \\
-\pi - 1 & < x \leq \pi - 1
\end{align*}
is the interval $-\pi + 1 < x \leq \pi - 1$.
For any $x$ in this interval, we will have \emph{both} $\ln( e^{i (x - 1)} ) = i (x - 1)$
and $\ln( e^{i (x + 1)} ) = i (x + 1)$.

\textbf{A. The interval $-\pi + 1 < x \leq \pi - 1$}.

For this interval, we use $z = e^{i(x + 1)}$ in Equation \eqref{E:Li2Inversion} to evaluate the first two $\text{Li}_2( \cdot )$'s in \eqref{E:gxDilogs}:
\begin{align} \label{E:Li2-FirstTwoTerms}
\text{Li}_2 \left( -e^{i (x+1)} \right) + \text{Li}_2 \left( -e^{-i (x+1)} \right)
& = -\frac{\pi^2}{6} - \frac{1}{2} \left( \ln(e^{i(x + 1)}) \right)^2
= -\frac{\pi^2}{6} - \frac{1}{2} (i(x + 1))^2 \notag \\
& = \frac{x^2}{2}+x-\frac{\pi ^2}{6}+\frac{1}{2} \, .
\end{align}

Next, we use $z = e^{i(x - 1)}$ to evaluate the last two $\text{Li}_2( \cdot )$ terms in \eqref{E:gxDilogs}, to get
\begin{align} \label{E:Li2-LastTwoTerms}
\text{Li}_2 \left( -e^{i (x-1)}\right) + \text{Li}_2 \left( -e^{-i (x-1)} \right)
& = -\frac{\pi^2}{6} - \frac{1}{2} \left( \ln(e^{i(x - 1)}) \right)^2
= -\frac{\pi^2}{6} - \frac{1}{2} (i(x - 1))^2 \notag \\
& = \frac{x^2}{2}-x-\frac{\pi ^2}{6}+\frac{1}{2} \, .
\end{align}
Subtracting the result in \eqref{E:Li2-LastTwoTerms} from that in \eqref{E:Li2-FirstTwoTerms}, the result is $2 x$, so \eqref{E:gxDilogs} becomes
\[
g(x) = \frac{1}{4}(2x) = \frac{x}{2} \, .
\]
This holds for $-\pi + 1 < x \leq \pi - 1$.

\textbf{B. The interval $\pi - 1 < x \leq \pi + 1$}.

Now, we will get an expression for $g(x)$ over the interval $\pi - 1 < x \leq \pi + 1$.
For this interval, $\ln( e^{i (x - 1)} )$ is still $i (x - 1)$, so \eqref{E:Li2-LastTwoTerms} is unchanged.

However, we now have $\pi < x + 1 \leq \pi + 2 < 3 \pi$, so $\ln( e^{i (x + 1)} ) = i (x + 1) - 2 \pi i$.
Therefore, Equation \eqref{E:Li2-FirstTwoTerms} becomes
\begin{align} \label{E:Li2-FirstTwoTerms2PiI}
\text{Li}_2 \left( -e^{i (x+1)} \right) + \text{Li}_2 \left( -e^{-i (x+1)} \right)
& = -\frac{\pi^2}{6} - \frac{1}{2} \left( \ln(e^{i(x + 1)}) \right)^2
= -\frac{\pi^2}{6} - \frac{1}{2} (i(x + 1) - 2 \pi i)^2 \notag \\
& = \frac{x^2}{2}-2 \pi  x+x+\frac{11 \pi ^2}{6}-2 \pi +\frac{1}{2} \, .
\end{align}
Subtracting the result in \eqref{E:Li2-LastTwoTerms} from that in \eqref{E:Li2-FirstTwoTerms2PiI}, the result is
\[
2 x - 2 \pi x + 2 \pi^2 - 2 \pi = 2 x + 2 \pi (\pi - 1 - x) \, .
\]
Finally,
\[
g(x) = \frac{1}{4} \left( 2 x + 2 \pi (\pi - 1 - x) \right) = \frac{x}{2} + \frac{\pi}{2} (\pi - 1 - x) \, .
\]
This is the expression for the series \eqref{E:altSignG2-copy} for $g(x)$ over the interval $\pi - 1 < x \leq \pi + 1$.

As we have seen, there are different intervals which have different polynomial expressions.
This is caused by one or more $\ln( \cdot )$ values requiring a different multiple of $2 \pi i$, as in Equation \eqref{E:Li2-FirstTwoTerms2PiI}.

\textbf{Summary.}

To summarize what we have done:
\begin{itemize}
 \item We converted the Fourier series to  an expression involving polylogarithms;
 \item Based on the interval $a < x < b$, we determined $k$ such that $\ln(e^{i(x + c)}) = i(x + c) + 2 \pi i k$;
 \item We used these logarithms in the inversion formula to produce a polynomial that equals the series over $a < x < b$;
 \item Each interval having a different $k$ produces a different polynomial expression for the series.
\end{itemize}


Equation \eqref{E:altSignG2-copy} gives polynomial expressions for $g(x)$ over an interval \emph{larger} than $\pi$.
This series involves only sines.
Therefore, if we know $g(x)$  only for $0 < x < \pi$, then:
(a) $\sin(-x) = -sin(x)$, so we can also compute $g(x)$ for $-\pi < x < 0$;
(b) $\sin(x) = -\sin(2 \pi - x)$, so we can also compute $g(x)$ for $\pi < x < 2 \pi$.

These observations will may some of the calculations in the next few sections unnecessary.
However, it is worth going through these calculations to see how
polylogarithms produce different polynomials over different intervals.

\subsection{Trilogarithms}

$\text{Li}_3(z)$ is called the \textit{trilogarithm}.

This inversion formula lets us evaluate sums involving trilogarithms \cite[Eq. A.2.6(5), p.\ 296]{Lewin}:
\begin{equation} \label{E:Li3Inversion}
\text{Li}_{3} (-z) - \text{Li}_{3} (-1/z) = -\frac{\pi^{2} }{6} \ln (z) - \frac{1}{6} \ln^{3} (z) \, .
\end{equation}

This holds unless $z$ is a real number, $-1 < z \leq 0$.
This inversion formula can also be found as Equation 62 in \cite[p.\ 113]{Srivastava}.

We will use trilogarithms to evaluate the sum
\begin{equation}  \label{E:trilogExprSum}
\sum_{n=1}^{\infty} \frac{\sin(n)}{n^3} = \frac{i}{2} \cdot
 \left(
    \sum_{n=1}^{\infty} \frac{(e^{-i})^n}{n^3} - \sum_{n=1}^{\infty} \frac{(e^{i})^n}{n^3}
 \right)
= \frac{i}{2} \cdot \left( \text{Li}_{3} (e^{-i}) - \text{Li}_{3} (e^{i}) \right)
\, .
\end{equation}

We use Equation \eqref{E:Li3Inversion} with $-z = e^{-i}$ and $-1/z = e^i$.
With this $z$, we have $z = (-1)e^{-i} = e^{\pi i} \cdot e^{-i} = e^{i(\pi - 1)}$.
For this $z$, we confirm that $\ln(z) = i(\pi - 1)$.
Then the sum is (note: $i^3 = -i$)
\begin{align} \label{E:trilogExpr}
\sum_{n=1}^{\infty} \frac{\sin(n)}{n^{3} }
  & = \frac{i}{2} \left(
         -\frac{\pi ^{2} }{6} i(\pi - 1) - \frac{1}{6} ( i(\pi - 1) )^3
       \right)  \notag \\
  & = \frac{i}{2} \left( (-i) \frac{\pi^3 - \pi^2}{6} - (-i) \frac{\pi^3 - 3 \pi^2 + 3 \pi - 1}{6} \right)
   = \frac{\pi ^{2} }{6} - \frac{\pi }{4} + \frac{1}{12} \, .
\end{align}

\subsection{Another example} \label{S:AnotherTrilogExample}  


Here, we present a simple proof that
\begin{equation} \label{E:ThirdSumOneHalf}
\sum_{n = 1}^{\infty} (-1)^{n+1} \left( \frac{\sin(n)}{n} \right)^3 = \frac{1}{2} \, .
\end{equation}

Because $\sin(n) = (e^{i n} - e^{-i n})/(2 i)$, we can write this sum as
\begin{equation*}  
   \sum_{n=1}^{\infty} (-1)^{n+1} \left( \frac{\sin(n)}{n} \right)^3
 = \sum_{n=1}^{\infty} (-1)^{n+1} \frac{1}{n^3} \left( \frac{1}{2 i} \cdot ( e^{i n} - e^{-i n} ) \right)^3
 = \frac{i}{8} \sum_{n=1}^{\infty} (-1)^n \frac{1}{n^3} ( e^{-i n} - e^{i n} )^3 \, ,
\end{equation*}
where we simplified $1/(2 i)^3 = i/8$.

The product of the exponentials in the numerator of the $n$\textsuperscript{th} term is
\begin{equation*}
  ( e^{-i n} - e^{i n} )^3 = e^{-3 i n} - 3 \left( e^{-i n} - e^{i n} \right) - e^{3 i n}
  = (e^{-3 i})^n - 3 \left( (e^{- i})^n - (e^{i})^n \right) - (e^{3 i})^n \, .
\end{equation*}
Multiplying this last line through by $(-1)^n$, the numerator of the $n$\textsuperscript{th} term is
\begin{equation*}
(-e^{-3 i})^n  - (-e^{3 i})^n - 3 \left( (-e^{- i})^n - (-e^{i})^n \right) \, .
\end{equation*}
So, the sum in Equation \eqref{E:ThirdSumOneHalf} is
\begin{equation*}
\frac{i}{8} \big( 
  \sum_{n=1}^{\infty} \frac{ (-e^{-3 i})^n }{n^3} - \sum_{n=1}^{\infty} \frac{(-e^{3 i})^n}{n^3}
  - 3 \big( \sum_{n=1}^{\infty} \frac{(-e^{- i})^n}{n^3} - \sum_{n=1}^{\infty} \frac{(-e^{i})^n}{n^3} \big)
    \big) \, .
\end{equation*}

Therefore, we can write
\begin{equation} \label{E:thirdSumPolyLogSimple} 
\sum_{n=1}^{\infty} (-1)^{n+1} \left( \frac{\sin(n)}{n} \right)^3
= \frac{i}{8} \;
   \big\{ \;
      \text{Li}_3( -e^{-3 i} ) - \text{Li}_3( -e^{3 i} )
     - 3 \big( \text{Li}_3( -e^{-i} ) - \text{Li}_3( -e^{i} ) \big)
      \big \} \, .
\end{equation}

We use the inversion formula, Equation \eqref{E:Li3Inversion}.
We will apply this formula twice, with $z = e^{-3i}$, then with $z = e^{-i}$.
The logarithms needed in Equation \eqref{E:Li3Inversion} are $\ln( e^{-3 i} ) = -3 i$ and $\ln( e^{-i} ) = -i$.

If we use these logarithms and Equation \eqref{E:Li3Inversion}, then Equation \eqref{E:thirdSumPolyLogSimple} becomes
\begin{align*} 
\sum_{n=1}^{\infty} (-1)^{n+1} \left( \frac{\sin(n)}{n} \right)^3
= & \frac{i}{8} \;
   \big\{ \;
      \text{Li}_3( -e^{-3 i} ) - \text{Li}_3( -e^{3 i} )          - 3 \big( \text{Li}_3( -e^{-i} ) - \text{Li}_3( -e^{i} ) \big)
   \big \} \\
= & \frac{i}{8} \;
   \big\{ \;
         -\frac{\pi^{2} }{6} \cdot (-3 i) - \frac{1}{6} (-3 i)^3  - 3 \big(  -\frac{\pi^{2} }{6} \cdot (-i) - \frac{1}{6} (-i)^3  \big)
   \big \} \\
= & \frac{i}{8} \;
   \big\{ \;
         \frac{3 i \pi^{2} }{6} - \frac{27 i}{6} - 3 \big(  \frac{i \pi^{2} }{6} - \frac{i}{6} \big)
   \big \} \\
      = & \frac{1}{2} \, .
\end{align*}

\subsection{Using polylogarithms to evaluate more sums involving sines} \label{S:PolylogApplicationX}


\mbox{}  

In previous sections, we were able to evaluate sums in terms of polynomials.
However, it is difficult to fit this to a polynomial:
\[
  \sum_{n=1}^{\infty} (-1)^{n+1} \frac{\sin(n x)}{n} \left( \frac{\sin(n)}{n} \right) ^2 \, .
\]
From $x = 0$ to $x = 1$, we get a good fit of this function to $x/2$.
Unfortunately, we cannot get a good fit of this function to a polynomial over the interval from $x = 1$ to $x = 2$.
As we will discover below, that is because this function equals $x/2$ for $0 \leq x \leq \pi - 2$ $(\approx 1.14)$, but equals a different polynomial for $x > \pi - 2$.

However, we will see that polylogarithms allow us to evaluate this sum.

Theorem \ref{T:TwoAlternatingSums} stated that, for $0 \leq x < \pi$,
\[
\sum_{n=1}^{\infty} (-1)^{n+1} \frac{\sin(nx)}{n} = \frac{x}{2} \, ,
\]
and that, for $0 \leq x \leq \pi - 1$ (in fact, for $-\pi + 1 < x \leq \pi - 1$ (see Section \ref{S:Dilogarithms})),
\[
\sum_{n=1}^{\infty} (-1)^{n+1} \frac{\sin(nx)}{n} \cdot \frac{\sin(n)}{n} = \frac{x}{2} \, .
\]
We will now use polylogarithms to prove that, for $0 \leq x \leq \pi - 2$, this sum is \emph{also} $x/2$:
\begin{equation} \label{E:alternatingSum3x}
  \sum_{n=1}^{\infty} (-1)^{n+1} \frac{\sin(nx)}{n} \cdot \left( \frac{\sin(n)}{n} \right) ^2 = \frac{x}{2} \, .
\end{equation}

\begin{figure}[ht]
  \mbox{\includegraphics[width=\picSingleWidth]{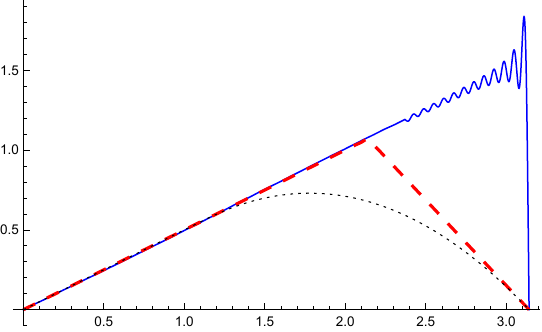}}
  \caption{$\sum_{n=1}^{100} (-1)^{n+1} \frac{\sin(nx)}{n} \left( \frac{\sin(n)}{n} \right)^k$ for $k = 0, 1, 2$}
  \label{fig:threeSumsAltSigns2}
\end{figure}

Figure \ref{fig:threeSumsAltSigns2} shows the graphs of
\[
  \sum_{n=1}^{\infty} (-1)^{n+1} \frac{\sin(nx)}{n} \cdot \left( \frac{\sin(n)}{n} \right)^k \, .
\]
$k = 0$ is the solid graph, $k = 1$ is the dashed graph, and $k = 2$ is the dotted graph.
All three curves appear to equal $x/2$ over some interval from $x = 0$ to somewhere around $x = 1$.

Using Equations \eqref{E:EulerFormForSinx} and \eqref{E:EulerFormForSin}, we can write the left side of Equation \eqref{E:alternatingSum3x} as
\begin{equation} \label{E:AltSum3x10}
  \sum_{n=1}^{\infty} (-1)^{n+1} \frac{\sin(nx)}{n} \cdot \left( \frac{\sin(n)}{n} \right) ^2
  = \frac{i}{8} \sum_{n=1}^{\infty} \frac{1}{n^3} (-1)^n ( e^{-i n x} - e^{i n x} )( e^{-i n} - e^{i n} )^2 \, .
\end{equation}
(Here, we have removed one factor of $(-1)$, leaving $(-1)^n$ inside the sum; also,  $(-1)(i/2)^3 = i/8$.)
The product of the exponentials in the numerator of the $n$\textsuperscript{th} term is
\begin{align*}
  & (-1)^n ( e^{-i n x} - e^{i n x} )( e^{-i n} - e^{i n} )^2 \\
  = & (-1)^n \left( e^{2 i n-i n x} - e^{i n x-2 i n} - 2 ( e^{-i n x} - e^{i n x} ) + e^{-i n x-2 i n} - e^{i n x+2 i n} \right) \\
  = & (-1)^n \left( \, (e^{-i(x - 2)})^n - (e^{i(x - 2)})^ n - 2 ( \, (e^{-i x})^n - (e^{i x})^n \, ) + (e^{-i( x+2)} )^n - (e^{i(x+2)} )^n \, \right) \, .
\end{align*}
Multiplying this last line through by $(-1)^n$, the numerator of the $n$\textsuperscript{th} term is
\begin{equation} \label{E:MultThrough} 
(-e^{-i(x - 2)})^n - (-e^{i(x - 2)})^n \; - \;  2  ( (-e^{-i x})^n - (-e^{i x})^n )  \;  + \;  (-e^{-i( x+2)} )^n - (-e^{i(x+2)} )^n \, .
\end{equation}

We can now separate the terms in the complicated numerator and write our expression in terms of polylogarithms.
The right side of Equation \eqref{E:AltSum3x10} is
\begin{align*}
\frac{i}{8} & \big( 
  \sum_{n=1}^{\infty} \frac{(-e^{-i(x - 2)})^n}{n^3} - \sum_{n=1}^{\infty} \frac{(-e^{i(x - 2)})^n}{n^3} \\
  - & 2 \big( \sum_{n=1}^{\infty} \frac{(-e^{-ix})^n}{n^3} - \sum_{n=1}^{\infty} \frac{(-e^{ix})^n}{n^3} \big) \\
   + & \sum_{n=1}^{\infty} \frac{(-e^{-i(x + 2)})^n}{n^3} - \sum_{n=1}^{\infty} \frac{(-e^{i(x + 2)})^n}{n^3}
    \big) \, .
\end{align*}
Finally, we have
\begin{align} \label{E:AltSum3x20}
& \sum_{n=1}^{\infty} (-1)^{n+1} \frac{\sin(nx)}{n} \cdot \left( \frac{\sin(n)}{n} \right) ^2 = \notag  \\
& \frac{i}{8} \;
   \big\{ \; 
      \text{Li}_3( -e^{-i(x-2)} ) - \text{Li}_3( -e^{i (x-2)} ) \notag  \\    
    & - 2 \big( \text{Li}_3( -e^{-i x} ) - \text{Li}_3( -e^{i x} ) \big) \notag  \\
    & + \text{Li}_3( -e^{-i (x+2)} ) - \text{Li}_3( -e^{i (x+2)} ) \;
      \big \} \, .
\end{align}
This holds for \emph{all} $x$, not merely for those $x$ in a short interval.


To evaluate these $\text{Li}_3( \cdot )$ expressions, we will use Equation \eqref{E:Li3Inversion}, which states that
\[
\text{Li}_{3} (-z) - \text{Li}_{3} (-1/z) = -\frac{\pi ^{2} }{6} \ln (z) - \frac{1}{6} \ln ^{3} (z) \, .
\]
We will apply this inversion identity three times, with $z = e^{-i(x-2)}$, with $z = e^{-ix}$, and with $z = e^{-i(x+2)}$.
Then the three logarithms we need are
\begin{align}
 \ln( e^{-i x} )       & = -i x \, , \label{E:AltSum3xSimpleLog1} \\
 \ln( e^{-i (x - 2)} ) & = -i (x - 2) \, , \label{E:AltSum3xSimpleLog2} \\
 \ln( e^{-i (x + 2)} ) & = -i (x + 2) \, . \label{E:AltSum3xSimpleLog3}
\end{align}


From Equation \eqref{E:SimpleLog}, $\ln( e^{i x} ) = i x$ if $-\pi < x \leq \pi$.
Therefore, \eqref{E:AltSum3xSimpleLog1} is valid if $-\pi < -x \leq \pi$, which is the same as $-\pi \leq x < \pi$.
So, Equation \eqref{E:AltSum3xSimpleLog1} holds for $-\pi \leq x < \pi$.

Replacing $x$ with $x - 2$, we see that \eqref{E:AltSum3xSimpleLog2} is valid for $-\pi \leq x - 2 < \pi$, that is, for $-\pi + 2 \leq x < \pi + 2$.
Similarly, \eqref{E:AltSum3xSimpleLog3} is valid for $-\pi \leq x + 2 < \pi$, that is, $-\pi - 2 \leq x < \pi - 2$.

The intersection of these three intervals is $-\pi + 2 \leq x < \pi - 2$.
Over this interval, all three logarithms in \eqref{E:AltSum3xSimpleLog1} through \eqref{E:AltSum3xSimpleLog3} hold without any extra $\pm 2 \pi i$.

For $x \geq \pi - 2$, one or more of these logarithms will change by a multiple of $2 \pi i$.
For example, if $x = \pi - 2$, then $\ln( e^{-i (x + 2)} ) = -i (x + 2) + 2 \pi i$.
This ``extra'' $2 \pi i$ when $x \geq \pi - 2$ gives rise to a different polynomial than the $x/2$ in Equation \eqref{E:alternatingSum3x}.

If we substitute the logarithms in \eqref{E:AltSum3xSimpleLog1} through \eqref{E:AltSum3xSimpleLog3} into Equation \eqref{E:Li3Inversion}, we get
\begin{align} \label{E:AltSum3xFirstLog}
  \text{Li}_3( -e^{-i(x-2)} ) - \text{Li}_3( -e^{i (x-2)} )
  & = -\frac{\pi ^{2} }{6} (-i (x - 2)) - \frac{1}{6} (-i (x - 2))^3   \notag   \\
  & = -\frac{i x^3}{6}+i x^2+\frac{1}{6} i \pi ^2 x-2 i x-\frac{i \pi^2}{3}+\frac{4 i}{3} \, ,
\end{align}
\begin{align} \label{E:AltSum3xSecondLog}
  \text{Li}_3( -e^{-i x} ) - \text{Li}_3( -e^{i x} )  
  = -\frac{\pi ^{2} }{6} (-i x) - \frac{1}{6} (-i x))^3
  = \frac{1}{6} i \pi ^2 x-\frac{i x^3}{6} \, ,
\end{align}
\begin{align} \label{E:AltSum3xThirdLog}
  \text{Li}_3( -e^{-i(x+2)} ) - \text{Li}_3( -e^{i (x+2)} )
  & = -\frac{\pi ^{2} }{6} (-i (x + 2)) - \frac{1}{6} (-i (x + 2))^3   \notag   \\
  & = -\frac{i x^3}{6}-i x^2+\frac{1}{6} i \pi ^2 x-2 i x+\frac{i \pi^2}{3}-\frac{4 i}{3} \, .
\end{align}
We now substitute the expressions in Equations \eqref{E:AltSum3xFirstLog}, \eqref{E:AltSum3xSecondLog}, and \eqref{E:AltSum3xThirdLog}
into the right side of Equation \eqref{E:AltSum3x20}. We get
\begin{align} \label{E:AltSum3x30}
\frac{i}{8}
 \big(
  & -\frac{i x^3}{6}+i x^2+\frac{1}{6} i \pi ^2 x-2 i x-\frac{i \pi^2}{3}+\frac{4 i}{3}   \notag  \\
  & -2\left(  \frac{1}{6} i \pi ^2 x-\frac{i x^3}{6}  \right)   \notag  \\
  &  -\frac{i x^3}{6}-i x^2+\frac{1}{6} i \pi ^2 x-2 i x+\frac{i \pi^2}{3}-\frac{4 i}{3}
 \big)   \notag \\
   & = \frac{i}{8} (-4 i x)   \notag \\
   & = \frac{x}{2}
 \, .
\end{align}
Finally, from Equations \eqref{E:AltSum3x20} and \eqref{E:AltSum3x30}, we have our result:
\begin{equation} \label{E:AltSum3xResult1}
\sum_{n=1}^{\infty} (-1)^{n+1} \frac{\sin(nx)}{n} \cdot \left( \frac{\sin(n)}{n} \right) ^2 = \frac{x}{2} \, .
\end{equation}
As stated just after Equation \eqref{E:AltSum3xSimpleLog3}, this is subject to the condition that $-\pi + 2 \leq x < \pi - 2$,
in order that the logarithms will have the expressions given in Equations \eqref{E:AltSum3xSimpleLog1} through \eqref{E:AltSum3xSimpleLog3}.
We will prove below that this equation also holds for $x = \pi - 2$; see Equation \eqref{E:AltSum3xSecondExpr}.

A nice special case:
Substituting $x = 1$ into Equation \eqref{E:AltSum3xResult1}, and recalling Equation \eqref{E:alternatingSigns}, we get
\begin{equation} \label{E:ThreeAlternatingSincSums}
\sum_{n=1}^{\infty} (-1)^{n+1} \frac{\sin(n)}{n} = 
\sum_{n=1}^{\infty} (-1)^{n+1} \left( \frac{\sin(n)}{n} \right) ^2 =
\sum_{n=1}^{\infty} (-1)^{n+1} \left( \frac{\sin(n)}{n} \right) ^3 =
 \frac{1}{2} \, .
\end{equation}

Table \ref{Ta:SumSumSquaredAndCubed} shows the first 10 terms of each of these series, along with the sums of the first 10 terms.

\begin{table}[ht]
 \begin{center}
  \begin{tabular}{ r r r r r }
  \rule[-8pt]{0pt}{22pt}  
  $n$  &  $\sin(n)$   &   $(-1)^{n+1}\frac{\sin(n)}{n}$&   $(-1)^{n+1}\left(\frac{\sin(n)}{n}\right)^2$   &   $(-1)^{n+1}\left(\frac{\sin(n)}{n}\right)^3$  \\ \hline  
   1  & $  0.841471 $ &  $  0.841471 $ & $  0.7080734 $ & $  0.5958232 $  \\
   2  & $  0.909297 $ &  $ -0.454649 $ & $ -0.2067055 $ & $ -0.0939784 $  \\
   3  & $  0.141120 $ &  $  0.047040 $ & $  0.0022128 $ & $  0.0001041 $  \\
   4  & $ -0.756803 $ &  $  0.189201 $ & $ -0.0357969 $ & $  0.0067728 $  \\
   5  & $ -0.958924 $ &  $ -0.191785 $ & $  0.0367814 $ & $ -0.0070541 $  \\
   6  & $ -0.279416 $ &  $  0.046569 $ & $ -0.0021687 $ & $  0.0001010 $  \\
   7  & $  0.656997 $ &  $  0.093855 $ & $  0.0088088 $ & $  0.0008268 $  \\
   8  & $  0.989358 $ &  $ -0.123670 $ & $ -0.0152942 $ & $ -0.0018914 $  \\
   9  & $  0.412118 $ &  $  0.045791 $ & $  0.0020968 $ & $  0.0000960 $  \\
  10  & $ -0.544021 $ &  $  0.054402 $ & $ -0.0029596 $ & $  0.0001610 $  \\
  sum & $           $ &  $  0.548226 $ & $  0.4950484 $ & $  0.5009610 $  \\
  \end{tabular}
  \caption{The first 10 terms of the 3 series in \eqref{E:ThreeAlternatingSincSums}, and their sums}
  \label{Ta:SumSumSquaredAndCubed}
 \end{center}
\end{table}

\textbf{A series where $\sum a_n = \sum (a_n)^3$.}

Equation \eqref {E:sincn} and Theorem \ref{T:Thm1} provided an example of a series where $\sum a_n = \sum (a_n)^2$ :
\begin{equation*}
\sum_{n=1}^{\infty} \frac{\sin(n)}{n} = \sum_{n=1}^{\infty} \left(\frac{\sin(n)}{n} \right) ^{2} \, .
\end{equation*}

Here is a series where $\sum a_n = \sum (a_n)^3$ :
Regardless of whether $n$ is even or odd, $(-1)^{n+1} = ((-1)^{n+1})^3$.
Therefore, the first and third sums in \eqref{E:ThreeAlternatingSincSums} can be written
\begin{equation} \label{E:SumAndSumCubed}
\sum_{n=1}^{\infty} (-1)^{n+1} \frac{\sin(n)}{n} = 
\sum_{n=1}^{\infty} \left( (-1)^{n+1} \frac{\sin(n)}{n} \right ) ^3 =
 \frac{1}{2} \, .
\end{equation}
See Equation \eqref{E:Ballif3} for another series with this property.

\ifthenelse {\boolean{BKMRK}}
  { \subsection{Extending equation \ref{E:AltSum3xResult1} to \texorpdfstring{$x \geq \pi - 2$}{x > pi - 2}} \label{S:GeqPiMinus2}  }
  { \subsection{Extending equation \ref{E:AltSum3xResult1} to $x \geq \pi - 2$} \label{S:GeqPiMinus2}  }

\mbox{}  

Equation \eqref{E:AltSum3x20} holds for \emph{all} $x$.
However, for $x$ in different intervals, the resulting \emph{polynomial} expressions for the sum will be different.
This happens because for different $x$, the principal branch of the logarithm function may require different multiples of $2 \pi i$.
We will now show that, for $\pi - 2 \leq x < \pi$, the sum in \eqref{E:AltSum3xResult1} is not $x/2$, but is instead $x/2 - (\pi/8)(\pi - 2 - x)^2$.


Let's revisit the earlier expressions for the logarithms.
For somewhat larger values of $x$, the expressions for the logarithms are:
\begin{equation} \label{E:AltSum3xLogExpr1}
  \ln( e^{-i x} ) =
  \begin{cases}
    -i x  &\text{for $-\pi \leq x < \pi$,}\\
    -i x + 2 \pi i  &\text{for $\pi \leq x < 3\pi$.}
  \end{cases}
\end{equation}
\begin{equation} \label{E:AltSum3xLogExpr2}  
  \ln( e^{-i (x - 2)} ) = 
  \begin{cases}
    -i (x - 2)  &\text{for $-\pi + 2 \leq x < \pi + 2$,}\\
    -i (x - 2) + 2 \pi i  &\text{for $\pi + 2 \leq x < 3\pi + 2$.}
  \end{cases}
\end{equation}
\begin{equation} \label{E:AltSum3xLogExpr3}
  \ln( e^{-i (x + 2)} ) = 
  \begin{cases}
    -i (x + 2)  &\text{for $-\pi - 2 \leq x < \pi - 2$,}\\
    -i (x + 2) + 2 \pi i  &\text{for $\pi - 2 \leq x < 3\pi - 2$.}
  \end{cases}
\end{equation}
Let's derive an expression for the sum in Equation \eqref{E:AltSum3xResult1} for $\pi - 2 \leq x < \pi$.
In this interval, we must use the second expression for the logarithm given in Equation \eqref{E:AltSum3xLogExpr3}.
Note that when $x = \pi$,  we must also use the second expression for $\ln( e^{-i x} )$ in Equation \eqref{E:AltSum3xLogExpr1}.
But as long as $\pi - 2 \leq x < \pi$, Equations \eqref{E:AltSum3xFirstLog} and \eqref{E:AltSum3xSecondLog} remain unchanged from the previous section,
while Equation \eqref{E:AltSum3xThirdLog} becomes
\begin{align*}
  & \text{Li}_3( -e^{-i(x+2)} ) - \text{Li}_3( -e^{i (x+2)} ) \\
  & = -\frac{1}{6} \pi ^2 (-i (x+2) + 2 i \pi) - \frac{1}{6} (-i (x+2) + 2 i \pi)^3 \\
  & = -\frac{i x^3}{6}+i \pi  x^2-i x^2-\frac{11}{6} i \pi ^2 x+4 i \pi  x-2 i x+i \pi ^3-\frac{11 i \pi ^2}{3}+4 i \pi-\frac{4 i}{3}
\end{align*}

\begin{figure}[ht]
  \mbox{\includegraphics[width=\picSingleWidth]{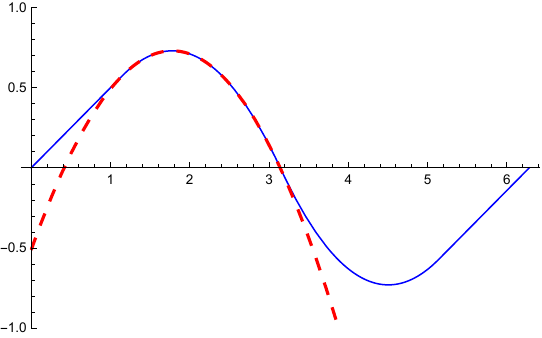}}
  \caption{$\sum_{n=1}^{100} (-1)^{n+1} \frac{\sin(nx)}{n} \left( \frac{\sin(n)}{n} \right)^2 $ (solid) and the quadratic in Equation \eqref{E:AltSum3xSecondExpr} (dashed)}
  \label{fig:altSum3Second}
\end{figure}

With this new expression, Equation \eqref{E:AltSum3x30} becomes
\begin{align*} 
\frac{i}{8}
 \big(
  & -\frac{i x^3}{6}+i x^2+\frac{1}{6} i \pi ^2 x-2 i x-\frac{i \pi^2}{3}+\frac{4 i}{3}   \notag  \\
  & -2\left(  \frac{1}{6} i \pi ^2 x-\frac{i x^3}{6}  \right)   \notag  \\
  & -\frac{i x^3}{6}+i \pi  x^2-i x^2-\frac{11}{6} i \pi ^2 x+4 i \pi x-2 i x+i \pi ^3-\frac{11 i \pi ^2}{3}+4 i \pi -\frac{4 i}{3}
 \big)   \notag \\
   & = \frac{i}{8} ( i \pi  x^2-2 i \pi ^2 x+4 i \pi  x-4 i x+i \pi ^3-4 i \pi ^2+4 i \pi )   \notag \\
   & = -\frac{\pi  x^2}{8}+\frac{\pi ^2 x}{4}-\frac{\pi x}{2}+\frac{x}{2}-\frac{\pi ^3}{8}+\frac{\pi ^2}{2}-\frac{\pi}{2}   \notag \\
   & = \frac{x}{2} - \frac{\pi}{8} (\pi - 2 - x)^2
 \, .
\end{align*}
So, for $\pi - 2 \leq x < \pi$,
\begin{equation} \label{E:AltSum3xSecondExpr}
  \sum_{n=1}^{\infty} (-1)^{n+1} \frac{\sin(nx)}{n} \cdot \left( \frac{\sin(n)}{n} \right) ^2
  = \frac{x}{2} - \frac{\pi}{8} (\pi - 2 - x)^2
  \, .
\end{equation}
At $x = \pi - 2$, the right hand side reduces to $x/2$, so  for this $x$, the sum in Equations \eqref{E:alternatingSum3x} and \eqref{E:AltSum3xSecondExpr} reduces to $x/2$.
At $x = \pi$, both sides of \eqref{E:AltSum3xSecondExpr} equal 0, so this equation is valid for $\pi - 2 \leq x \leq \pi$.
Figure \ref{fig:altSum3Second} shows the graphs of the left and right sides of Equation \eqref{E:AltSum3xSecondExpr}.

Putting Equations \eqref{E:AltSum3xResult1} and \eqref{E:AltSum3xSecondExpr} together, we have proved that
\[
\sum_{n=1}^{\infty} (-1)^{n+1} \frac{\sin(nx)}{n} \left( \frac{\sin(n)}{n} \right)^2 =
  \begin{dcases}  
    \frac{x}{2}  & \text{ for $-\pi + 2 \leq x < \pi - 2 \, ,$}\\
    \frac{x}{2} - \frac{\pi}{8} ( \pi - 2 - x )^2  &\text{ for $\pi - 2 \leq x \leq \pi \, .$}
  \end{dcases}
\]

\ifthenelse {\boolean{BKMRK}}
  { \subsection{A fourth sum that equals \texorpdfstring{$x/2$}{x/2}} }
  { \subsection{A fourth sum that equals $x/2$} }

\mbox{}  

We have already seen that, for $k = 0$, $k = 1$, and $k = 2$, the following holds for $0 \leq x < \pi - k$:
\begin{equation*}
  \sum_{n=1}^{\infty} (-1)^{n+1} \frac{\sin(nx)}{n} \cdot \left( \frac{\sin(n)}{n} \right) ^k = \frac{x}{2} \, .
\end{equation*}
(For $k = 1$ and $k = 2$ this actually holds for $0 \leq x \leq \pi - k$).
Can we extend this to higher values of $k$?
$\pi - 4 < 0$, but $\pi - 3 > 0$, so the case $k = 3$ is worth checking.
Numeric experiments suggest that, for $0 \leq x \leq \pi - 3$, this sum is \emph{also} $x/2$:
\begin{equation} \label{E:AltSum4x}
  \sum_{n=1}^{\infty} (-1)^{n+1} \frac{\sin(nx)}{n} \cdot \left( \frac{\sin(n)}{n} \right) ^3 = \frac{x}{2} \, .
\end{equation}
For various values of $x$ in this interval, we can calculate that the sum of the first million terms
is roughly $x/2$.
Also, for $x = 1/10$, $x = 1/8$, or $x = 2/15$, \textit{Mathematica} returns the exact value $x/2$.

We can prove Equation \eqref{E:AltSum4x} using the same techniques we used to prove Equation \eqref{E:AltSum3xResult1}.
We won't go through all the steps in as much detail, but here is an outline of the proof.

First, the sum in Equation \eqref{E:AltSum4x} can be expressed in terms of the polylogarithm $\text{Li}_4( \cdot )$.
For this sum, the equivalent of Equation \eqref{E:AltSum3x20} is
\begin{align} \label{E:AltSum4xPolylogs}
& \sum_{n=1}^{\infty} (-1)^{n+1} \frac{\sin(nx)}{n} \cdot \left( \frac{\sin(n)}{n} \right) ^3 = \notag  \\
& \frac{1}{16}
   \big \{ \;
      \text{Li}_4( -e^{-i(x-3)} ) + \text{Li}_4( -e^{i (x-3)} ) \notag  \\
    & - 3 \big(  \text{Li}_4( -e^{-i(x-1)} ) + \text{Li}_4( -e^{i (x-1)} )  \big) \notag  \\
    & + 3 \big(  \text{Li}_4( -e^{-i(x+1)} ) + \text{Li}_4( -e^{i (x+1)} )  \big) \notag  \\
    & -   \big( \text{Li}_4( -e^{-i (x+3)} ) + \text{Li}_4( -e^{i (x+3)} )  \big)
      \; \big \}
    \, .
\end{align}

The inversion identity for fourth order polylogarithms is (see Equations A.2.7(4) and A.2.7(6) on \cite[p.\ 298]{Lewin}):
\begin{equation} \label{E:Li4Inversion}
  \text{Li}_{4} (-z) + \text{Li}_{4} (-1/z) = -\frac{7 \pi^4}{360} - \frac{\pi^2 \ln^2(z)}{12} - \frac{\ln^{4}(z)}{24} \, .
\end{equation}
We will use this with four values of $z$: $z = e^{-i(x-3)}$, $z = e^{-i(x-1)}$, $z = e^{-i(x+1)}$, and $z = e^{-i(x+3)}$.

The four logarithms we need in Equation \eqref{E:Li4Inversion} are:
\begin{align}
\ln( e^{-i (x - 3)} ) & = -i (x - 3) \, , \label{E:AltSum4xSimpleLog1} \\
\ln( e^{-i (x - 1)} ) & = -i (x - 1) \, , \label{E:AltSum4xSimpleLog2} \\
\ln( e^{-i (x + 1)} ) & = -i (x + 1) \, , \label{E:AltSum4xSimpleLog3} \\
\ln( e^{-i (x + 3)} ) & = -i (x + 3) \, . \label{E:AltSum4xSimpleLog4}
\end{align}
Each $\ln( e^{-i (x - c)} ) = -i (x - c)$ holds (without an extra $2 \pi i$) for $-\pi + c \leq x < \pi + c$.
They hold for positive $x$ less than $\pi + 3$, $\pi + 1$, $\pi - 1$, and $\pi - 3$, respectively.
Therefore, \emph{all four} of these hold for $0 \leq x < \pi - 3$.

The first pair of $\text{Li}_{4}( \cdot )$ values on the right hand side of Equation \eqref{E:AltSum4xPolylogs} is
\begin{align}  \label{E:Sum4xExMinus3}
  & \text{Li}_4( -e^{-i(x-3)} ) + \text{Li}_4( -e^{i (x-3)} )   \notag \\
  & = -\frac{7 \pi^4}{360} - \frac{\pi^2 (-i (x - 3))^2}{12} - \frac{ (-i (x - 3))^4 }{24}   \notag \\
  & = -\frac{1}{24} (x-3)^4+\frac{1}{12} \pi ^2 (x-3)^2-\frac{7 \pi ^4}{360}
  \, .
\end{align}
The \emph{fourth} pair of $\text{Li}_{4}( \cdot )$ values in Equation \eqref{E:AltSum4xPolylogs} is
\begin{align} \label{E:Sum4xExPlus3}
  & \text{Li}_4( -e^{-i(x+3)} ) + \text{Li}_4( -e^{i (x+3)} )   \notag \\
  & = -\frac{7 \pi^4}{360} - \frac{\pi^2 (-i (x + 3))^2}{12} - \frac{ (-i (x + 3))^4 }{24}   \notag \\
  & = -\frac{1}{24} (x+3)^4+\frac{1}{12} \pi ^2 (x+3)^2-\frac{7 \pi ^4}{360}
  \, .
\end{align}
If we expand these expressions and subtract, there's a lot of cancellation and we get
\begin{align} \label{E:Sum4xOneMinusFour}
& \text{Li}_4( -e^{-i(x-3)} ) + \text{Li}_4( -e^{i (x-3)} )
 - \big( \text{Li}_4( -e^{-i (x+3)} ) + \text{Li}_4( -e^{i (x+3)} )  \big)   \notag \\
& = 9 x - \pi^2 x + x^3
\, .
\end{align}

The \emph{second} pair of $\text{Li}_{4}( \cdot )$ values in Equation \eqref{E:AltSum4xPolylogs} is
\begin{align*}
  & \text{Li}_4( -e^{-i(x-1)} ) + \text{Li}_4( -e^{i (x-1)} )  \\
  & = -\frac{7 \pi^4}{360} - \frac{\pi^2 (-i (x - 1))^2}{12} - \frac{ (-i (x - 1))^4 }{24}  \\
  & = -\frac{1}{24} (x-1)^4+\frac{1}{12} \pi ^2 (x-1)^2-\frac{7 \pi ^4}{360}
  \, .
\end{align*}
The \emph{third} pair of $\text{Li}_{4}( \cdot )$ values in Equation \eqref{E:AltSum4xPolylogs} is
\begin{align*}
  & \text{Li}_4( -e^{-i(x+1)} ) + \text{Li}_4( -e^{i (x+1)} )  \\
  & = -\frac{7 \pi^4}{360} - \frac{\pi^2 (-i (x + 1))^2}{12} - \frac{ (-i (x + 1))^4 }{24}  \\
  & = -\frac{1}{24} (x+1)^4+\frac{1}{12} \pi ^2 (x+1)^2-\frac{7 \pi ^4}{360}
  \, .
\end{align*}
Note that, in Equation \eqref{E:AltSum4xPolylogs}, we subtract the third sum minus the second, then multiply the result by 3.
So, combining the last two results, we get
\begin{align} \label{E:Sum4xThreeMinusTwo}
  & \big(  \text{Li}_4( -e^{-i(x+1)} ) + \text{Li}_4( -e^{i (x+1)} )  \big)
  - \big(  \text{Li}_4( -e^{-i(x-1)} ) + \text{Li}_4( -e^{i (x-1)} )  \big)   \notag \\
  & = \frac{1}{24} (x-1)^4-\frac{1}{12} \pi ^2 (x-1)^2-\frac{1}{24} (x+1)^4+\frac{1}{12} \pi ^2 (x+1)^2   \notag \\
  & = -\frac{x^3}{3}+\frac{\pi ^2 x}{3}-\frac{x}{3}
   \, .
\end{align}
We now combine the results of Equations \eqref{E:Sum4xOneMinusFour} and \eqref{E:Sum4xThreeMinusTwo}:
\begin{align} \label{E:AltSum4xAllPolylogSums}
& \text{Li}_4( -e^{-i(x-3)} ) + \text{Li}_4( -e^{i (x-3)} )
 - \big( \text{Li}_4( -e^{-i (x+3)} ) + \text{Li}_4( -e^{i (x+3)} )  \big)   \notag \\
& + 3 \big\{ 
         \big(  \text{Li}_4( -e^{-i(x+1)} ) + \text{Li}_4( -e^{i (x+1)} )  \big)
       - \big(  \text{Li}_4( -e^{-i(x-1)} ) + \text{Li}_4( -e^{i (x-1)} )  \big)  \big\}   \notag \\
& = 9 x - \pi^2 x + x^3 + 3 \left( -\frac{x^3}{3}+\frac{\pi ^2 x}{3}-\frac{x}{3} \right)   \notag \\
& = 8x \, .
\end{align}
Accounting for the factor of $1/16$ on the right side of Equation \eqref{E:AltSum4xPolylogs}, we have, for $0 \leq x < \pi - 3$,
\begin{equation} \label{E:AltSum4xResult1}
  \sum_{n=1}^{\infty} (-1)^{n+1} \frac{\sin(nx)}{n} \cdot \left( \frac{\sin(n)}{n} \right) ^3 = \frac{x}{2} \, .
\end{equation}
We will prove below that this equation also holds for $x = \pi - 3$; see Equation \eqref{E:AltSum4xSecondExpr}.

\ifthenelse {\boolean{BKMRK}}
  { \subsection{Extending equation \ref{E:AltSum4xResult1} to \texorpdfstring{$x \geq \pi - 3$}{x > pi - 3}} \label{S:AltSum4xResult2} }
  { \subsection{Extending equation \ref{E:AltSum4xResult1} to $x \geq \pi - 3$}  \label{S:AltSum4xResult2} }

\mbox{}  

What happens to the sum if $x \geq \pi - 3$?
For somewhat larger $x$, we will need to use this expression instead of the one in Equation \eqref{E:AltSum4xSimpleLog4}:
\begin{equation} \label{E:AltSum4xLogExpr4}
  \ln( e^{-i (x + 3)} ) = -i (x + 3) + 2 \pi i \, .
\end{equation}
This holds for $\pi - 3 \leq x < \pi - 3 + 2 \pi = 3 \pi - 3 \approx 6.42$.
As long as $x < \pi - 1$, the other logarithms in \eqref{E:AltSum4xSimpleLog1} through \eqref{E:AltSum4xSimpleLog3} do not need a $2 \pi i$ added to them.

We follow the same procedure as before.
The only difference is that Equation \eqref{E:Sum4xExPlus3} becomes
\begin{align} \label{E:Sum4xExPlus3B}
  & \text{Li}_4( -e^{-i(x+3)} ) + \text{Li}_4( -e^{i (x+3)} )   \notag \\
  & = -\frac{7 \pi^4}{360} - \frac{\pi^2 (-i (x + 3) + 2 \pi i)^2}{12} - \frac{ (-i (x + 3) + 2 \pi i)^4 }{24} \, .
\end{align}
Now, when we do this subtraction, there is less cancellation than we saw in Equation \eqref{E:Sum4xOneMinusFour}:
\begin{align} \label{E:Sum4xOneMinusFourB}
& \text{Li}_4( -e^{-i(x-3)} ) + \text{Li}_4( -e^{i (x-3)} )
 - \big( \text{Li}_4( -e^{-i (x+3)} ) + \text{Li}_4( -e^{i (x+3)} )  \big)   \notag \\
& = -\frac{\pi  x^3}{3}+x^3+\pi ^2 x^2-3 \pi  x^2-\pi ^3 x+5 \pi ^2 x - 9 \pi x + 9x+\frac{\pi ^4}{3}-3 \pi ^3+9 \pi ^2-9 \pi
\, .
\end{align}

We now combine the results of Equations \eqref{E:Sum4xOneMinusFourB} and \eqref{E:Sum4xThreeMinusTwo}:
\begin{align*}
& \text{Li}_4( -e^{-i(x-3)} ) + \text{Li}_4( -e^{i (x-3)} )
 - \big( \text{Li}_4( -e^{-i (x+3)} ) + \text{Li}_4( -e^{i (x+3)} )  \big)   \notag \\
& + 3 \big( 
         \big(  \text{Li}_4( -e^{-i(x+1)} ) + \text{Li}_4( -e^{i (x+1)} )  \big)
       - \big(  \text{Li}_4( -e^{-i(x-1)} ) + \text{Li}_4( -e^{i (x-1)} )  \big)  \big)   \notag \\
& = -\frac{\pi  x^3}{3}+x^3+\pi ^2 x^2-3 \pi  x^2-\pi ^3 x+5 \pi ^2 x - 9 \pi x + 9x+\frac{\pi ^4}{3}-3 \pi ^3+9 \pi ^2-9 \pi   \notag \\
& + 3 \left( -\frac{x^3}{3}+\frac{\pi ^2 x}{3}-\frac{x}{3} \right) \\
& = -\frac{\pi  x^3}{3}+\pi ^2 x^2-3 \pi  x^2-\pi ^3 x+6 \pi ^2 x-9 \pi  x+8x+\frac{\pi ^4}{3}-3 \pi ^3+9 \pi ^2-9 \pi   \notag \\
& = 8x + \frac{\pi}{3} (\pi - 3 - x)^3
\, .
\end{align*}

Accounting for the factor of $1/16$ on the right side of Equation \eqref{E:AltSum4xPolylogs}, we have, for $\pi - 3 \leq x < \pi - 1$,
\begin{equation} \label{E:AltSum4xSecondExpr}
  \sum_{n=1}^{\infty} (-1)^{n+1} \frac{\sin(nx)}{n} \cdot \left( \frac{\sin(n)}{n} \right) ^3
   = \frac{x}{2} + \frac{\pi}{48} (\pi - 3 - x)^3
   \, .
\end{equation}
Testing this for $x = \pi - 3$, the cubic term on the right is 0, so the right side reduces to $x/2$.
Therefore, Equation \eqref{E:AltSum4xResult1} actually holds for $0 \leq x \leq \pi - 3$.

But note that, for $x = 1$, the right-hand side of \eqref{E:AltSum4xSecondExpr} is no longer $1/2$,
unlike Equation \eqref{E:ThreeAlternatingSincSums}, where the sums with $\sin(n)/n$ to the first, second, and third powers, are exactly $1/2$.

For $x = 1/7$, which is slightly more than $\pi - 3$, the sum of the series in Equation \eqref{E:AltSum4xSecondExpr} is
\[
\frac{x}{2} + \frac{\pi}{48} (\pi - (3 + x))^3 = 
  \frac{1}{14} + \frac{\pi}{48} \left(\pi -\frac{22}{7} \right)^3
  \approx \frac{1}{14} - 1.32 \cdot 10^{-10}
   \, .
\]
Note the occurrence of $22/7$, a popular approximation to $\pi$.

Using $355/113$, another good approximation to $\pi$, we can take $x = 16/113$, which is slightly more than $\pi - 3$.
The right-hand side of \eqref{E:AltSum4xSecondExpr} is
\[
  \frac{8}{113} + \frac{\pi}{48} \left(\pi -\frac{355}{113} \right)^3
  = \frac{8}{113} + \frac{\pi}{48} \cdot (-1.898 \cdot 10^{-20}) \, .
\]


Let's collect into a Theorem the related results that we've proven.
\begin{theorem} \label{T:AltSumkxTheorem}
\begin{align*}
\sum_{n=1}^{\infty} (-1)^{n+1} \frac{\sin(nx)}{n}                                          & = \frac{x}{2}  \quad \text{ for $0 \leq x < \pi$} \, .        & \eqref{E:xOver2ab} \\
\sum_{n=1}^{\infty} (-1)^{n+1} \frac{\sin(nx)}{n} \cdot \frac{\sin(n)}{n}                  & = \frac{x}{2}  \quad \text{ for $0 \leq x \leq \pi - 1$} \, . & \eqref{E:xOver2ab} \\
\sum_{n=1}^{\infty} (-1)^{n+1} \frac{\sin(nx)}{n} \cdot \left( \frac{\sin(n)}{n} \right)^2 & = \frac{x}{2}  \quad \text{ for $0 \leq x \leq \pi - 2$} \, . & \eqref{E:AltSum3xResult1} \\
\sum_{n=1}^{\infty} (-1)^{n+1} \frac{\sin(nx)}{n} \cdot \left( \frac{\sin(n)}{n} \right)^3 & = \frac{x}{2}  \quad \text{ for $0 \leq x \leq \pi - 3$} \, . & \eqref{E:AltSum4xResult1}
\end{align*}
\end{theorem}

Note that, as a special case ($x = 1$), we have
\begin{equation} \label{E:3sumsOneHalf}
\sum_{n = 1}^{\infty} (-1)^{n+1} \frac{\sin(n)}{n} =
\sum_{n = 1}^{\infty} (-1)^{n+1} \left( \frac{\sin(n)}{n} \right)^2 =
\sum_{n = 1}^{\infty} (-1)^{n+1} \left( \frac{\sin(n)}{n} \right)^3 = \frac{1}{2} \, .
\end{equation}

For somewhat larger values of $x$, the sums are:
\begin{align*}
\sum_{n=1}^{\infty} (-1)^{n+1} \frac{\sin(nx)}{n} \cdot \frac{\sin(n)}{n} & =
    \frac{x}{2} + \frac{\pi}{2}(\pi - 1 - x)  \quad \text{ for $\pi - 1 \leq x \leq \pi + 1$} \, . & \eqref{E:altSignG2} \\
\sum_{n=1}^{\infty} (-1)^{n+1} \frac{\sin(nx)}{n} \cdot \left( \frac{\sin(n)}{n} \right)^2 & =
    \frac{x}{2} - \frac{\pi}{8} (\pi - 2 - x)^2  \quad \text{ for $\pi - 2 \leq x \leq \pi$} \, . & \eqref{E:AltSum3xSecondExpr} \\
\sum_{n=1}^{\infty} (-1)^{n+1} \frac{\sin(nx)}{n} \cdot \left( \frac{\sin(n)}{n} \right)^3 & =
    \frac{x}{2} + \frac{\pi}{48} (\pi - 3 - x)^3  \quad \text{ for $\pi - 3 \leq x < \pi - 1$} \, . & \eqref{E:AltSum4xSecondExpr}
\end{align*}


\subsection{\textit{Mathematica} code to evaluate polylogarithms} 

\mbox{}  

\textit{Mathematica} implements the $n$\textsuperscript{th} order polylogarithm function $\text{Li}_{n} (z)$ as \verb+PolyLog[n, z]+.
However, to express a sum of polylogs as a polynomial, we must do some of the work ourselves.

The $n$\textsuperscript{th} order polylogarithm has the following inversion equation
\begin{equation} \label{E:LinInversion}
  \text{Li}_{n} (-z) + (-1)^n \text{Li}_{n} (-1/z)
  = -\frac{\ln^n(z)}{n!} + 2 \sum_{k=1}^{\left \lfloor n/2 \right \rfloor} \frac{\ln^{n-2k}(z)}{(n-2k)!} \text{Li}_{2k}(-1) \, .
\end{equation}
This is Equation A.2.7(6) from Lewin's book \cite[p.\ 299]{Lewin}.
This is also Equation 89 in \cite[p.\ 116]{Srivastava}.
For $n = 2, 3,$ and 4, this gives Equations \eqref{E:Li2Inversion}, \eqref{E:Li3Inversion}, and \eqref{E:Li4Inversion}.
Note that the right side is a polynomial in $\ln(z)$.

Equation \eqref{E:LinInversion} uses the values of $\text{Li}_{2}(-1)$, $\text{Li}_{4}(-1)$, etc.
For even $n > 0$, these can be computed with this \textit{Mathematica} code \cite[Eq. 94, p.\ 117]{Srivastava}:
\begin{verbatim}
  liEvenNegOne[n_?(# > 0 && EvenQ[#] &)] :=
    (-1)^(n/2) * (2^(n - 1) - 1) * BernoulliB[n] * Pi^n/n! ;
\end{verbatim}

Caution: Lewin \cite[Eq. A.2.7(4), p.\ 299]{Lewin} indexes the Bernoulli numbers differently than \textit{Mathematica}; he omits the odd-numbered Bernoulli numbers (these are all zero), and considers all the even-numbered Bernoulli numbers to be positive.

The polynomial on the right of Equation \eqref{E:LinInversion} can be computed with this \textit{Mathematica} code:
\begin{verbatim}
  inversionOrderN[n_?(# > 1 && IntegerQ[#] &), t_] := -(t^n/n!) +
    2 * Sum[ t^(n - 2 k)/((n - 2 k)!) * liEvenNegOne[2 k], {k, 1, Floor[n/2]}] ;
\end{verbatim}

With this code, \verb+inversionOrderN[2, t]+, \verb+inversionOrderN[3, t]+ and \verb+inversionOrderN[4, t]+
give the inversion equations that we used in Equations \eqref{E:Li2Inversion}, \eqref{E:Li3Inversion} and \eqref{E:Li4Inversion}.

For example, this converts the right-hand side of \eqref{E:trilogExprSum} into the right-hand side of \eqref{E:trilogExpr}:
\begin{verbatim}
  Expand[ (I/2) * inversionOrderN[3, Log[-E^(-I)]] ]
\end{verbatim}

Using the logarithms given in Equations \eqref{E:AltSum4xSimpleLog1} through \eqref{E:AltSum4xSimpleLog4}, we show how the above \textit{Mathematica} code can produce some of the results stated above.

For example, this sum in Equation \eqref{E:Sum4xExMinus3}
\[
  \text{Li}_4( -e^{-i(x-3)} ) + \text{Li}_4( -e^{i (x-3)} )
\]
can be computed with \verb+inversionOrderN[4, -I(x-3)]+ .

It will be useful here to define a shorthand version of \verb+inversionOrderN[4, t]+:
\begin{verbatim}
  inv4[t_] := inversionOrderN[4, t] ;
\end{verbatim}

\textbf{Checking Equations \eqref{E:AltSum4xAllPolylogSums} and \eqref{E:AltSum4xSecondExpr}}.

The combination of the four polylog sums in Equation \eqref{E:AltSum4xAllPolylogSums} can be computed with
\begin{verbatim}
  poly = inv4[-I (x - 3)] - inv4[-I (x + 3)] +
     3*( inv4[-I (x + 1)] - inv4[-I (x - 1)] )
  Expand[poly]
\end{verbatim}
\verb+Expand+ expands the products, then cancels as many terms as possible.
The result is $8x$.

To get the result in Equation \eqref{E:AltSum4xSecondExpr}, we add $2 \pi i$ to the logarithm in the seond term:
\begin{verbatim}
  poly2 = inv4[-I (x - 3)] - inv4[-I (x + 3) + 2 Pi I] +
      3*( inv4[-I (x + 1)] - inv4[-I (x - 1)] )
  Expand[poly2]
\end{verbatim}

Equation \eqref{E:AltSum4xSecondExpr} is \verb+poly2+ divided by 16.
At this point, \verb+FullSimplify[poly2/16]+ will give an expression equivalent to the right side of Equation \eqref{E:AltSum4xSecondExpr}.

\textbf{Extending Equation \eqref{E:AltSum4xSecondExpr} to larger $x$ and verifying our work}.

Just for fun, let's get a polynomial representation for the sum in Equation \eqref{E:AltSum4xSecondExpr}, this time, for $\pi - 1 < x < \pi + 1$.
For this interval, \emph{two} of the logarithms require an extra $2 \pi i$.
We can extend the above calculation to this interval as follows:
\begin{verbatim}
  poly3 = inv4[-I (x - 3)] - inv4[-I (x + 3) + 2 Pi I] +
      3*( inv4[-I (x + 1) + 2 Pi I] - inv4[-I (x - 1)] )
  FullSimplify[poly3/16]
\end{verbatim}
This gives
\[
-\frac{1}{24} \left(12 + \pi(-9 + (\pi - x)^2) \right) (\pi-x) =
\frac{x}{2} + \frac{\pi}{24} (x - \pi) \left( (x - \pi)^2 - 9 \right) - \frac{\pi}{2}
\]

So, for $\pi - 1 < x < \pi + 1$, we have
\[
\sum_{n=1}^{\infty} (-1)^{n+1} \frac{\sin(nx)}{n} \cdot \left( \frac{\sin(n)}{n} \right)^3 =
\frac{x}{2} + \frac{\pi}{24} (x - \pi) \left( (x - \pi)^2 - 9 \right) - \frac{\pi}{2} \, .
\]

Now let's define $f(x)$ over $0 \leq x \leq \pi$ to consist of the following polynomials:
\begin{equation*}
  f(x) =
  \begin{dcases}  
   \frac{x}{2}  &\text{if $0 \leq x \leq \pi - 3$,} \\
   \frac{x}{2} + \frac{\pi}{48} (\pi - 3 - x)^3  &\text{if $\pi - 3 < x  \leq \pi - 1$,} \\
   \frac{x}{2} + \frac{\pi}{24} (x - \pi) \left( (x - \pi)^2 - 9 \right) - \frac{\pi}{2}   &\text{if $\pi - 1 < x \leq \pi$.}
  \end{dcases}
\end{equation*}
(The third expression for $f(x)$ also works for $\pi < x < \pi+1$, but here, we don't care about $x > \pi$.)
Let's also define $f(x)$ for $-\pi \leq x < 0$ by taking the odd periodic extension of $f(x)$.

Now, let's compute the Fourier sine coefficients of $f(x)$.
We claim that $f(x)$ is represented by the Fourier sine series
\[
  \sum_{n=1}^{\infty} (-1)^{n+1} \frac{\sin(nx)}{n} \cdot \left( \frac{\sin(n)}{n} \right) ^3
= \sum_{n=1}^{\infty}  (-1)^{n+1} \frac{\sin^3(n)}{n^4} \cdot \sin(nx)
\]
over (at least) the interval $0 < x < \pi$, so the $n$\textsuperscript{th} coefficient should be
\[
(-1)^{n+1} \frac{\sin^3(n)}{n^4} \, .
\]

This \textit{Mathematica} code uses Equation \eqref{E:bnOdd} to compute the $n$\textsuperscript{th} sine coefficient of $f(x)$:
\begin{verbatim}
  i1 = Integrate[Sin[n x] * (x/2), {x, 0, Pi - 3}]
  i2 = Integrate[Sin[n x] * (x/2 + (Pi/48)*(Pi - 3 - x)^3), {x, Pi - 3, Pi - 1}]
  i3 = Integrate[Sin[n x] * (x/2 + (Pi/24)(x - Pi)((x - Pi)^2 - 9) - Pi/2),
          {x, Pi - 1, Pi}]
  FullSimplify[(2/Pi) * (i1 + i2 + i3), Assumptions -> Element[n, Integers]]
\end{verbatim}
The result is $(-1)^{n+1} \sin^3(n)/n^4$.
It is gratifying to see that all of these polynomials over their respective intervals do, indeed, have as their Fourier series the sum on the left side of Equation \eqref{E:AltSum4xSecondExpr}.
This provides additional confirmation that our calculations are correct.

\section{Opportunities for further exploration}\label{S:FurtherExploration}

We certainly have not exhausted the possibilities for exploration.
The reader might enjoy graphing
\[
\sum_{n=1}^{\infty} \frac{\sin ^{j} (nx)}{n ^{k}}
\]
for small $j$ and $k$.
When $j$ is odd and $k$ is 1, we get graphs like the solid graphs in Figures \ref{fig:hpPlot}, \ref{fig:hpPlot56}, and \ref{fig:plot78} (pages \pageref{fig:hpPlot} and \pageref{fig:plot78}).
For a given $j$, what are the discrete $y$ values on the graphs of these sums?
How many are there?  What are the $x$ values at the endpoints of the ``pieces''?
The graphs also illustrate why Equations \eqref{E:Prob6241a}, \eqref{E:SumSinCubedOverN}, and \eqref{E:SumSinFifthOverN} hold, and why higher-power versions do not.

Besides substituting $\pi - kx$ for $x$ as we did above, there are other ways to manipulate Fourier series to get interesting results.
For example, one may often differentiate or integrate a Fourier series term by term to produce a new Fourier series.
See \cite[pp.\ 125-129]{Tolstov}.

One can also consider the "integral analogues" of the sums here.
Notice that the values of the sums in \eqref{E:Prob6241a} are exactly 1/2 less than the values of these ``corresponding'' integrals:
\[
\int_{0}^{\infty } \frac{\sin(x)}{x} dx = \int_{0}^{\infty } \left(\frac{\sin(x)}{x} \right) ^{2} dx = \frac{\pi }{2} \, .
\]

After dividing through by 3, we can do the same for \eqref{E:FourSumsWithN}.
That is, for $k = 0$, 1, 2, and 3,
\[
\sum_{n=1}^{\infty} \left(\frac{\sin(n)}{n} \right) ^{k} \frac{\sin(3n)}{3n} =
 -\frac{1}{2} + \int_{0}^{\infty } \left(\frac{\sin(x)}{x} \right) ^{k} \frac{\sin(3x)}{3x} dx = \frac{\pi - 3}{6} \, .
\]

It is known \cite{MathWorldSincFunction} that, if $m$ is a positive integer, then
\[
\int_{0}^{\infty } \left(\frac{\sin(x)}{x} \right) ^{m} dx
\]
is a rational multiple of $\pi$.
Experiments with both numeric and symbolic calculations suggested that, for $m = 1$, 2, 3, 4, 5, and 6, we have
\[
\sum_{n=1}^{\infty} \left(\frac{\sin(n)}{n} \right) ^{m} = -\frac{1}{2} + \int_{0}^{\infty } \left(\frac{\sin(x)}{x} \right) ^{m} dx \, .
\]

However, for $m = 7$, the integral is
\[
\int_{0}^{\infty } \left(\frac{\sin(x)}{x} \right) ^{7} dx = \frac{5887\pi }{23040} ,
\]
but the sum now has a completely different form:
\[\sum_{n=1}^{\infty} \left(\frac{\sin(n)}{n} \right) ^{7} =\]
\[-\frac{1}{2} + \frac{129423 \pi - 201684 \pi^{2} + 144060 \pi^{3} - 54880 \pi^{4} + 11760 \pi^{5} - 1344 \pi^{6} + 64 \pi^{7} }{46080} \, .\]

The surprising appearance of this 7\textsuperscript{th}-degree polynomial and the relationships between these sums and integrals are explored in \cite{BBB}. 

We experimented mostly with the sine function here.
What results could we get by experimenting with the cosine function?


We have found numerous examples of series in which we can multiply the $n$\textsuperscript{th} term by $\text{sinc}(n)$, \emph{without changing the sum}.
Why does the $\text{sinc}$ function have this property?
Do other functions have this property?

\medskip

Here are some Fourier series with a slightly different form.
Most of the Fourier coefficients we've encountered so far have had some power of $n$ in the denominator.
By contrast, in these Fourier series, the denominators are $2^n$:
\begin{align}
\sum_{n = 1}^{\infty} \frac{\sin(nx)}{2^n} & = \frac{2 \sin(x)}{5 - 4 \cos(x)} \, , \label{E:Sin2ToN} \\
\sum_{n = 1}^{\infty} \frac{\cos(nx)}{2^n} & = \frac{-1 + 2 \cos(x)}{5 - 4 \cos(x)} \, . \label{E:Cos2ToN}
\end{align}

\begin{figure}[ht] 
\centering
\begin{minipage}[t]{\figureMinipageWidth}
\centering
\includegraphics[width=\picDblWidth]{{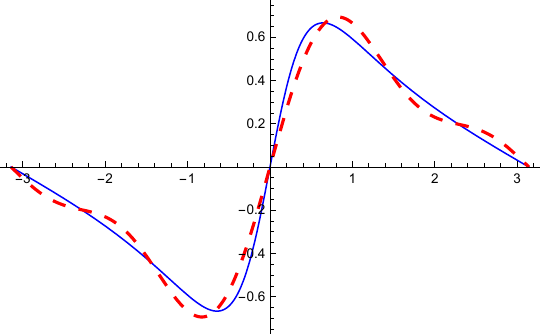}}
\caption{$\frac{2\sin(x)}{5 - 4 \cos(x)}$ (solid); sum of 3 terms of its Fourier series (dashed)}
\label{fig:sin2ToNPlot}
\end{minipage}\hfill
\begin{minipage}[t]{\figureMinipageWidth}
\centering
\includegraphics[width=\picDblWidth]{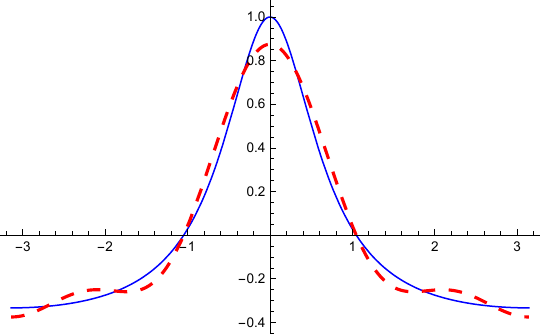}
\caption{$\frac{-1 + 2\cos(x)}{5 - 4 \cos(x)}$ (solid); sum of 3 terms of its Fourier series (dashed)} 
\label{fig:cos2ToNPlot}
\end{minipage}
\end{figure}

Starting with the right-hand sides of Equations \eqref{E:Sin2ToN} and \eqref{E:Cos2ToN}, \textit{Mathematica} seems unable to calculate the Fourier coefficients $1/2^n$ that appear on the left.
However, if you start with the left-hand sides and simplify,
\begin{verbatim}
  FullSimplify[ Sum[Sin[n x]/2^n, {n, 1, Infinity}] ]
  Sum[Cos[n x]/2^n, {n, 1, Infinity}]
\end{verbatim}
\textit{Mathematica} returns the expressions on the right sides of \eqref{E:Sin2ToN} and \eqref{E:Cos2ToN}.

What happens if we replace $2^n$ with $3^n$ in the above series?
\textit{Mathematica} claims that
\begin{align}
\sum_{n = 1}^{\infty} \frac{\sin(nx)}{3^n} & = \frac{3 \sin(x)}{10 - 6 \cos(x)} \, , \label{E:Sin3ToN} \\
\sum_{n = 1}^{\infty} \frac{\cos(nx)}{3^n} & = \frac{-1 + 3 \cos(x)}{10 - 6 \cos(x)} \, . \label{E:Cos3ToN}
\end{align}

A little more experimentation shows that \textit{Mathematica} can often evaluate similar sums where $2^n$ is replaced with $k^n$ where $k > 1$.
Here we will prove that, for any real number $k > 1$,
\begin{align}
\sum_{n = 1}^{\infty} \frac{\sin(nx)}{k^n} & = \frac{k \sin(x)}{k^2 + 1 - 2k \cos(x)} \, , \label{E:SinKToN} \\
\sum_{n = 1}^{\infty} \frac{\cos(nx)}{k^n} & = \frac{-1 + k \cos(x)}{k^2 + 1 - 2k \cos(x)} \, . \label{E:CosKToN}
\end{align}
In fact, we will do something even better: starting with the sums, we will \emph{derive} the expressions on the right sides of Equations \eqref{E:SinKToN} and \eqref{E:CosKToN} and show where these expressions come from.

The trick is to combine the sums on the left into a single sum and use the fact that, for all $n$ and all $x$, $\cos(nx) + i \sin(nx) = e^{inx}$:
\[
  \sum_{n = 1}^{\infty} \frac{\cos(nx) + i \sin(nx)}{k^n}
= \sum_{n = 1}^{\infty} \frac{e^{i n x}}{k^n}
= \sum_{n = 1}^{\infty} \left( \frac{e^{i x}}{k} \right)^n \, .
\]
The sum on the right side of this expression is a geometric series with first term $a = e^{ix}/k$ and common ratio $r = e^{ix}/k$.
If $k > 1$, then $\lvert r \rvert = 1/k < 1$ for all $x$, so the geometric series converges for all $x$.
The sum equals $a/(1 - r)$, which is
\[
\frac{ \frac{e^{i x}}{k} }{1 - \frac{e^{i x}}{k} }
 = \frac{ \frac{e^{i x}}{k} }{1 - \frac{e^{i x}}{k} } \cdot \frac{k}{k}
 = \frac{ e^{i x} }{ k - e^{i x} }
 = \frac{ \cos(x) + i \sin(x) }{ k - \cos(x) - i \sin(x) } \, .
\]
We can get the $i$ out of the denominator by multiplying both numerator and denominator by the conjugate of the denominator:
\[
\frac{ \cos(x) + i \sin(x) }{ k - \cos(x) - i \sin(x) }
 = \frac{ \cos(x) + i \sin(x) }{ k - \cos(x) - i \sin(x) } \cdot \frac{ k - \cos(x) + i \sin(x) }{ k - \cos(x) + i \sin(x) } \, .
\]
The numerator is
\begin{align*}
   & \left ( \cos(x) + i \sin(x) \right ) \cdot \Bigl( (k - \cos(x)) + i \sin(x) \Bigr) \\
 = & k \cos(x) - \cos^2(x) + i^2 \sin^2(x) + i \sin(x) \Bigl((k - \cos(x)) + \cos(x) \Bigr) \\
 = & (k \cos(x) - 1) + i \cdot k \sin(x) \, .
\end{align*}
The denominator is
\begin{align*}
   & \left (  k - \cos(x) - i \sin(x) \right ) \cdot \left ( k - \cos(x) + i \sin(x) \right ) \\
 = & (k -\cos(x))^2 - i^2 \sin^2(x) = k^2 - 2 k \cos(x) + \cos^2(x) + \sin^2(x) \\
 = & k^2 + 1 - 2 k \cos(x) \, .
\end{align*}
Finally, we have
\[
  \sum_{n = 1}^{\infty} \frac{\cos(nx)}{k^n} + i \cdot \sum_{n = 1}^{\infty} \frac{\sin(nx)}{k^n}
  = \frac{ (k \cos(x) - 1) + i \cdot k \sin(x) }{ k^2 + 1 - 2 k \cos(x) } \, ,
\]
and, separating the real and imaginary parts, Equations \eqref{E:SinKToN} and \eqref{E:CosKToN} follow.

The reader may enjoy experimenting with term-by-term integration and differentiation of the above series, and of the corresponding right-hand sides of Equations \eqref{E:SinKToN} and \eqref{E:CosKToN}.


\section{A Fourier series package in Mathematica} \label{S:Package}

Here we discuss \textit{Mathematica} package that computes Fourier series coefficients of a wide variety of functions.
This package makes it easy (and fun!) to experiment with different functions and see what their Fourier series look like.
We show how to use this package to compute the Fourier series of functions where the only information you need to enter are the pairs of $x$ and $y$ values at endpoints, corners, or jump discontinuities.

The package is in (text) file \verb+FS.m+, which has been uploaded to math arXiv as an ``ancillary file''.
The link to that file is given in Section \ref{S:WhatsNew}.


If you download the file to, say, \verb+C:/MATH/FS.m+, then you can read that file into \textit{Mathematica}
with (something like) the following \textit{Mathematica} command:
\begin{verbatim}
  << C:/MATH/FS.m
\end{verbatim}

The main function in this package is \verb+fSeries[ ]+.
An easy way to use this function is to specify these 6 parameters:

\begin{itemize}
 \item Parameter 1 is  is the name of the variable used in the function, typically, \verb+x+;
 \item parameter 2 is the function of \verb+x+;
 \item parameter 3 is $a$, the lefthand endpoint of the interval;
 \item parameter 4 is $b$, the righthand endpoint (note: the period is $b - a$);
 \item parameter 5 is the variable to be used in the coefficients, typically, \verb+n+;
 \item parameter 6 is the number of terms in the partial sum to include in the graph.
\end{itemize}

\begin{figure}[ht]
  \mbox{\includegraphics[width=\picDblWidth]{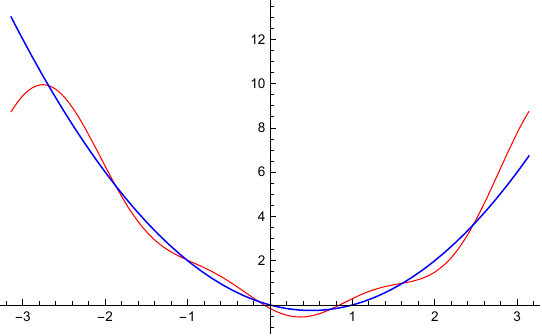}}
  \caption{$x^2 - x$ (blue) and the sum of 3 terms of its Fourier series (red)}
  \label{fig:quadratic3}
\end{figure}

Example:
\begin{verbatim}
  fSeries[x, x^2 - x, -Pi, Pi, n, 3]
\end{verbatim}
computes the coefficients, as functions of \verb+n+, of the Fourier series for $f(x) = x^2 - x$.
The interval is $(-\pi, \pi)$, so the period is $\pi - (-\pi) = 2 \pi$.
It also plots $f(x)$ along with the sum of 3 terms of its Fourier series; see Figure \ref{fig:quadratic3}.

In addition to drawing the graphs, \verb+fSeries+ returns a list with three expressions:
\begin{equation} \label{E:QuadraticList}
\left \{ \frac{\pi ^2}{3}, \quad \frac{4 (-1)^n}{n^2}, \quad \frac{2 (-1)^n}{n} \right \} \, .
\end{equation}
These three expressions are: the constant term $a_0/2$; and, for $n \geq 1$, the coefficients of $\cos(nx)$ and $\sin(nx)$.
So, the Fourier series for $x^2 - x$ over $(-\pi, \pi)$ is
\begin{equation} \label{E:QuadraticListSeries}
x^2 - x = \frac{\pi ^2}{3} + \sum_{n=1}^{\infty} \left( \frac{4 (-1)^n}{n^2} \cdot \cos(nx) + \frac{2 (-1)^n}{n} \cdot \sin(nx) \right) \, .
\end{equation}

We can also use \verb+fSeries[ ]+ to specify $k$ functions of $x$ over $k$ subintervals,
which are specified using $k+1$ endpoints.

Example: Here is Equation \eqref{E:gxTentative2} again, but extended to ($-\pi$, $\pi$), to make it be an odd function:
\begin{equation} \label{E:ThreeSegsAndOdd}
g(x) =
  \begin{cases}
             - (\pi + x)/2   &\text{for $-\pi \leq x < -1$,}\\
    \phantom{-}x(\pi - 1)/2  &\text{for $-1 \leq x < 1$,}\\
    \phantom{-}(\pi - x)/2   &\text{for $\phantom{-}1 \leq x < \pi \, .$}
  \end{cases}
\end{equation}

$g(x)$ is composed of \emph{three} functions.
There are either endpoints or ``corners'' at \emph{four} points: $x = -\pi$, $-1$, $1$, and $\pi$.
The following call to \verb+fSeries+ will compute the Fourier coefficients for $g(x)$.
\begin{verbatim}
  fSeries[x, { -(Pi + x)/2, x(Pi - 1)/2, (Pi - x)/2 }, {-Pi, -1, 1, Pi}, n, 3,
    AspectRatio -> Automatic]
\end{verbatim}

\begin{figure}[ht]
  \mbox{\includegraphics[width=\picSingleWidth]{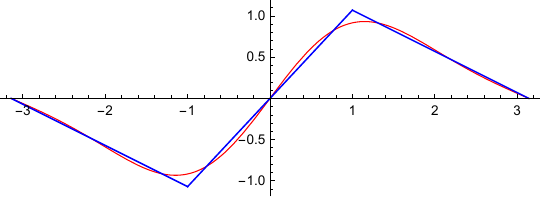}}
  \caption{$g(x)$ in \eqref{E:ThreeSegsAndOdd} with 3 terms of its Fourier series}
  \label{fig:ThreeSegsAndOdd}
\end{figure}

The result is
\[
  \left\{0, \quad 0, \quad \frac{\sin (n)}{n^2}\right\} \, ,
\]
so the Fourier series is
\[
\sum_{n = 1}^{\infty} \frac{\sin (n)}{n^2} \cdot \sin(n x) \, ,
\]
as we have seen before.
The graph is in Figure \ref{fig:ThreeSegsAndOdd}.

You will notice that this call to \verb+fSeries[ ]+ has an extra parameter at the end.
This optional parameter specifies a \verb+PlotStyle+ directive to \textit{Mathematica}'s \verb+Plot[ ]+ function.

In this case, \verb+AspectRatio -> Automatic+ tells \textit{Mathematica}
to make one unit on the $x$ axis be the same size as one unit on the $y$ axis.
You can also include several \textit{Mathematica} \verb+Plot[ ]+ options, such as \verb+ AspectRatio -> Automatic, PlotStyle -> Green+.


Another example: The square wave in Figure \ref{fig:nonSymmSquareWave} consists of two pieces between 0 and $\pi$:
\begin{equation*}
f(x) =
  \begin{cases}
    -1            &\text{for $0 < x < \pi/2$,}\\
    \phantom{-}2  &\text{for $\pi/2 < x < \pi \, .$}
  \end{cases}
\end{equation*}
To compute the Fourier series for $f(x)$, we enter the two functions  of $x$, namely, the (constant) functions $-1$ and 2,
along with the $x$ values of the three endpoints (0, $\pi/2$, and $\pi$):
\begin{verbatim}
  fSeries[x, {-1, 2}, {0, Pi/2, Pi}, n, 10]
\end{verbatim}

\begin{figure}[ht] 
\centering
\begin{minipage}[t]{\figureMinipageWidth}
\centering
\includegraphics[width=\picDblWidth]{{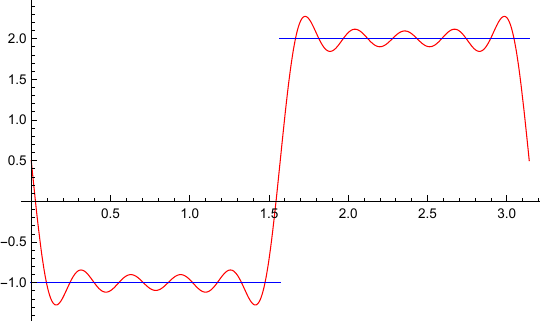}}
\caption{Square wave with 10 terms of its Fourier series}
\label{fig:nonSymmSquareWave}
\end{minipage}\hfill
\begin{minipage}[t]{\figureMinipageWidth}
\centering
\includegraphics[width=\picDblWidth]{{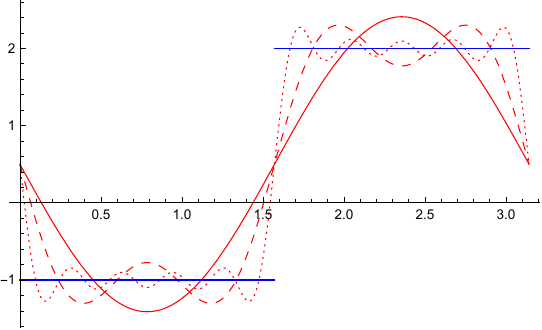}}
\caption{Square wave with sums of 2, 4, and 10 terms of its Fourier series}
\label{fig:squareWave3Sums}
\end{minipage}
\end{figure}

This returns
\[
\left\{\frac{1}{2}, \quad 0, \quad \frac{3 \left((-1)^n-1\right)}{\pi  n}\right\} \, .
\]
Notice that the interval we specified in the call to \verb+fSeries+ is $0 < x < \pi$.
\verb+fSeries+ always takes the period $p$ to be the length of the entire interval: here, this is $p = \pi - 0 = \pi$.
Therefore, because of Equation \eqref{E:FxRepPeriodP}, we need to include a factor of $2 \pi / p = 2$ in the arguments of the $\sin$ and $cos$ functions in the Fourier series.
So, the Fourier series for $f(x)$ is
\[
\frac{1}{2} + \frac{3}{\pi}\sum_{n = 1}^{\infty} \frac{ \left( (-1)^n - 1 \right) }{n} \cdot \sin(2 n x)
= \frac{1}{2} - \frac{6}{\pi} \left( \sin(2 x) + \frac{\sin(6 x)}{3} + \frac{\sin(10 x)}{5} + \ldots   \right)
 \, .
\]

\textbf{An easier way to get the Fourier series for a piecewise linear function.}

There is an easier way to get the Fourier series of a function that consists of linear or constant pieces.
Suppose you didn't know the \emph{functions} that made up $g(x)$, but you did know
the $(x, y)$ \emph{values} at the endpoints and corners.
\verb+fSeries[ ]+ can use those values to compute the functions for you, and then compute the Fourier series.

This version of \verb+fSeries+ makes it easy to experiment with new series by entering interesting-looking combinations of $(x, y)$ pairs.

Even if the function has a jump discontinuity, we can specify the functions using only the $(x, y)$ pairs.
For example, at $x = \pi/2$, the square wave in Figure \ref{fig:nonSymmSquareWave} jumps from $y = -1$ to $y = 2$.
We can treat this function as if it has \emph{two} distinct $y$ values at $x =  \pi/2$.
So, we simply enter \emph{both} $(x, y)$ pairs, $(\pi/2, -1)$ and $(\pi/2, 2)$:
\begin{verbatim}
  fSeries[x, { {0, -1}, {Pi/2, -1}, {Pi/2, 2}, {Pi, 2} }, n, 10]
\end{verbatim}

This gives the same result as in Figure \ref{fig:nonSymmSquareWave}.

This code
\begin{verbatim}
  fSeries[x, { {0, -1}, {Pi/2, -1}, {Pi/2, 2}, {Pi, 2} }, n, {2, 4, 10}]
\end{verbatim}
plots the square wave in Figure \ref{fig:nonSymmSquareWave}, plus \emph{three}
Fourier series approximations to the square wave: namely the sums of 2, 4, and 10 terms of the Fourier series.
This graph is in Figure \ref{fig:squareWave3Sums}.


For another example, suppose we know that our function consists of \emph{linear} pieces between these $(x, y)$ pairs:
$(-\pi, 0)$, $(-1, -\frac{\pi - 1}{2})$, $(1, \frac{\pi - 1}{2})$, and $(\pi, 0)$.
(Note: this is $g(x)$ in Equation \eqref{E:ThreeSegsAndOdd}.)

The \verb+fSeries+ function lets you enter just the list of $(x, y)$ pairs at the endpoints and corners
without having to first compute the linear functions yourself.
So,
\begin{verbatim}
  fSeries[x, { {-Pi, 0}, {-1, -(Pi - 1)/2}, {1, (Pi - 1)/2}, {Pi, 0} }, n, 3]
\end{verbatim}
gives the same graph as in Figure \ref{fig:ThreeSegsAndOdd}).

\textbf{Two ways to count the terms of a Fourier series.}

\begin{figure}[ht] 
\centering
\begin{minipage}[t]{\figureMinipageWidth}
\centering
\includegraphics[width=\picDblWidth]{{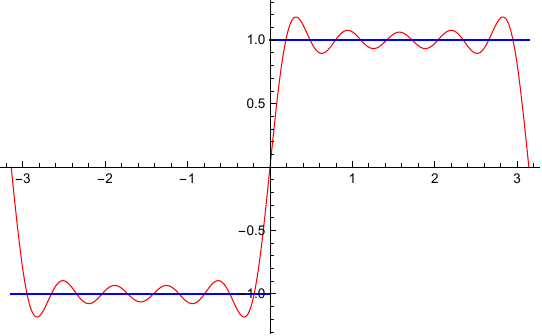}}
\caption{Square wave with the sum of the first 5 terms of the Fourier series \eqref{E:unitSquareWaveCopy}}
\label{fig:unitSquareWave5Peaks}
\end{minipage}\hfill
\begin{minipage}[t]{\figureMinipageWidth}
\centering
\includegraphics[width=\picDblWidth]{{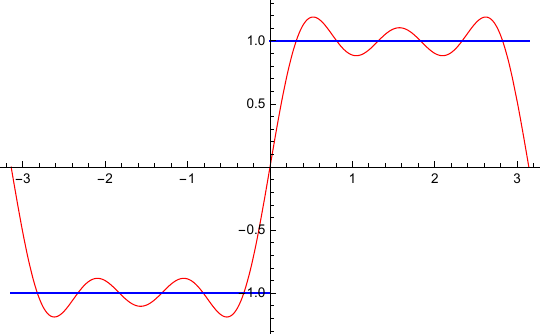}}
\caption{fSeries plot of the sum of the first 5 terms of the Fourier series \eqref{E:unitSquareWave135}}
\label{fig:unitSquareWave3Peaks}
\end{minipage}
\end{figure}

Consider a simple square wave, whose series, \eqref{E:unitSquareWave}, we reproduce here for convenience:
\begin{equation}
\frac{4}{\pi} \left( \sin(x) +  \frac{\sin(3x)}{3} + \frac{\sin(5x)}{5} \ldots \right) =
\frac{4}{\pi} \sum_{n=1}^{\infty} \frac{\sin((2n-1)x)}{2n-1} =
  \begin{cases}
   -1  &\text{if $-\pi < x < 0$,} \\
    0  &\text{if $x = 0$,} \label{E:unitSquareWaveCopy} \\
    1  &\text{if $0 < x < \pi$.}
  \end{cases}
\end{equation}
If we plot the sum of the first 5 terms of this series, we get the graph in Figure \ref{fig:unitSquareWave5Peaks}.

This call to \verb+fSeries[x, {-1, 1}, {-Pi, 0, Pi}, n, 5]+ displays the graph in Figure \ref{fig:unitSquareWave3Peaks}.

We told \verb+fSeries[ ]+ to graph the sum of the first 5 terms of the Fourier series.
So, why aren't these graphs the same?
This call to \verb+fSeries[ ]+ returns
\[
\left\{0,0,-\frac{2 \left((-1)^n - 1\right)}{\pi  n}\right\}
= \left\{ 0, 0,\frac{2}{\pi} \cdot \frac{(1 - (-1)^n)}{n}\right\} \, ,
\]
so \verb+fSeries[ ]+ tells us that the Fourier series for the square wave is
\begin{equation} \label{E:unitSquareWave135}
\frac{2}{\pi} \sum_{n = 1}^{\infty} \frac{(1 - (-1)^n)}{n} \sin(n x) \, .
\end{equation}
When $n$ is odd, $1 - (-1)^n = 2/n$.
However, when $n$ is even, $1 - (-1)^n = 0$.

In \eqref{E:unitSquareWaveCopy}, the first 5 terms with $n = 1$ through 5 are $\sin(x)$, $\sin(3x)/3$, $\sin(5x)/5$, $\sin(7x)/7$, and $\sin(9x)/9$.
But in \eqref{E:unitSquareWave135}, the first 5 terms with $n = 1$ through 5 are $\sin(x)$, 0, $\sin(3x)/3$, 0, and $\sin(5x)/5$.
So, although \verb+fSeries[ ]+ does add 5 terms, only 3 of them contribute
anything to the sum the way we wrote it in Equation \eqref{E:unitSquareWaveCopy}.

In the above example, the even-numbered terms are 0.
But in this example,
\begin{equation*}
\frac{2}{\pi} \sum_{n = 1}^{\infty} \frac{(1 + (-1)^n)}{n} \sin(n x) =
\frac{2}{\pi} \sum_{n=1}^{\infty} \frac{\sin(2 n x)}{n} =
  \begin{cases}
            -1 - (2 x/\pi)  &\text{if $-\pi < x < 0$,} \\
  \phantom{-}1 - (2 x/\pi)  &\text{if $0 < x < \pi$} \, ,
  \end{cases}
\end{equation*}
the \emph{odd}-numbered terms are 0.
Here is a function where every \emph{third} term is 0:
\begin{equation*}
f(x) = 
  \begin{cases}
      0  &\text{if $0 < x < \pi/3$,} \\
  \pi/3  &\text{if $\pi/3 < x < \pi$} \, .
  \end{cases}
\end{equation*}
If we compute the sine series of $f(x)$, we find that the $n^{\text{th}}$ sine coefficient is
\[
\frac{2 \left(\cos \left(\frac{\pi  n}{3}\right)-(-1)^n\right)}{3 n} =
\left\{1,-\frac{1}{2},0,-\frac{1}{4},\frac{1}{5},0,\frac{1}{7},-\frac{1}{8},0,-\frac{1}{10},\frac{1}{11},0\right\}
\]
for $n = 1, 2, 3, \ldots, 12$.

Finally, in $g(x)$ in Equations \eqref{E:g3x} and \eqref{E:gx3Guess}, every term is 0 \emph{except} for terms $n = 3, 6, 9, \ldots$.

\textbf{Parseval's equation.}

This package also has a function to evaluate Parseval's equation, Equation \eqref{E:GenParsevalEqn}.

Given the function
\[
f(x) =
  \begin{cases}
             - (\pi + x)/2   &\text{for $-\pi \leq x \leq 0$ ,}\\
    \phantom{-}(\pi - x)/2   &\text{for $\phantom{-}0 < x < \pi \, ,$}
  \end{cases}
\]
we can compute its Fourier series with
\begin{verbatim}
  fSeries[x, {-(Pi + x)/2, (Pi - x)/2}, {-Pi, 0, Pi}, n, 5]
\end{verbatim}
The Fourier series over $-\pi \leq x \leq \pi$ is
\[
f(x) = \sum_{n = 1}^{\infty} \frac{1}{n} \cdot \sin(n x) \, .
\]
The $n$\textsuperscript{th} coefficient is $1/n$.
By integrating the square of $f(x)$, Parseval's equation gives the value of the sum
\[
\sum_{n = 1}^{\infty} \frac{1}{n^2} = \frac{\pi^2}{6} \, .
\]
This call to \verb+parseval[ ]+ does, indeed, return $\pi^2/6$:
\begin{verbatim}
  parseval[x, {-(Pi + x)/2, (Pi - x)/2}, {-Pi, 0, Pi}]
\end{verbatim}

We saw in Equations \eqref{E:Zeta2} and \eqref{E:FS-SpecialQuadratic} that, for $0 \leq x \leq 2 \pi$,
\[
\frac{3 x^2 - 6 \pi x + 2 \pi^2}{12} = \sum_{n=1}^{\infty} \frac{\cos(n x)}{n^2} \, .
\]
If we apply Parseval's equation to this polynomial,
\begin{verbatim}
  parseval[x, (3 x^2 - 6 Pi x + 2 Pi^2)/12, 0, 2 Pi]
\end{verbatim}
we find that
\[
\sum_{n=1}^{\infty} \frac{1}{n^4} = \frac{\pi^4}{90} \, .
\]

As another example, the function in Equation \eqref{E:ThreeSegsAndOdd} has the Fourier series
\[
\sum_{n = 1}^{\infty} \frac{\sin(n)}{n^2} \sin(nx)
\]
which can be computed with either the functions, or the coordinates of the corners and endpoints:
\begin{verbatim}
  fSeries[x, {-(Pi + x)/2, x (Pi - 1)/2, (Pi - x)/2}, {-Pi, -1, 1, Pi}, n, 3]
  fSeries[x, { {-Pi, 0}, {-1, -(Pi - 1)/2}, {1, (Pi - 1)/2}, {Pi, 0} }, n, 3]
\end{verbatim}

Either of these calls to \verb+parseval[ ]+
\begin{verbatim}
  parseval[x, {-(Pi + x)/2, x (Pi - 1)/2, (Pi - x)/2}, {-Pi, -1, 1, Pi}]
  parseval[x, { {-Pi, 0}, {-1, -(Pi - 1)/2}, {1, (Pi - 1)/2}, {Pi, 0} }]
\end{verbatim}
will compute the sum in Equation \eqref{E:ThreePartParseval}:
\[
\sum_{n = 1}^{\infty} \left( \frac{\sin(n)}{n^2} \right)^2 =
\sum_{n = 1}^{\infty} \frac{\sin^2(n)}{n^4} = \frac{(\pi - 1)^2}{6} \, .
\]

\textbf{\textit{Mathematica}'s built-in functions.}

Most of this package was written before \textit{Mathematica} introduced built-in functions such as\\
\verb+FourierSeries+ and \verb+FourierCoefficient+,
which express coefficients using complex numbers.


For example, the list in \eqref{E:QuadraticList} gives the Fourier series coefficients of $x^2 - x$ over $(-\pi, \pi)$.
Equation \eqref{E:QuadraticListSeries} shows the series written out.

In Version 13 of \textit{Mathematica}, \verb+FourierCoefficient[x^2 - x, x, n]+ gives the coefficients as:
\begin{equation*}
\begin{dcases}  
 \frac{\pi ^2}{3} & \text{ if } n = 0 \, , \\
 \frac{(-1)^n (2-i n)}{n^2} \, .
\end{dcases}
\end{equation*}
With $n = 0$, 1, 2, and 3, the first few coefficients are
\[
 \frac{\pi ^2}{3}, \text{ } i-2, \text{ } \frac{1}{4} (2-2 i), \text{ and } \frac{1}{9} (3 i-2) \, .
\]
\verb+ser3 = FourierSeries[x^2 - x, x, 3]+ shows the first 3 terms of the Fourier series in the form
\[
\frac{\pi ^2}{3} +
(-2-i) e^{-i x}-(2-i) e^{i x}+\left(\frac{1}{2}+\frac{i}{2}\right) e^{-2 i
   x}+\left(\frac{1}{2}-\frac{i}{2}\right) e^{2 i x}-\left(\frac{2}{9}+\frac{i}{3}\right) e^{-3 i
   x}-\left(\frac{2}{9}-\frac{i}{3}\right) e^{3 i x}  \, .
\]
We can use \verb+ExpToTrig[ser3]+ to convert this to a series with sines and cosines.
We get
\[
\frac{\pi^2}{3}
-2 \sin (x)+\sin (2 x)-\frac{2}{3} \sin (3 x)-4 \cos (x)+\cos (2 x)-\frac{4}{9} \cos (3 x) \, .
\]
These agree with the first few terms in Equation \eqref{E:QuadraticListSeries}.

We can also use the built-in function \verb+FourierTrigSeries[x^2 - x, x, 3]+ to get the sine and cosine coefficients directly.
This gives
\[
\frac{\pi ^2}{3} + 
4 \left(-\cos (x)+\frac{1}{4} \cos (2 x)-\frac{1}{9} \cos (3 x)\right) +
2 \left(-\sin (x)+\frac{1}{2} \sin (2 x)-\frac{1}{3} \sin (3 x)\right) \, .
\]

\newpage

\medskip


\noindent\textnormal{rjbaillie2.718@gmail.com; Middletown, NY}

\newpage

\pagestyle{empty} 











\ifthenelse {\boolean{BKMRK}}
  { \section{Appendix: A collection of interesting formulas involving \texorpdfstring{$\pi$}{pi}}\label{S:CollectionOfFormulas} }
  { \section{Appendix: A collection of interesting formulas involving pi}\label{S:CollectionOfFormulas} }
For any $0 < k \leq \pi$ and any $x$ in $1 \leq x \leq \frac{2\pi}{k} - 1$,
\[
\sum_{n=1}^{\infty} \frac{\sin(knx)}{kn} = \sum_{n=1}^{\infty} \frac{\sin(knx)}{kn} \cdot \frac{\sin(kn)}{kn} = \frac{\pi - kx}{2k}
 \quad \quad \text{(Equation \ref{E:twoGeneralSumsWithSinkx})}
\]
%
%
Interesting special cases are $x = 1$, with $k = 1$, 2, and 3.
We get four series with $\sum a_n = \sum (a_n)^2$:
\[
 \sum_{n=1}^{\infty} \frac{\sin(n)}{n} = \sum_{n=1}^{\infty} \left( \frac{\sin(n)}{n} \right) ^{2} = \frac{\pi - 1}{2}
  \quad \quad \text{(Equation \ref{E:Prob6241a})}
\]
\[
 \sum_{n=1}^{\infty} \frac{\sin(2n)}{2n} = \sum_{n=1}^{\infty} \left( \frac{\sin(2n)}{2n} \right) ^{2} = \frac{\pi - 2}{4}
 \quad \quad \text{(Equation \ref{E:SincEvenN})}
\]
\[
 \sum_{n=1}^{\infty} \frac{\sin(3n)}{3n} = \sum_{n=1}^{\infty} \left( \frac{\sin(3n)}{3n} \right) ^{2} = \frac{\pi - 3}{6}
 \quad \quad \text{(Equation \ref{E:sin3n})}
\]
Because of Equations \ref{E:Prob6241a} and \ref{E:SincEvenN}, the sum over \emph{odd} $n$ also has the same property:
\[
 \sum_{n=1}^{\infty} \frac{\sin(2n-1)}{2n-1} = \sum_{n=1}^{\infty} \left(\frac{\sin(2n-1)}{2n-1} \right) ^{2} = \frac{\pi }{4}
 \quad \quad \text{(Equation \ref{E:SincOddN})}
\]
Parseval's equation applied to Equation \eqref{E:gkDef2}: For $0 < k \leq \pi$, we have
\[
 \sum_{n=1}^{\infty} \left( \frac{\sin(kn)}{(kn)^2} \right)^2 = \frac{(\pi-k)^2}{6k^2}
 \quad \eqref{E:gkParseval} \text{; so with $k = 1$,} \quad
 \sum_{n=1}^{\infty} \frac{\sin^{2}(n)}{ n^{4} } = \frac{(\pi - 1)^{2} }{6} \quad \eqref{E:Prob6241b}
\]
A variation on the classic Gregory-Leibniz series for $\pi $/4:
\[
 \frac{\pi }{4} = 1- \frac{1}{3} + \frac{1}{5} - \frac{1}{7} + \dots =
 1\cdot \frac{\sin(1)}{1} - \frac{1}{3} \cdot \frac{\sin(3)}{3} + \frac{1}{5} \cdot \frac{\sin(5)}{5} - \frac{1}{7} \cdot \frac{\sin(7)}{7} + \dots
 \quad \eqref{E:GregoryGen1}
\]
A generalization of Equation \eqref{E:GregoryGen1}: for all non-zero $x$ in $-\pi /2 \leq x \leq \pi /2$,
\[
 \frac{\pi }{4} = 1-\frac{1}{3} + \frac{1}{5} - \frac{1}{7} + \cdots =
 1\cdot \frac{\sin(x)}{x} - \frac{1}{3} \cdot \frac{\sin(3x)}{3x} + \frac{1}{5} \cdot \frac{\sin(5x)}{5x} -\frac{1}{7} \cdot \frac{\sin(7x)}{7x} + \dots
 \quad \eqref{E:GregoryGen2}
\]
For $0 \leq x \leq \pi - 2$ (including $x = 1$), these \emph{three} sums equal $x/2$ (Thm. \ref{T:AltSumkxTheorem}; Eqs. \ref{E:xOver2ab}, \ref{E:AltSum3xResult1}, \ref{E:ThreeAlternatingSincSums}):
\[
  \sum_{n=1}^{\infty} (-1)^{n+1} \frac{\sin(nx)}{n} =
  \sum_{n=1}^{\infty} (-1)^{n+1} \frac{\sin(nx)}{n} \frac{\sin(n)}{n} =
  \sum_{n=1}^{\infty} (-1)^{n+1} \frac{\sin(nx)}{n} \left( \frac{\sin(n)}{n} \right)^2 =
  \frac{x}{2}
\]
\[
 \text{A series where } \sum a_n = \sum (a_n)^3 \text{ : } \quad
 \sum_{n=1}^{\infty} (-1)^{n+1} \frac{\sin(n)}{n} = 
 \sum_{n=1}^{\infty} \left( (-1)^{n+1} \frac{\sin(n)}{n} \right ) ^3 =
 \frac{1}{2} \quad \eqref{E:SumAndSumCubed} \qquad \qquad
\]
All \emph{four} of these sums equal $(\pi - 3)/2$ (Equation \ref{E:FourSumsWithN}):
\[
 \sum_{n=1}^{\infty} \frac{\sin(3n)}{n}  =
 \sum_{n=1}^{\infty} \frac{\sin(n)}{n} \frac{\sin(3n)}{n} =
 \sum_{n=1}^{\infty} \left(\frac{\sin(n)}{n} \right) ^{2} \frac{\sin(3n)}{n} =
 \sum_{n=1}^{\infty} \left(\frac{\sin(n)}{n} \right) ^{3} \frac{\sin(3n)}{n} =
 \frac{\pi - 3}{2}
\]
\[  
x = 1 \text{ in Equation \eqref{E:thmHP1} gives }
 \sum_{n=1}^{\infty} \frac{\sin ^{3}(n)}{n} =
 \sum_{n=1}^{\infty} \frac{\sin ^{4}(n)}{n^2} = \frac{\pi }{4}
\]
\[
x = 1 \text{ in Equation \eqref{E:thmHP2} gives }
 \sum_{n=1}^{\infty} \frac{\sin ^{5} (n)}{n} =
 \sum_{n=1}^{\infty} \frac{\sin ^{6} (n)}{n^{2} } = \frac{3 \pi }{16}
\]

\end{document}